\documentclass[12pt]{article}
 
\pdfoutput=1
\usepackage[small]{titlesec}
\usepackage[margin=1in]{geometry} 
\usepackage[fleqn]{amsmath}
\usepackage{mathtools}
\usepackage{graphicx}
\usepackage{epstopdf}
\usepackage{caption}
\usepackage{subcaption}
\usepackage{amscd,amssymb,stmaryrd}
\usepackage{amsthm}
\usepackage{enumerate}
\usepackage{color}
\usepackage{tikz,mathpazo}
\usetikzlibrary{shapes.geometric, arrows}
\tikzstyle{startstop} = [rectangle, rounded corners, minimum width=3cm, minimum height=1cm,text centered, draw=black, fill=red!30]
\tikzstyle{io} = [trapezium, trapezium left angle=70, trapezium right angle=110, minimum width=3cm, minimum height=1cm, text centered, draw=black, fill=blue!30]
\tikzstyle{process} = [rectangle, minimum width=3cm, minimum height=1cm, text centered, draw=black, fill=orange!30]
\tikzstyle{decision} = [diamond, minimum width=3cm, minimum height=1cm, text centered, draw=black, fill=green!30]
\tikzstyle{coord}=[coordinate, node distance=6mm and 25mm]
\tikzstyle{arrow} = [thick,->,>=stealth]

\newcommand{\R}{\mathbb{R}}

\begin{document}
 
% --------------------------------------------------------------
%                         Start here
% --------------------------------------------------------------
 
\title{Numerical inversion of 3D geodesic X-ray transform arising from traveltime tomography}
\author{Tak Shing Au Yeung\thanks{Department of Mathematics, The Chinese University of Hong Kong, Hong Kong SAR.} \and
Eric T. Chung\thanks{Department of Mathematics, The Chinese University of Hong Kong, Hong Kong SAR.} \and 
Gunther Uhlmann\thanks{Department of Mathematics and Institute of Advanced Study, The University of Science and Technology of Hong Kong, Hong Kong SAR;
and Department of Mathematics, University of Washington, Seattle, USA.}
} 
 
\maketitle
 
 \begin{abstract}
In this paper, 
we consider the inverse problem of determining an unknown function defined in three space dimensions
from its geodesic X-ray transform. 
The standard X-ray transform is defined on the Euclidean metric and is given by the integration of a function along straight lines.
The geodesic X-ray transform is the generalization of the standard X-ray transform in Riemannian manifolds and is defined by integration of a function along geodesics.
This paper is motivated by Uhlmann and Vasy's theoretical reconstruction algorithm for geodesic X-ray transform and mathematical formulation for traveltime tomography to develop a novel numerical algorithm
for the stated goal. Our numerical scheme is based on a Neumann series approximation and a layer stripping approach.
In particular, we will first reconstruct the unknown function by using a convergent Neumann series
for each small neighborhood near the boundary. Once the solution is constructed on a layer near the boundary, 
we repeat the same procedure for the next layer, and continue this process until the unknown function is recovered on the whole domain.
One main advantage of our approach is that the reconstruction is localized, and is therefore very efficient,
compared with other global approaches for which the reconstructions are performed on the whole domain. 
We illustrate the performance of our method by showing some test cases including the Marmousi model. 
Finally, we apply this method to a travel time tomography in 3D, in which the inversion of the
geodesic X-ray transform is one important step, and present several numerical results to validate the scheme. 

\end{abstract}

\section{Introduction}
We address in this paper the inverse problem of determining an unknown function from its geodesic X-ray transform. The travel time information and geodesic X-ray data are encoded in the boundary distance function, which measures the distance between boundary points. The geodesic X-ray transform is associated with the boundary rigidity problem of determining a Riemannian metric on a compact manifold from its boundary distance function (see \cite{Ivanov,Stefanov1}). In this paper,  we consider the linearization of the boundary rigidity problem in compact Riemannian manifold with boundary in the isotropic case. Recent progress on the inverse problem of inverting geodesic X-ray transform \cite{Frigyik,Pestov1,Uhlmann1}, and the boundary rigidity problem \cite{Kurley,Pestov2,Stefanov2,Stefanov3,Stefanov4} has motivated us to transfer these theoretical results into numerical methods for determining an unknown function from its travel time information and geodesic X-ray data. This paper is mostly motivated by Uhlmann and Vasy's theoretical reconstruction algorithm for geodesic X-ray transform \cite{Uhlmann1}  and mathematical formulation for traveltime tomography in \cite{Chung1,Chung2} to develop a novel numerical reconstruction algorithm.

Traveltime tomography handles the problem of determining the internal properties, such as index of refraction, of a medium by measuring the traveltimes of waves going through the medium. It arose in global seismology in determining the earth velocity models by measuring the traveltimes of seismic waves produced by earthquakes measured at different seismic stations. Compared to other methods such as full wave inversion, traveltime tomography is generally easier to implement and can be computed much efficiently.

The standard X-ray transform is defined on the Euclidean metric and is given by the integration of a function along straight lines. It is the mathematical model underlying medical imaging techniques such as CT and PET. The geodesic X-ray transform is the generalization of standard X-ray transform in Riemannian manifolds and is defined by integrating a function along geodesics. Since the speed of elastic waves increases with depth in the Earth, the rays curve back to the surface and thus  geodesic X-ray transform arises in geophysical imaging in determining the inner structure of the Earth .

In two dimensions, the uniqueness and stability for the geodesic X-ray transform was proved by Mukhometov \cite{Mukhometov1} on simple surface and also for a general family of curves. 
For simple manifolds, the case of geodesics was generalized to higher dimensions in \cite{Mukhometov2}.
When the metric is not simple, the results are proved for some geometries with some additional assumptions. In \cite{Frigyik}, the authors proved some generic injectivity results and stability estimates for the X-ray transform for a general family of curves and weights with additional assumptions (analog to simplicity in the geodesic X-ray transform)
of microlocal condition which includes real-analytic metrics for a class of non-simple manifolds in dimension $d \geq 3$.

Besides some injectivity and stability results, some explicit methods for reconstructing the unknown function have been found. Several methods have been introduced to recover the unknown function from X-ray transforms, including Fourier methods, backprojection and singular value decomposition \cite{Natterer}. In \cite{Guo}, the authors introduce a reconstruction algorithm to recover functions from their 3D X-ray transforms by shearlet and wavelet decompositions. This algorithm is effective and optimal in mean square error rate. On the other hand, under the afoliation condition, Uhlmann and Vasy showed that the geodesic X-ray transform is globally injective and can be inverted locally in a stable manner \cite{Uhlmann1}.  Also, the authors propose a layer stripping type algorithm for recovering the unknown function in the form of a convergent Neumann series. On the numerical side, the geodesic X-ray transform has been implemented numerically in two dimensions. In \cite{Monard}, Monard develops a numerical implementation of geodesic X-ray transform based on the fan-beam geometry and the Pestov-Uhlmann reconstruction formulas \cite{Pestov1}. 
%This implementation can be used in most general case and both simple and non-simple metrics. \\

In this paper we consider the inverse problem of computing the unknown function $f$ defined in a domain $\Omega$ from its geodesic X-ray transform in three dimensions. 
As far as we know, this is the first numerical implementation in 3D. 
Our numerical method is based on the theoretical results in Uhlmann and Vasy \cite{Uhlmann1}, which shows that
the unknown function can be represented by a convergent Neumann series and can be recovered by a layer stripping type algorithm. 
Based on this theoretical foundation, we develop a numerical approach using a truncated Neumann series
and layer stripping. A related method is developed for photoacoustic tomography in the work by Qian, Stefanov, Uhlmann, and Zhao \cite{Qian}.
%to the Neumann series based numerical method 
%which is developed for photoacoustic tomography in a paper by Qian, Stefanov, Uhlmann, and Zhao \cite{Qian} and also the idea of layer stripping algorithm by Uhlmann and Vasy \cite{Uhlmann1}. 
%Our algorithm is based on a convergent Neumann series, giving an efficient and convergent numerical scheme that recovers the initial condition of an acoustic wave equation with non-constant sound speeds by boundary measurements. 
The key ideas of our approach are discussed as follows. 
First, we use
a layer stripping algorithm to separate the domain into small regions and then reconstruct the unknown function on each small region by using a truncated Neumann series and the given geodesic X-ray data restricted to the region. We start this reconstruction process with the outermost region. Combining the above methodologies, the resulting algorithm can undergo the reconstruction layer by layer so to prevent full inversion in the whole domain. Thus, our approach is very efficient and consumes much less computational time and memory. Numerical examples including synthetic isotropic metrics and the Marmousi model are shown to validate our new numerical method.  

Our proposed numerical approach for geodesic X-ray transform is a foundation of a new numerical scheme for traveltime tomography in 3D. 
Our new method is based on the  Stefanov-€"Uhlmann identity which links two Riemannian metrics with their travel time information. The method first uses those geodesics that produce smaller mismatch with the travel time measurements and continues on in the spirit of layer stripping. We then apply the above reconstruction method to the Stefanov-€"Uhlmann identity in order to solve the inverse problem for reconstructing the index of refraction of a medium. %We demonstrate that our method can recover the isotropic metrics and the three dimensional Marmousi synthetic model. 
The new method is motivated by the methods in \cite{Chung1,Chung2}, but we use the above numerical X-ray transform for solving an integral equation involved
in the inversion process. We will present some numerical test cases including the Marmousi model to show the performance of the inversion. 

The paper is organized as follows. In Section 2, we will present some background materials and the reconstruction method, and in Section 3, we present in detail the numerical algorithm and implementation. Numerical results are shown in Section 4 to demonstrate the performance of our method. 
In Section \ref{sec:travel}, we present our inversion method for traveltime tomography and numerical results. 
The paper ends with a conclusion. 
%%%%%%%%%%%%%%%%%%%%%%%%%%%%%%%%%%%%%%%%%%%%%%
%2 Mathematical formulation for traveltime tomography                     				                %
%%%%%%%%%%%%%%%%%%%%%%%%%%%%%%%%%%%%%%%%%%%%%%
\section{Problem description}

In this section, we will present the mathematical formulation of the inverse problem of reconstructing an unknown function
using its geodesic X-ray transform, that is, integrals of the function along geodesics. 
We will then present the numerical reconstruction algorithm
based on the theoretical foundation in \cite{Uhlmann1}.
The algorithm is based on a representation of the unknown function by a convergent Neumann-series.

\subsection{Geodesics}
We first introduce the notion of geodesics. 
Let $\Omega$ be a  strictly convex bounded domain in $\R^d$ and let $(g_{ij})$ be a Riemannian metric defined on it. In the case of isotropic medium,
$$g_{ij}=\frac{1}{c^2}\delta_{ij},$$ 
where $c$ is a given function defined in $\R^d$. Following \cite{Chung1,Chung2,Uhlmann1}, we define the Hamiltonian $H_g$ by 
$$ H_g(x,\xi)=\frac{1}{2}(\sum^d_{i,j=1}g^{ij}(x)\xi_i\xi_j-1)$$
for each $x \in \Omega$ and $ \xi\in\R^d$, where $(g^{ij}) = (g_{ij})^{-1}$ is the inverse of $(g_{ij})$.
Let $X^{(0)}=(x(0),\xi(0))$ be a given initial condition, where $x(0)\in\partial\Omega$ and $\xi(0)\in\R^d$, such that the following inflow conditions hold:
$$H_g(x(0),\xi(0))=1, \qquad \sum_{i,j=1}^d g^{ij}(x(0))\xi_i\nu_j (x(0))<0$$
where $\nu(x)$ is the unit outward normal vector of $\partial\Omega$ at the point $x$. We define $X_g(s,X^{(0)})=(x(s),\xi(s))$ by the solution to the hamiltonian system defined by 
\begin{flalign}
\label{eq:ODE}
\noindent&\frac{dx}{ds}=\frac{\partial H_g}{\partial \xi} \qquad , \qquad \frac{d\xi}{ds}=-\frac{\partial H_g}{\partial x} \\
\text{with the initial} &\text{ condition} \nonumber &\\
&(x(0),\xi(0))=X^{(0)}. \nonumber
 \end{flalign}
The solution $X_g$ defines a geodesic in the physical space, parametrically via $x(s)$ with the co-tangent vector $\xi(s)$ at any point $x(s)$. The parameter $s$ denotes travel time. Thus, we denote the set of all geodesics $X_g$, which are contained in $\Omega$ with endpoints on $\partial\Omega$, by $\mathcal{M}_{\Omega}$.

\subsection{The reconstruction method}
Let $f$ be a smooth function from $\Omega$ to $\R$. The aim of
this paper is to develop a numerical scheme to 
determine the unknown function $f(x)$ from the geodesic X-ray transform of $f(x)$. Our method is based on a truncation of a convergent Neumann series. First, we define the geodesic X-ray transform. Given a function $f$ defined on $\Omega$, 
its geodesic X-ray transform is the
collection $\{ (If)(X_g) \}$ of integrals of $f$ along geodesics $X_g\in\mathcal{M}_{\Omega}$, where 
$$
(If)(X_g) := \int_{x(s)} f(s) \; ds,
$$
where $X_g(s) = (x(s),\xi(s))$.
We note that this is the measurement data, and we use this data to reconstruct $f(x)$. 
Let $\Lambda$ be the adjoint of the operator $I$. Then Uhlmann and Vasy \cite{Uhlmann1} show that there is an operator $R$ such that
$$R \Lambda (If) =f-Kf,$$
where the error operator $K$ is small in the sense of $||K|| < 1$ for an appropriate norm. 
Hence, the unknown function $f$ can be represented by
the following convergent Neumann series
\begin{flalign}
\label{eq:rep}
 f=\sum^{\infty}_{n=0}K^n R\Lambda (If).   
\end{flalign}
We remark that the operator $R$ is the inverse of the operator $\Lambda \circ I$. 
Moreover, results from \cite{Uhlmann1} suggest that the unknown function $f$ can be reconstructed locally in a layer by layer fashion. 
In particular, one can first reconstruct the unknown function $f$ using (\ref{eq:rep}) in small neighborhoods near the boundary of the domain,
and then repeat the procedure in the next inner layer of the domain, and so on. 
The challenges in the numerical computations of the unknown function $f$ using the above representation (\ref{eq:rep}) lie in the fact that
computing the operators $\Lambda$ and $R$ as well as the implementation of the layer stripping algorithm
needs some careful constructions. 
It is the purpose
of the rest of the paper to develop numerical procedures to compute these operators and implement the layer stripping algorithm. 

\subsection{The numerical procedure}

We will give a general overview of our numerical scheme, and present the implementation details in Section \ref{sec:num}. 
We will assume that the set of geodesic X-ray transforms $\{ (If)(X_g) \}$ is given (as a set of measurement data)
and is a finite set. 
Let $Z$ be a set of grid points, denoted as $\{ z_i\}$, in the domain $\Omega$. We will determine the values of the unknown function $f(x)$ 
at these grid points using the given data set $\{ (If)(X_g) \}$.
Next, we will present our approach to numerically construct the operators $\Lambda$ and $R$.

For a given point $x\in\Omega$, we define $\mathcal{M}_{\Omega}(x)$ as the set of all geodesics passing through the point $x$. 
%Notice that, in computations, the set $\mathcal{M}_{\Omega}(x)$ is a finite set. 
We use the notation $|\mathcal{M}_{\Omega}(x)|$ to represent the number of elements in the set $\mathcal{M}_{\Omega}(x)$.
Then, numerically, we can define the action of the operator $\Lambda \circ I$
as the average of the line integrals (back-projection) by 
\begin{align}
\label{eq:backp}
\Lambda (If):=\frac{\sum\limits_{j=1}^{|\mathcal{M}_{\Omega}(x)|} (If)(X_g^j) }{|\mathcal{M}_{\Omega}(x)|},
\end{align}
where we use the notations $\{ X_g^j\}$, $j=1,2,\cdots, |\mathcal{M}_{\Omega}(x)|$, to denote the elements in the set $\mathcal{M}_{\Omega}(x)$.
Notice that the above formula defines a function with domain $\Omega$ using the given geodesic X-ray transform data $(If)(X_g)$.
For standard X-ray transform, the above operator $\Lambda$ gives a good approximation to the unknown function $f(x)$, as for the standard back-projection operator. 
For geodesic X-ray transform, this operator provides an initial approximation, which 
is the initial term in a convergent Neumann series representation of $f(x)$.

Motivated by \cite{Frigyik,Uhlmann1}, we will use an operator $A$ to model the action of the operator $R$ presented above. In particular, we will construct
an operator $A$ such that 
\begin{align}
\label{eq:rep1}
(A^*A)^{-1}\Lambda \circ I=Id-K,
\end{align}
where K is an error operator with $||K||<1$ for some appropriate norm, and $Id$ is the identity operator. 
From \cite{Frigyik,Uhlmann1}, we know that the operator $A$ is essentially integrals along geodesics,
and the operator $A^*$ is the adjoint operator of $A$. 
Using the above, we can write the following Neumann series
$$f=\sum^{\infty}_{n=0}K^n (A^*A)^{-1}\Lambda (If)$$
for the unknown function $f(x)$. Note that $If$ is the given data. 
We remark that, in computations, 
the inverse of $A^*A$ is not the actual inverse of the operator $\Lambda \circ I$, but gives a kind of approximation to the operator $R$.
In the following, we will explain a numerical approximation to this operator $A$ and its adjoint $A^*$. 
We remark that both the operators $I$ and $A$ correspond to integrations along geodesics,
but we use different notations since they will be defined on two different grids, which will be explained in the next section.
The situation is similar for the operators $\Lambda$ and $A^*$, which correspond to taking averages of integrals on geodesics.

%%%%%%%%%%%%%%%%%%%%%%%%%%%%%%%%%%%%%%%%%%%%%%
%3 Numerical methods                            										       %
%%%%%%%%%%%%%%%%%%%%%%%%%%%%%%%%%%%%%%%%%%%%%%
\section{Numerical implementations}\label{sec:num}

In the previous section, we presented a general overview of our numerical procedure
and its corresponding theoretical motivations. 
In this section, we give details of the implementations of the numerical reconstruction algorithm. 
There are three important ingredients in our numerical algorithm. They are
\begin{enumerate}
\item Given the data $\{ (If)(X_g) \}$, compute $\Lambda (If)$.
\item Compute the action of $\Lambda \circ I$.
\item Compute the action of $K := Id - (A^*A)^{-1} (\Lambda \circ I)$.
\end{enumerate}

\subsection{Implementation details}
\begin{figure}[ht]
    \centering
    \begin{subfigure}[b]{0.6\textwidth}
        \includegraphics[width=\textwidth]{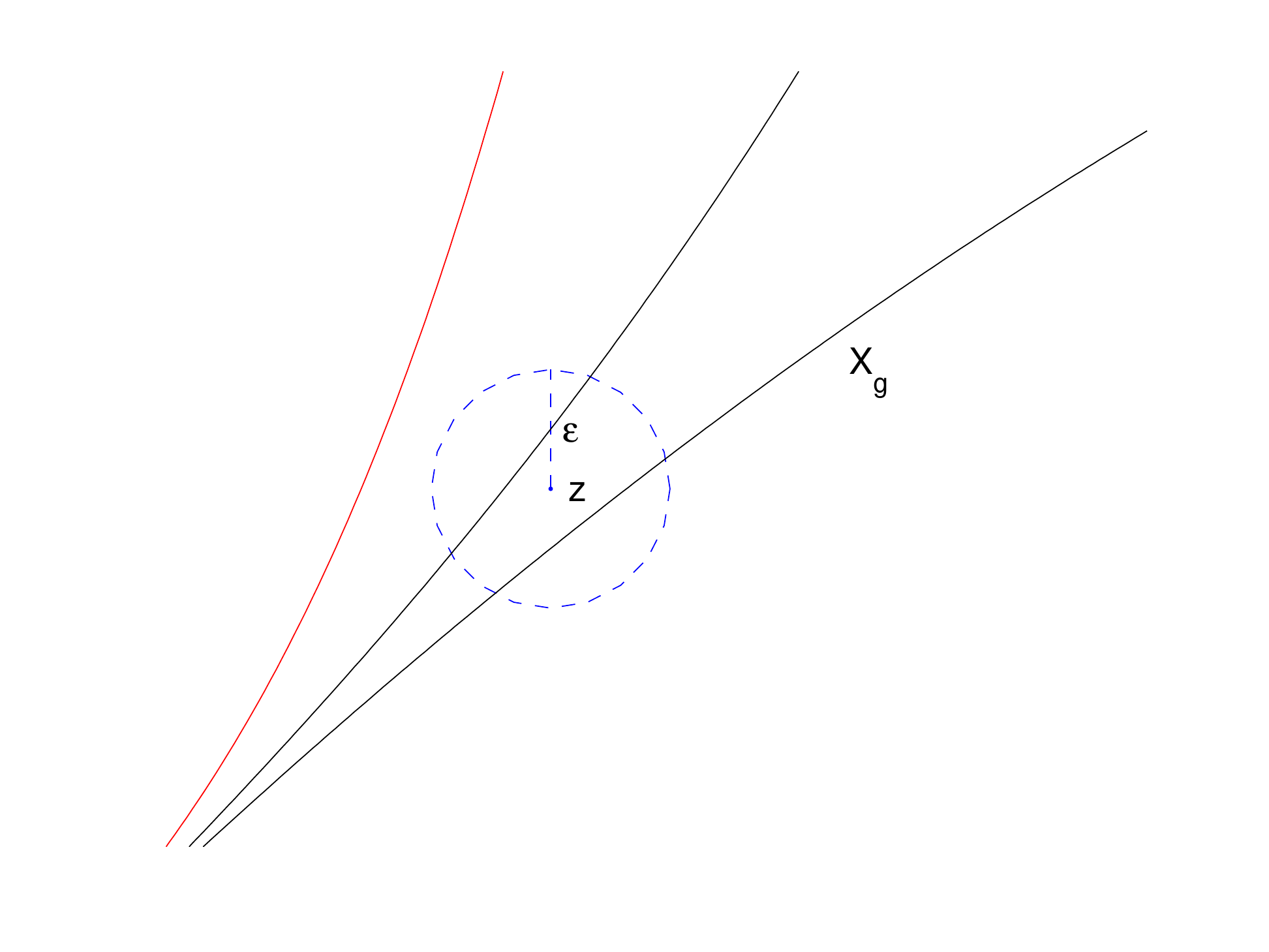}
  
    \end{subfigure}
     
        \caption{An example of choosing the geodesics passing through a given grid point $z$.}\label{fig:geodesic}
\end{figure}

First, we discuss the computation of $\Lambda (If)$ using the given data set $\{ (If)(X_g) \}$.
In fact, we will use the formula (\ref{eq:backp}) for this purpose. We emphasize that, since we will only recover the unknown function $f$
on the set of points $Z$, the output $\Lambda (If)$ is also defined only on the same set of grid points $Z$.
One difficulty is that one, in general, cannot find the geodesic that is passing through a given point $z_i \in Z$ since there are only finitely many geodesics in computations.
We will choose a small neighborhood of $z_i$ and find all the geodesics passing through the $\epsilon$-neighborhood of $z_i$. Then, we can apply the formula (\ref{eq:backp}) using this set of geodesics for each $z_i$. Figure \ref{fig:geodesic} shows how to choose the set of geodesics. The black geodesics shown in Figure \ref{fig:geodesic} are said to be the geodesic that is passing through the $\epsilon$-neighborhood of $z_i$ but the red geodesic is not.

Next, we will discuss the action of the operator $\Lambda \circ I$. Notice that this is needed in the computation
of the error operator $K$. We will discuss the computation of the operator $I$. The action of $\Lambda$
is similar to the above. 
For a given function $f$ whose values are defined only on the set of grid points $Z$,
we will evaluate $If$. 
Thus, we need to compute the integral of $f$ on a geodesic $X_g$. 
To compute a geodesic $X_g$, we
will solve the ray equation (\ref{eq:ODE}) in the phase space starting from a particular initial point $X^{(0)}$ by the 4th order Runge-Kutta Method. 
Then we will obtain a set of points $\{ x(s_i)\}$ defining the geodesic in the physical space. 
Next, we use a version of the trapezoidal rule to compute the line integral $(If)(X_g)$ of the function $f$ along the geodesic $X_g$. In particular, we have 
%$$ (If)(X_g) \approx \sum_i f(x(s_i))|x'(s_i)|(s_{i}-s_{i-1}) \approx \sum_i f(x(s_i))\times \text{dist}( x(s_i)-x(s_{i-1})) ,$$
$$ (If)(X_g) \approx \sum_i f(x(s_i)) \,(x'(s_i)) \, (s_{i}-s_{i-1}).$$
We remark that, in the above formula, $(s_i - s_{i-1})$ is the step size used in the 4th order Runge-Kutta Method. 
In addition, the expression $x'(s_i)$ can be computed using $\displaystyle x'= \frac{\partial H_g}{\partial \xi}$.
The term $f(x(s_i))$ is not well defined since $f$ is only defined on the set of grid points $Z$ and the point $x(s_i)$
may not be one of the grid points. To overcome this issue, we replace $f(x(s_i))$
by $\widehat{f}(x(s_i))$, which is the linear interpolation of $f$ using the grid points near $x(s_i)$. So, we will
use the following formula to compute the action of $f$:
\begin{align}
\label{eq:integral}
(If)(X_g) \approx \sum_i \widehat{f}(x(s_i)) \,(x'(s_i)) \, (s_{i}-s_{i-1}).
\end{align}

%where dist$(x,y)$ is distance function in the given Riemannian metric, which can be approximate by Euclidean norm when $x,y$ are close enough. Then we have the following approximation to $If$:
%$$ (If)(X_g) \approx \sum_i f(x(t_i))\times | x(t_i)-x(t_{i-1})|,$$
%where $|\cdot|$ is the Euclidean norm. 

\begin{figure}[ht]
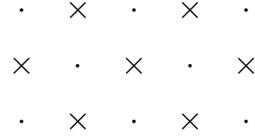

  $$\quad \cdot \quad \times \quad \cdot \quad \times \quad \cdot \quad $$
  $$\quad \times \quad \cdot \quad \times \quad \cdot \quad \times \quad$$
  $$\quad \cdot \quad \times \quad \cdot \quad \times \quad \cdot \quad $$
  \caption{Illustration of the two grids. All the $\{\times\}$ form the grid $Z$ and all $\{\times,\cdot\}$ form the finer grid $Z_f$.}
  \label{fig:grid}
\end{figure}

Finally, we discuss the action of the error operator $K$. 
In this part, the main ingredient is to compute the operator $(A^*A)^{-1}$, which approximates the action of the operator $R$ in (\ref{eq:rep}). 
From \cite{Frigyik,Uhlmann1}, we know that the operator $A$ is essentially integrals along geodesics,
and the operator $A^*$ performs average of line integrals passing through a given point. 
Notice that, the action of $A$ is similar to that of $I$, and the action of $A^*$ is similar to the action of $\Lambda$. 
In order to obtain a good approximation to the operator $R$ and hence a good reconstructed $f$, 
we will perform the action of $A$ and A$^*$
on a finer grid $Z_f$, which is a refinement of the grid $Z$. 
Figure \ref{fig:grid} illustrates the grids in 2D. 
In particular, the grid points in the set $Z$ are represented by $\{ \times\}$,
and the grid points in the set $Z_f$ are represented by both $\{\times, \cdot\}$. 
For the computation of $A$, we use the same formula as in (\ref{eq:integral})
but with the grid $Z_f$.
For the computation of $A^*$, we use the same formula as in (\ref{eq:backp}) but with the grid $Z_f$. 
Thus, the operator $(A^*A)$ and its inverse $(A^*A)^{-1}$ take inputs as functions defined on the grid $Z_f$,
and give outputs also as functions defined on the grid $Z_f$. 

%Since we need $||K||<1$ in the Neumann series, we use another grid to approximate operator $\Lambda$ which is different from $A^*A$. Then, two set of grid points are needed in this method. One is fine rectangular grid $Z_f$ and another is coarse cross-hatch grid $Z_c$. The following is the example of 2D grids:

To complete the steps, we need a projection operator $P$, which maps functions defined on the finer grid $Z_f$
to functions defined on the grid $Z$. In fact, the operator $P$ is the restriction operator. Then the operator $P^*$
maps functions defined on the grid $Z$ to functions defined on the finer grid $Z_f$. Now, we can write 
down the reconstruction formula:
$$f=\sum^{\infty}_{n=0}K^n P(A^*A)^{-1}P^*\Lambda (If)$$
where
$$K = Id - P (A^*A)^{-1} P^*(\Lambda \circ I).$$
To construct $f$, we will use a truncated version of the above infinite sum. Notice that, to regularize the problem,
we will replace the above sum by
\begin{equation}
\label{eq:Neumann}
f=\sum^{\infty}_{n=0}K^n P(A^*A- \delta \Delta)^{-1}P^*\Lambda (If)
\end{equation}
where $\delta >0$ is a regularization parameter and $\Delta$ is the Laplace operator.

 \begin{figure}[ht]
    \centering
    \begin{subfigure}[b]{0.4\textwidth}
        \includegraphics[width=\textwidth]{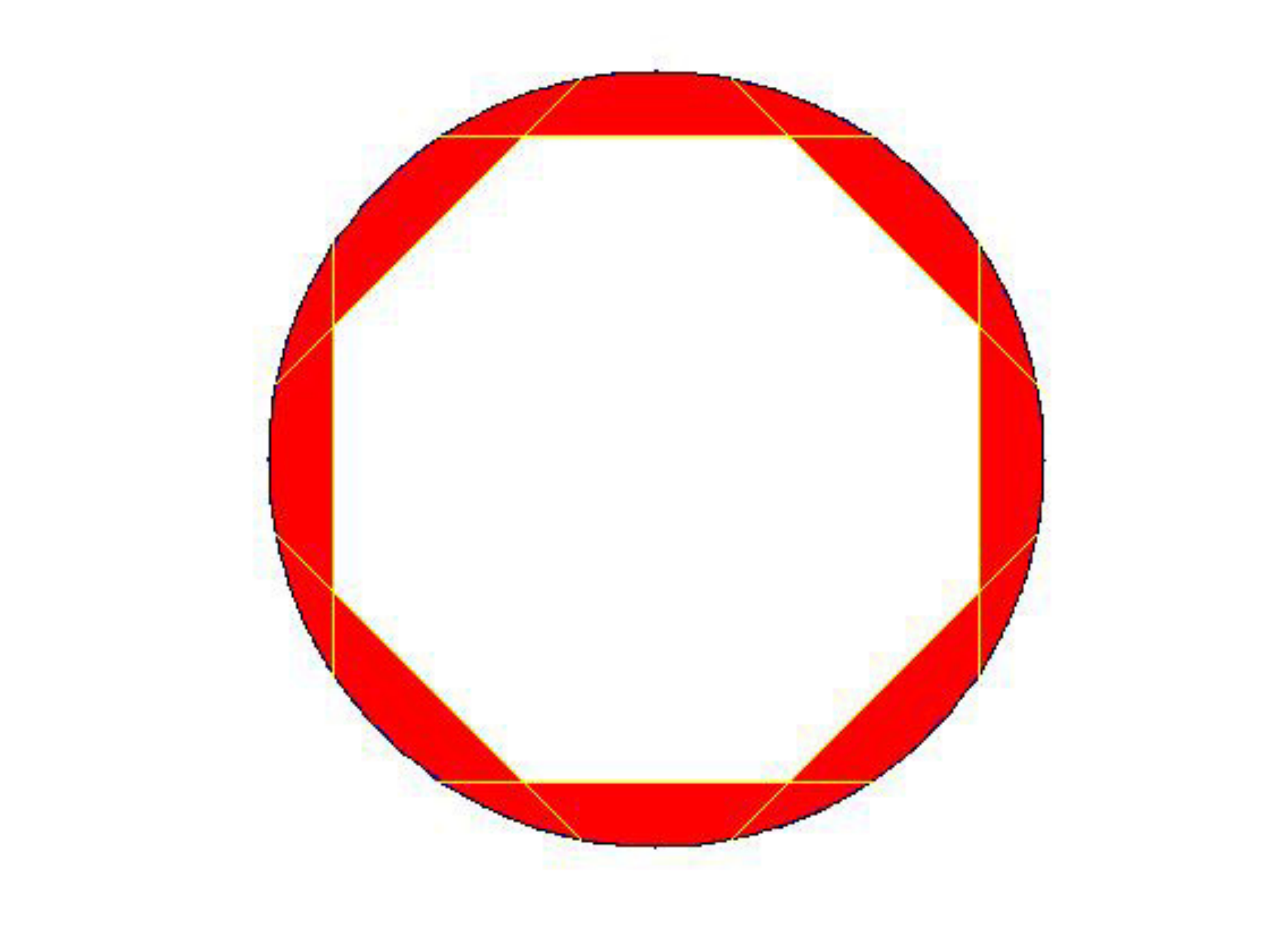}
  
    \end{subfigure}
     \begin{subfigure}[b]{0.4\textwidth}
        \includegraphics[width=\textwidth]{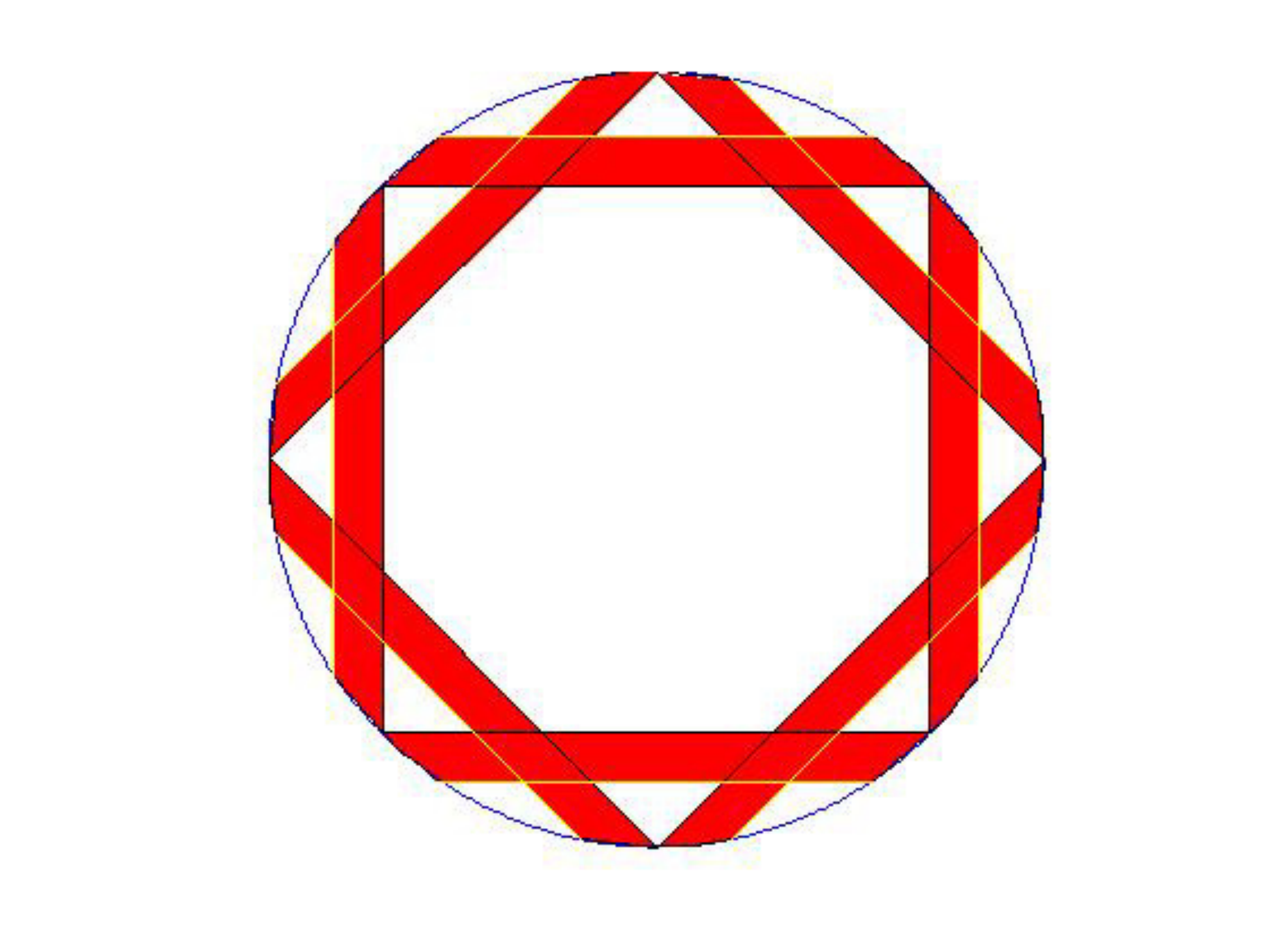}

    \end{subfigure}
        \caption{An example of dividing the first and second layers of 3D domain into disks.}\label{fig: layer}
\end{figure}

\subsection{A layer stripping algorithm for 3D models}

The algorithm can be implement on a 3D domain. We consider domains with compact closure $\overline{\Omega}$ equipped with a function $\rho : \overline{\Omega}\to [0,\infty)$ whose level sets $\Sigma_t=\rho^{-1}(t)$, $t<T$, are strictly convex, with $\Sigma_0=\partial X$ and $X\setminus \cup_{t\in[0,T)}\Sigma_t=\rho^{-1}([T,\infty))$ either having 0 measure of having empty interior. In general, one can take a finite partition of $[0,T]$ by such intervals $[t_i,t_{i+1})$, $i=1,\dots,k$ and proceed inductively to recover $f$ on $ \cup_{t\in[0,T)}\Sigma_t$.

We consider the sphere as our 3D domain. 
%First, we set the centre as $\textbf{c}$, radius as $r$. For any point $x\in\Omega$, fix some positive integer m, we define $\rho$ and $\rho_{lj}$ as follow:
%$$\rho(x)=k\{1-\sup_{l,j}[(\textbf{c}+(\sin(l\pi/m)\cos(j\pi/m),\sin(l\pi/m)\sin(j\pi/m),\cos(l\pi/m)^T)\cdot (x-\textbf{c}))/r]\},$$
%$$\rho_{lj}(x)=k\{1-(\textbf{c}+(\sin(l\pi/m)\cos(j\pi/m),\sin(l\pi/m)\sin(j\pi/m),\cos(l\pi/m)^T)\cdot (x-\textbf{c}))/r\}.$$\\
First, we divide the 3D domain into $k$ layers by the function $\rho$. We set $t_i=i-1$, $i=1,\dots,k+1$. Then $\{\sigma_i=\rho^{-1}([t_i,t_{i+1}))\}_{i=1}^k$ is partition of 3D domain. We then separate the layers into $m$ small domains which is similar to a 2D disc by the function $\rho_{ij}$. Similarly, $\{\sigma_i^j=\rho_{ij}^{-1}([t_i,t_{i+1}))\}_{j=1}^m$ is the set of the small domain which the union of the set returns the layer $\{\rho^{-1}([t_i,t_{i+1}))\}$. Figure \ref{fig: layer} shows an example of dividing the first and second layers of 3D domain into small disks by yellow lines. The colored layers are different 2D disks. This separation is done for several directions .The height of this disc is the grid size $h$ so each disc contains about one layer of grid points, that is, $k=r/h$.

\begin{figure}[h!]
\centering
\includegraphics[scale=0.72]{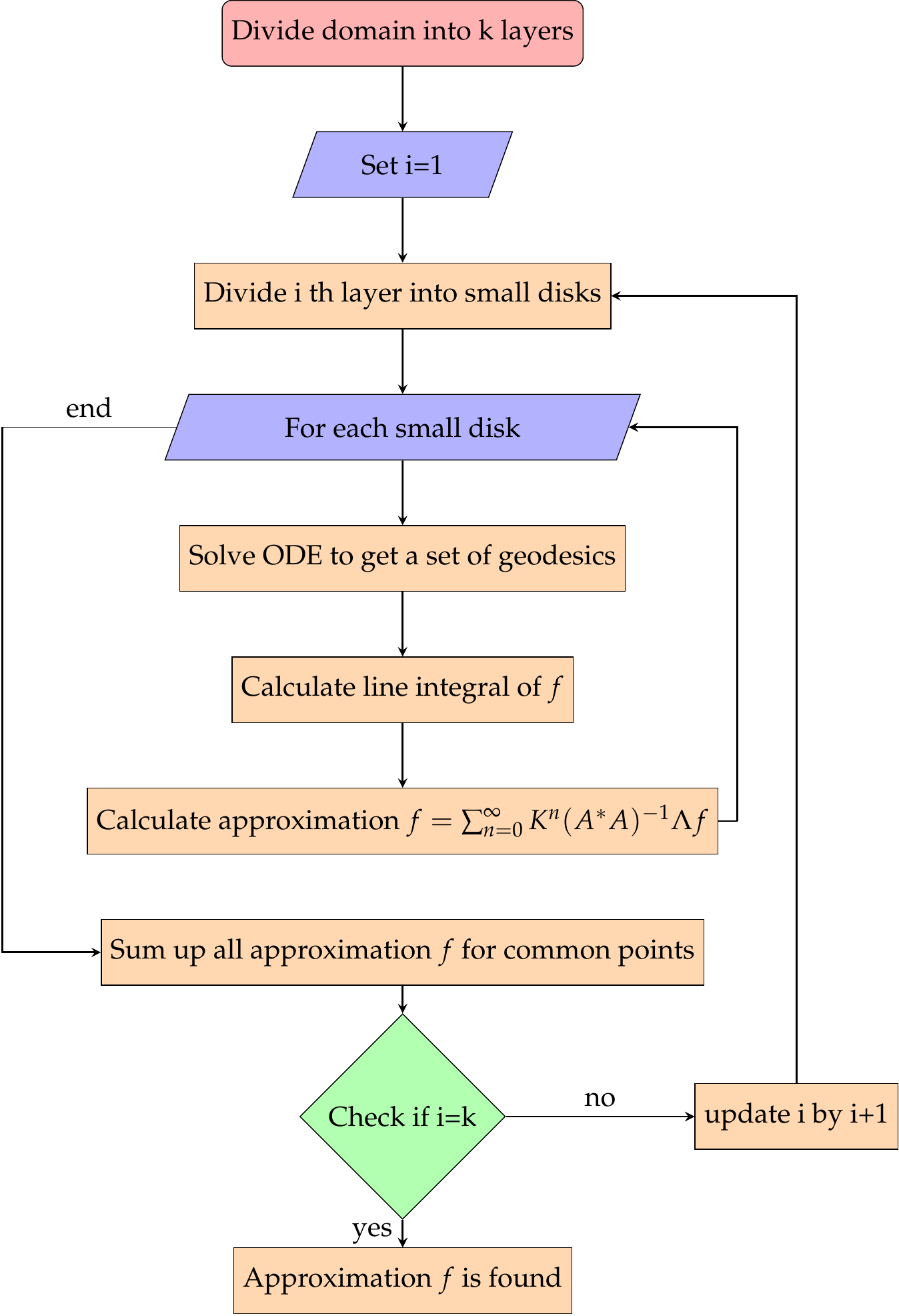}
\caption{The layer stripping algorithm.}
\label{flow}
\end{figure}

Second, we choose different initial points $X^{(0)}_i$ of the geodesics according to the discs. These initial points set at the middle of the boundary of the disc and also form a circle. Then we solve the system (\ref{eq:ODE}) to get the geodesic and calculate the line integral of $f$ for each geodesic $X^j_g$.
Third, for each small domains, we use the regularized version of the 
formula $f=\sum^{\infty}_{n=0}K^n (A^*A)^{-1}\Lambda f$ to calculate the approximation value of $f$ at each fine grid point $z_i$. Also, for each small domains, the calculation of $f$ in the inner layer use the result of that of outer layers.  
Finally, we add up all the approximation values of $f$ of each small domains. If there are some common points between different small domains, then we take the average of approximation values.
The whole algorithm is illustrated in the flow chart in Figure~\ref{flow}.

%\newpage
%%%%%%%%%%%%%%%%%%%%%%%%%%%%%%%%%%%%%%%%%%%%%%
%4           Numerical experiments                									      %
%%%%%%%%%%%%%%%%%%%%%%%%%%%%%%%%%%%%%%%%%%%%%%
\section{Numerical experiments}
In this section we demonstrate the performance of our method using several test examples.  The domain $\Omega$ is a sphere with center $(0.5,0.5,0.5)$ and radius $0.4$. To solve the 
system (\ref{eq:ODE}) to get the 
geodesic curves, we applied the classical Runge-Kutta method of 4th order. Also, for the calculations of the error operator $K$ in (\ref{eq:Neumann}), 
the regularization parameter is chosen as $\delta = 0.2$.

\subsection{Accuracy tests} \label{experiment}
In this section, we illustrate the performance of our reconstruction method using 
several unknown functions $f(x)$. 
The reconstruction is defined on a fine grid
with grid size $h=0.02$. We assume that the speed $c$ is chosen as  $$c(x,y,z)=1+0.3\times \cos(\sqrt{(x-0.5)^2+(y-0.5)^2+(z-0.5)^2}).$$
We consider the following five test cases for this experiment:\\
$f_1=0.01+\sin(2\pi(x+y+z)/10)$, \\$f_2=0.01+\sin(2 \pi(x+y)/10)+\cos(2\pi z/20)$, \\$f_3=x+y^2+z^2/2$, \\$f_4=1+6x+4y+9z+\sin(2 \pi(x+z))+\cos(2\pi y)$, \\$f_5=x+e^{y+z/2}$.

\begin{figure}[ht]
    \centering
    \begin{subfigure}[b]{0.45\textwidth}
        \includegraphics[width=\textwidth]{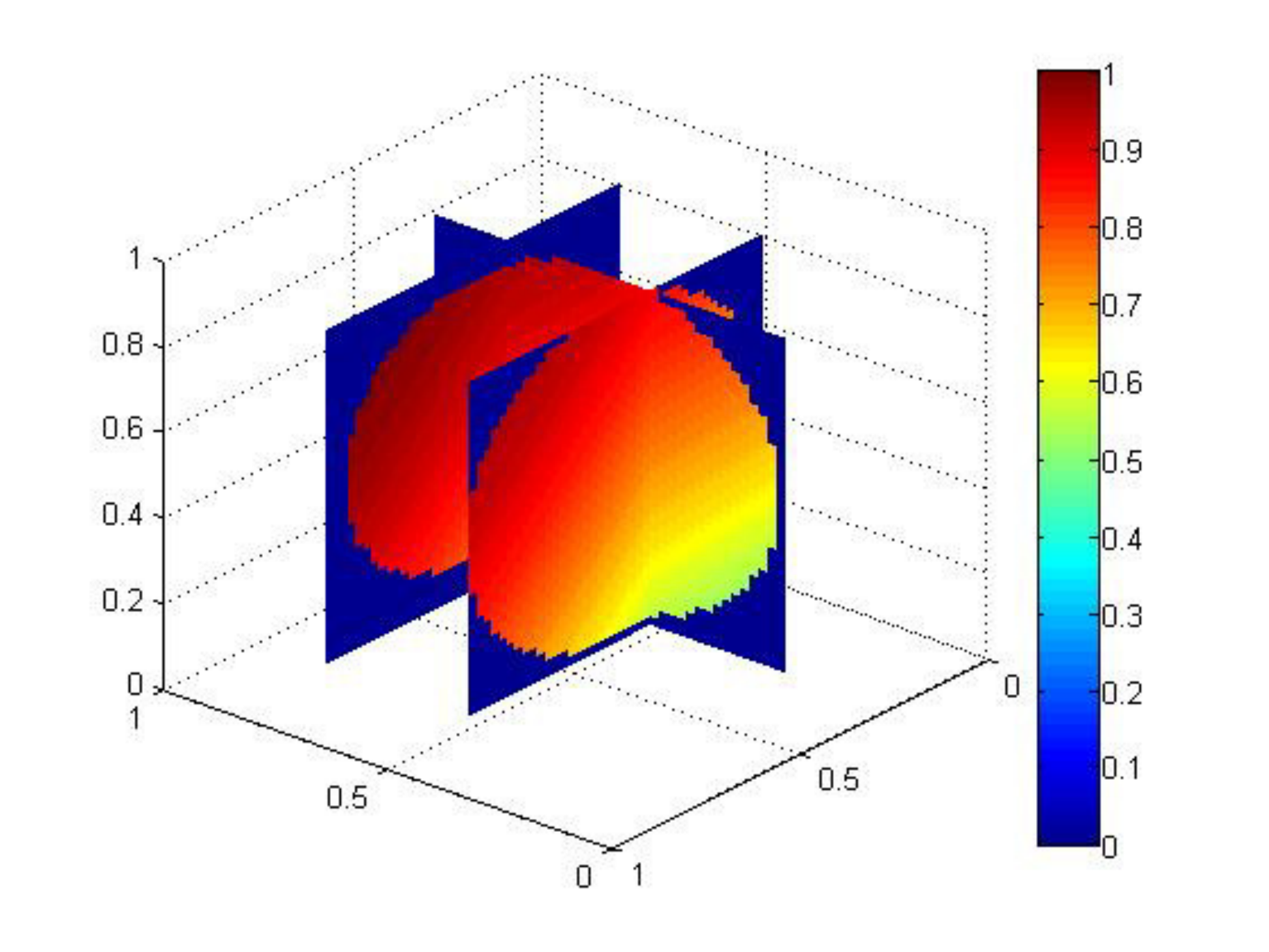}
        \caption{exact solution for $f_1$}
        \label{fig:true solution}
    \end{subfigure}
    ~ %add desired spacing between images, e. g. ~, \quad, \qquad, \hfill etc. 
        \begin{subfigure}[b]{0.45\textwidth}
        \includegraphics[width=\textwidth]{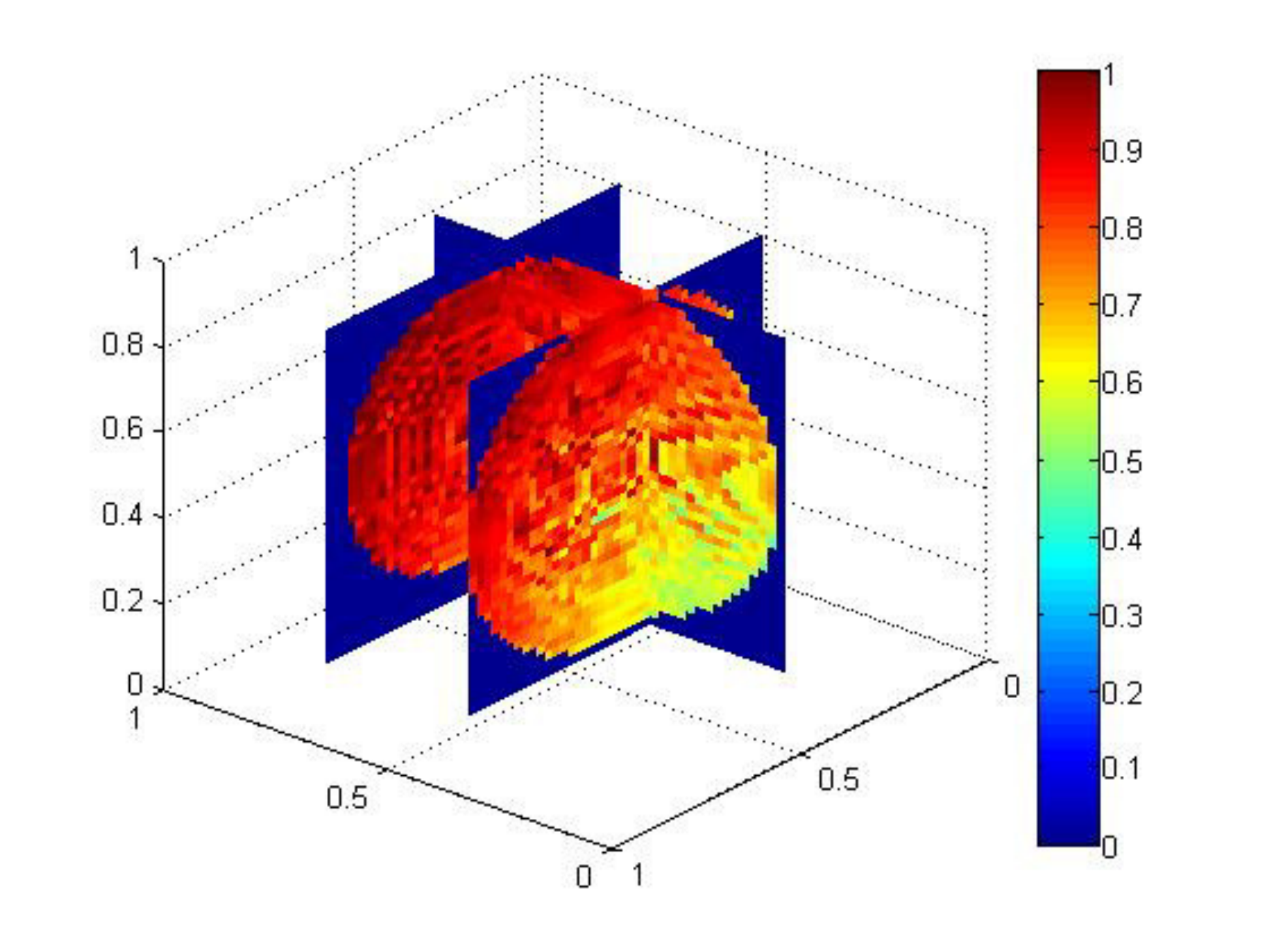}
        \caption{approximate solution for $f_1$}
        \label{fig:approximate solution}
    \end{subfigure}

    \begin{subfigure}[b]{0.45\textwidth}
        \includegraphics[width=\textwidth]{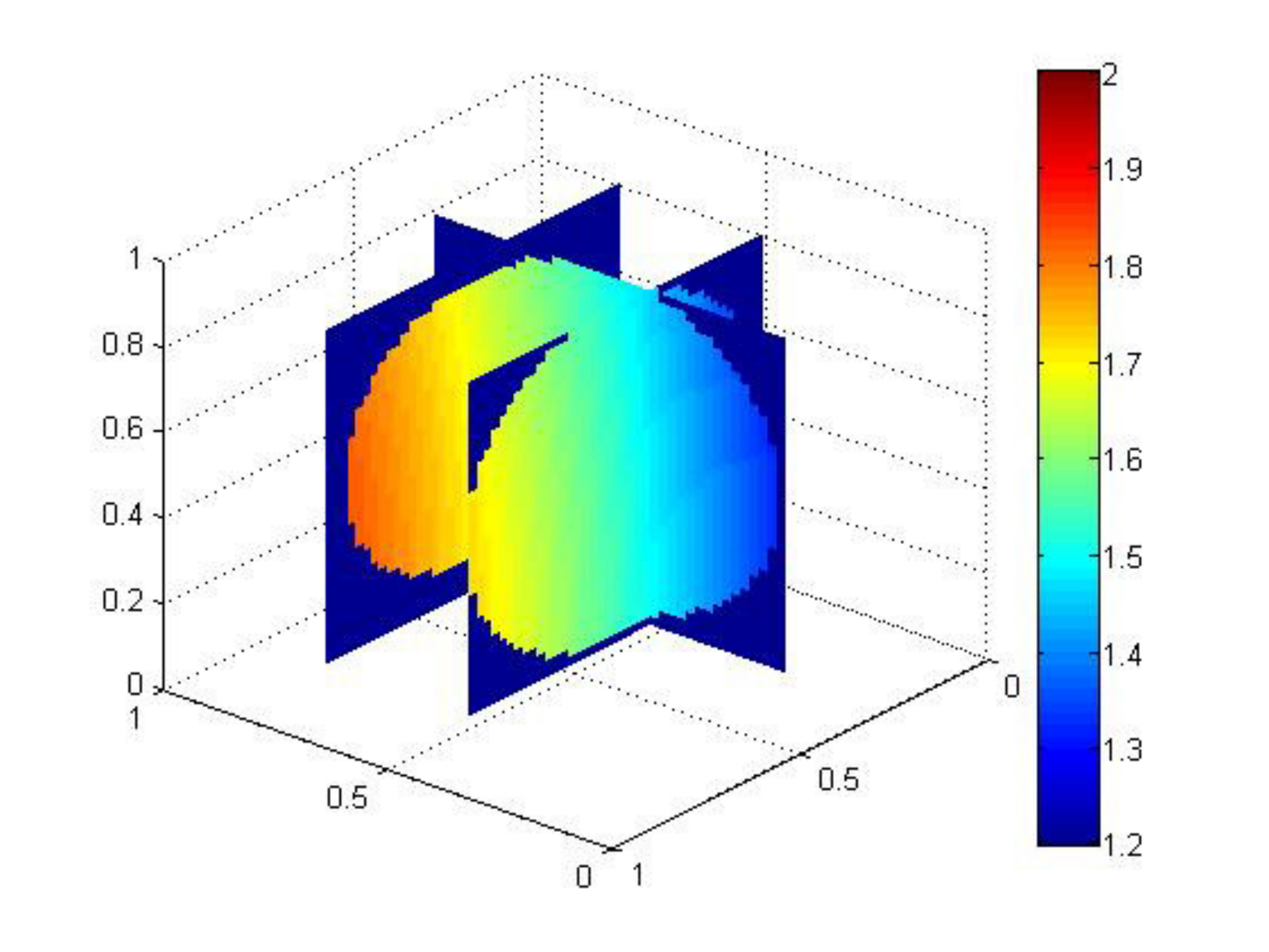}
        \caption{exact solution for $f_2$}
        \label{fig:true solution}
    \end{subfigure}
    ~ %add desired spacing between images, e. g. ~, \quad, \qquad, \hfill etc. 
        \begin{subfigure}[b]{0.45\textwidth}
        \includegraphics[width=\textwidth]{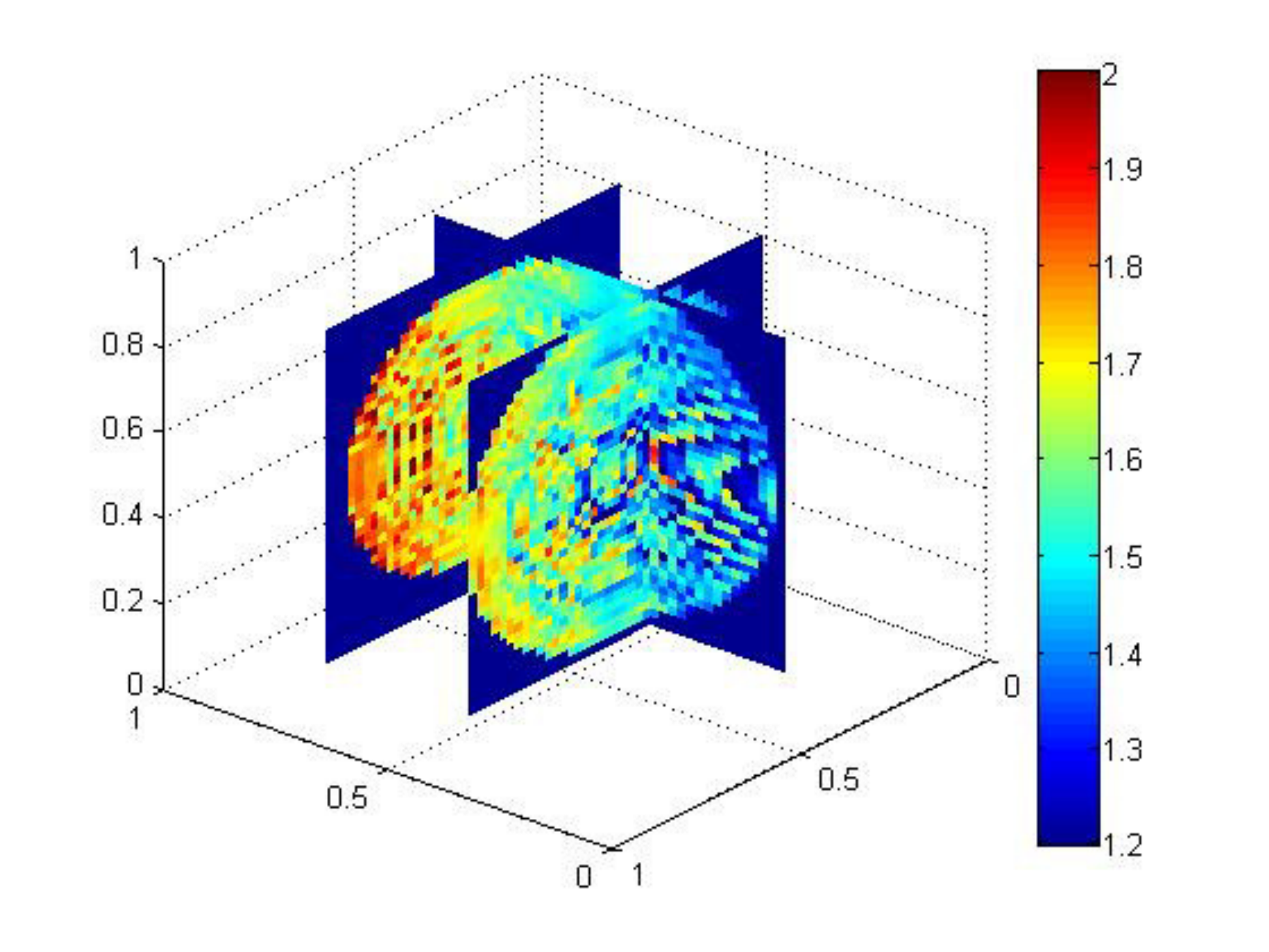}
        \caption{approximate solution for $f_2$}
        \label{fig:approximate solution}
    \end{subfigure}

    \begin{subfigure}[b]{0.45\textwidth}
        \includegraphics[width=\textwidth]{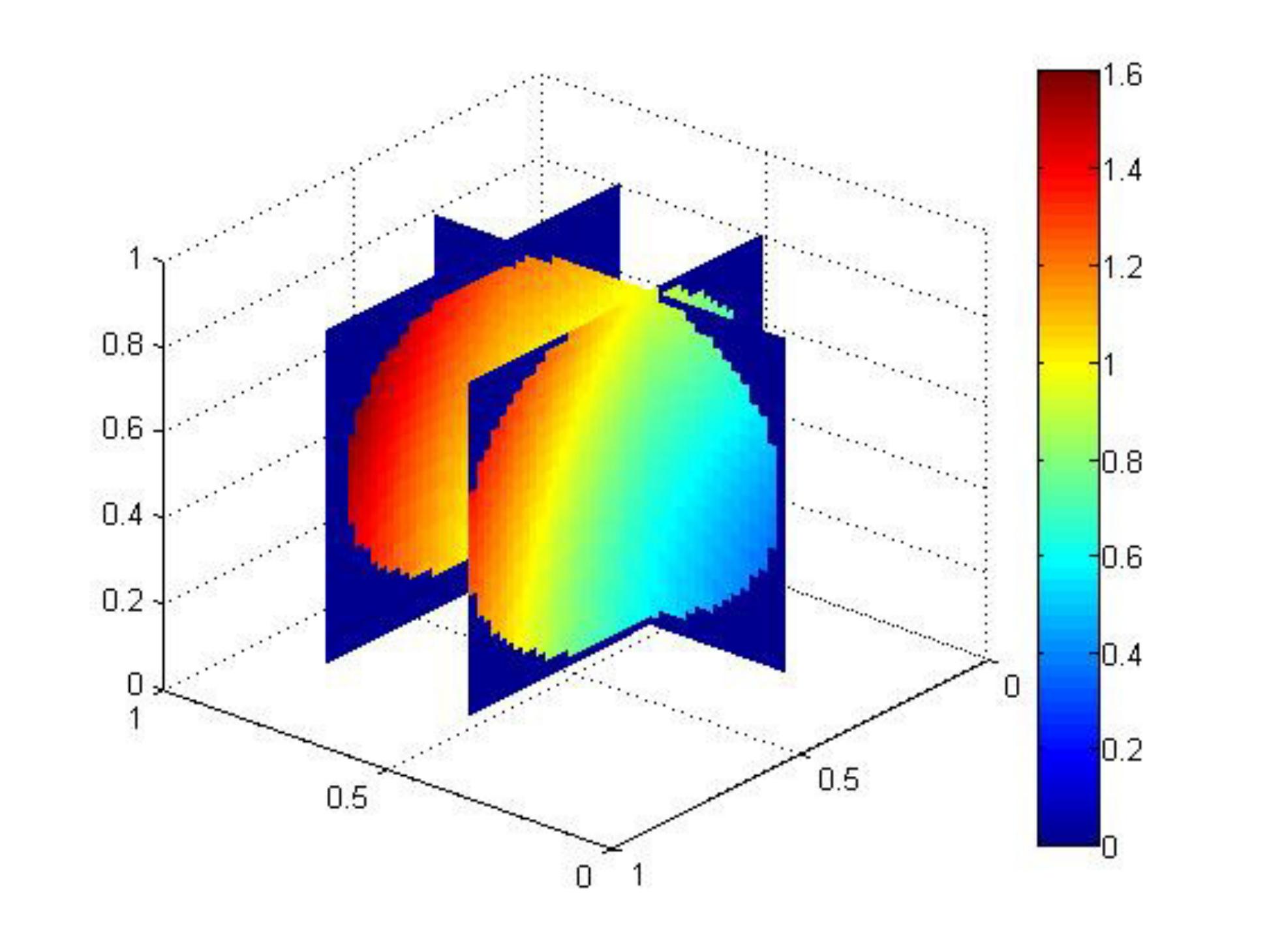}
        \caption{exact solution for $f_3$}
        \label{fig:true solution}
    \end{subfigure}
    ~ %add desired spacing between images, e. g. ~, \quad, \qquad, \hfill etc. 
        \begin{subfigure}[b]{0.45\textwidth}
        \includegraphics[width=\textwidth]{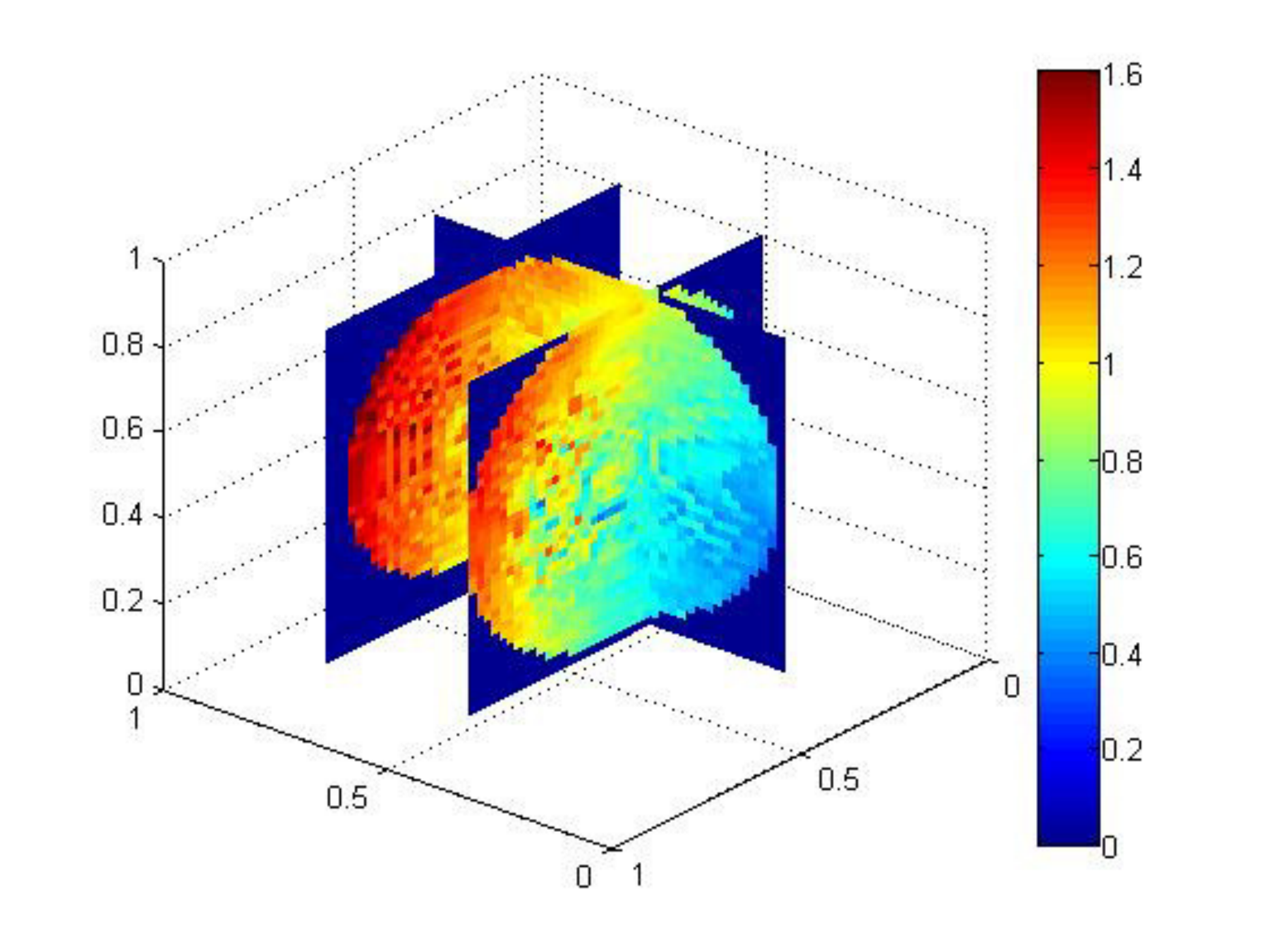}
        \caption{approximate solution for $f_3$}
        \label{fig:approximate solution}
    \end{subfigure}
     \caption{Graphs of true and reconstructed solutions for functions $f_1$, $f_2$, $f_3$.}\label{fig:solution1}
\end{figure}

\begin{figure}[ht]
    \centering
    \begin{subfigure}[b]{0.45\textwidth}
        \includegraphics[width=\textwidth]{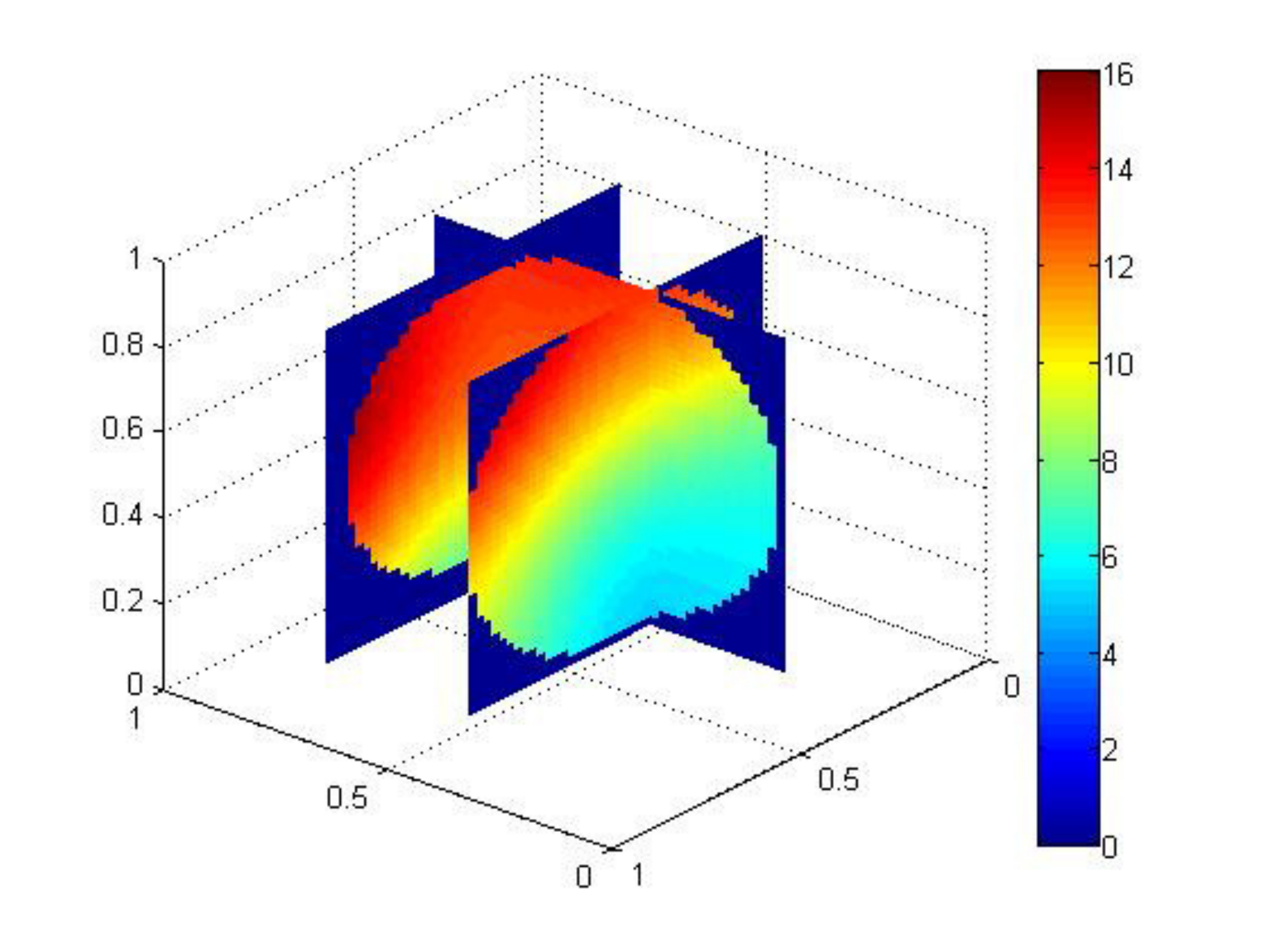}
        \caption{exact solution for $f_4$}
        \label{fig:true solution}
    \end{subfigure}
    ~ %add desired spacing between images, e. g. ~, \quad, \qquad, \hfill etc. 
        \begin{subfigure}[b]{0.45\textwidth}
        \includegraphics[width=\textwidth]{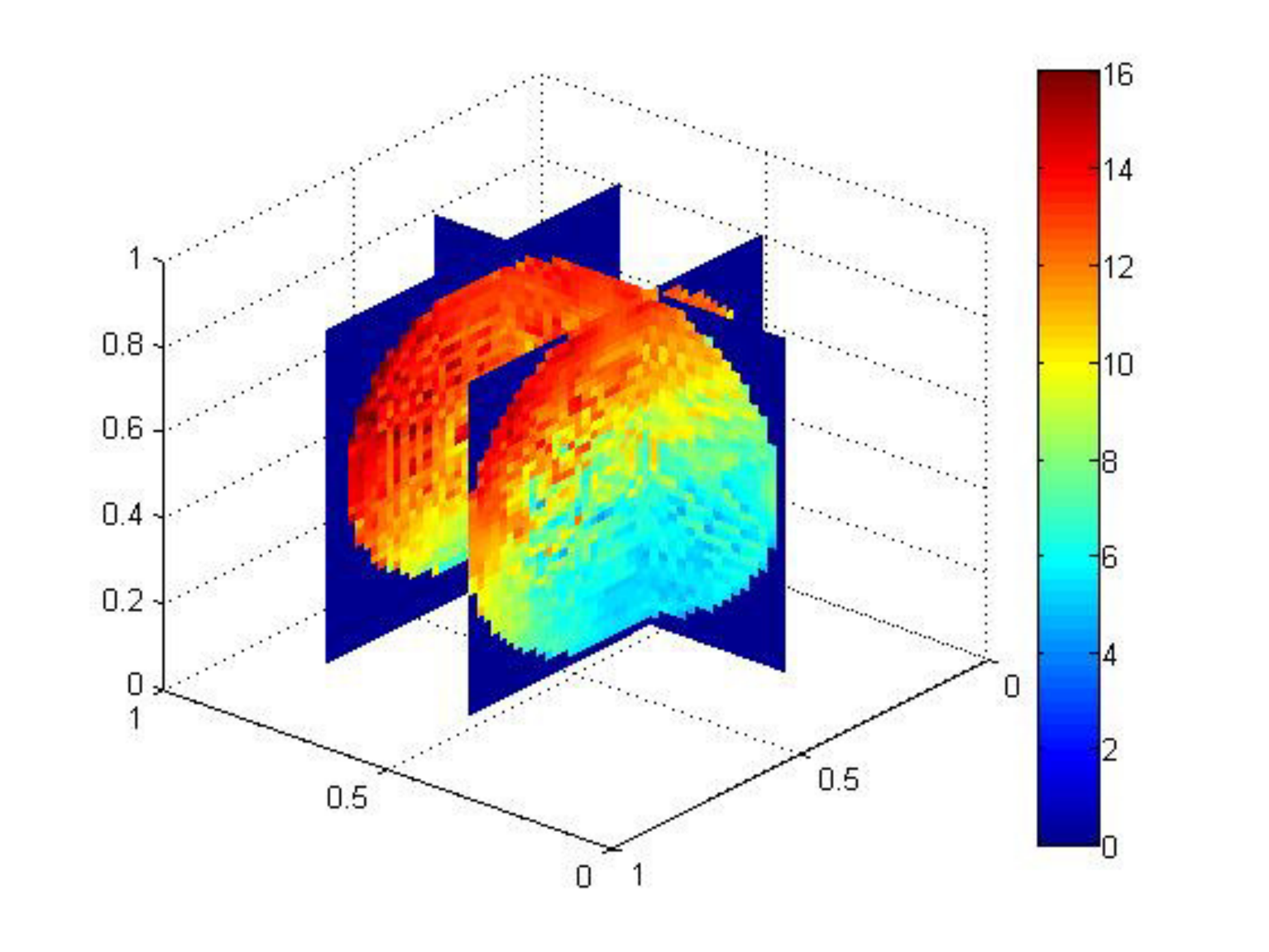}
        \caption{approximate solution for $f_4$}
        \label{fig:approximate solution}
    \end{subfigure}
    
       \begin{subfigure}[b]{0.45\textwidth}
        \includegraphics[width=\textwidth]{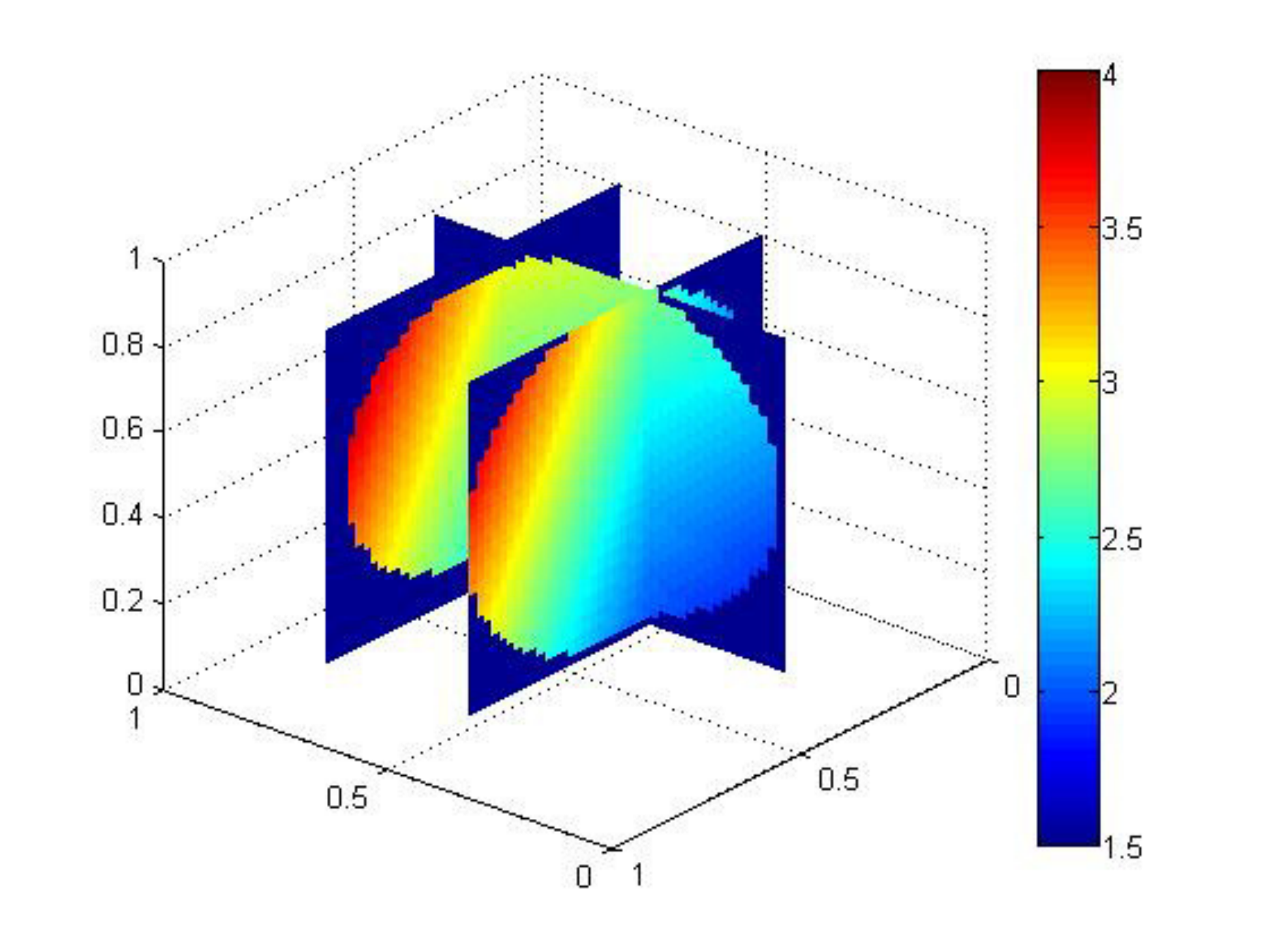}
        \caption{exact solution for $f_5$}
        \label{fig:true solution}
    \end{subfigure}
    ~ %add desired spacing between images, e. g. ~, \quad, \qquad, \hfill etc. 
        \begin{subfigure}[b]{0.45\textwidth}
        \includegraphics[width=\textwidth]{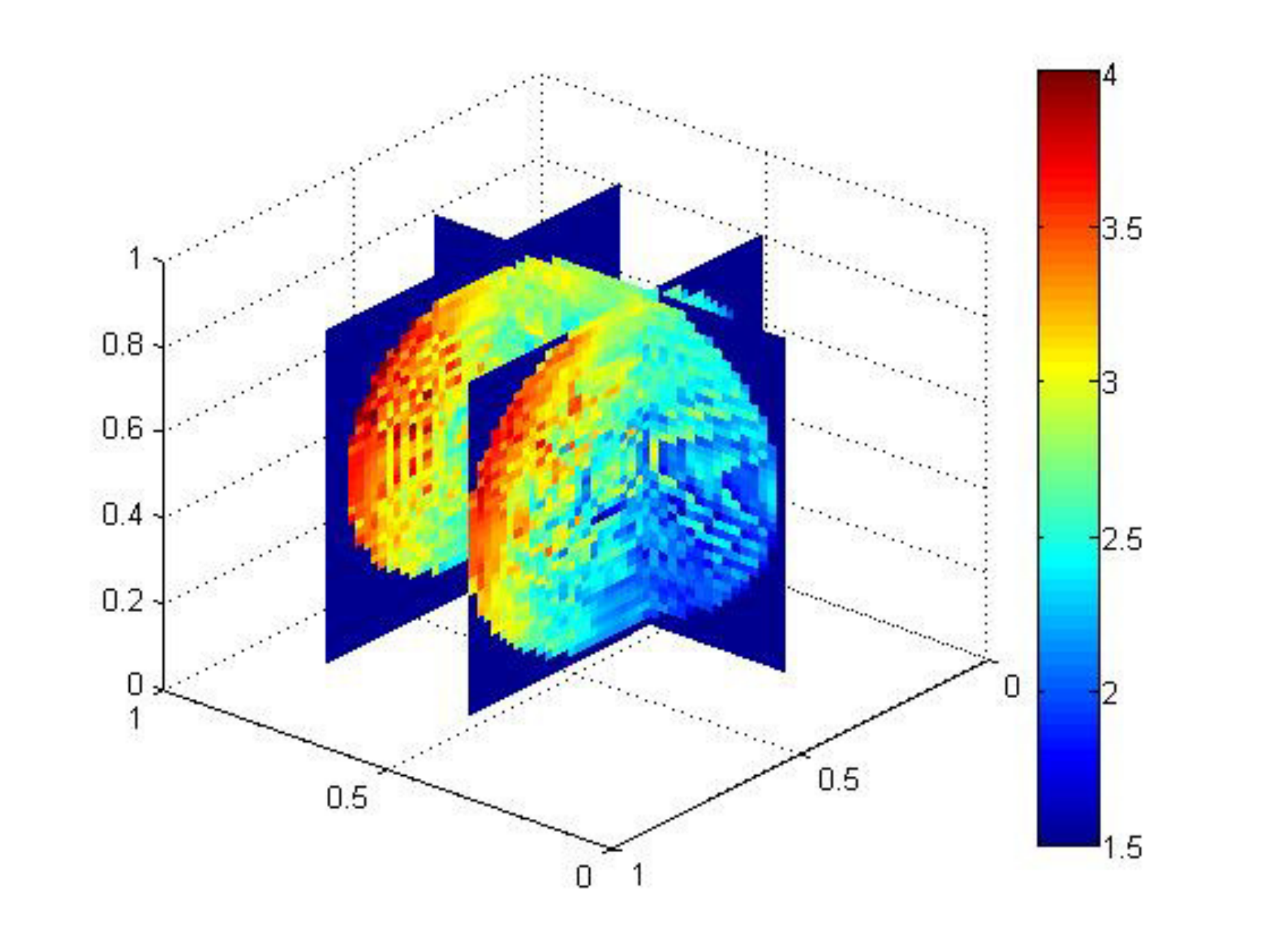}
        \caption{approximate solution for $f_5$}
        \label{fig:approximate solution}
    \end{subfigure}
    \caption{Graphs of exact and reconstructed solutions for functions $f_4$, $f_5$.}\label{fig:solution2}
\end{figure}

We will apply our algorithm to reconstruct each of these five functions. 
Table \ref{fig:relative error} shows the $L^2$-norm relative errors of these five test cases. 
For these cases, we compute the reconstructions using up to the first $5$ terms in the series (\ref{eq:Neumann}). 
We see clearly that the results are improved as we add more terms in the Neumann series (\ref{eq:Neumann}). 
%After calculating the third power of the error operator $K$, the relative error attained minimum. Also, figure 3 shows that if the test function  $f$ is about linear like $f_1$ and $f_2$, the error is smaller.\\
Figures \ref{fig:solution1} and \ref{fig:solution2} show the exact and reconstructed solutions of these test functions with $n=4$. We observe that 
our scheme is able to produce very good reconstructions. 
%that the boundary regions have better approximation than the central regions.

 \begin{table}[ht]
    \centering
 \begin{center}
  \begin{tabular}{| l | c | c | c | c | c |}
    \hline
    relative error & $f_1$ & $f_2$ & $f_3$ & $f_4$ & $f_5$ \\ \hline
    n=0 &  47.28\% &  47.45\% &  46.83\% &  46.97\% &  47.18\% \\ \hline
    n=1 &  23.89\% &  24.10\% &  23.43\% &  23.61\% &  23.76\% \\ \hline
    n=2 &  13.09\%  &  13.26\% &  13.03\% &  13.23\% &  13.00\%\\ \hline
    n=3 &  8.53\% &  8.52\% &  9.32\% &  9.48\% &  8.53\% \\ \hline
    n=4 &  6.99\% &  6.74\% &  8.71\% &  8.80\% &  7.14\% \\ \hline
   % n=5 &  6.70\% &  6.22\% &  9.10\% &  9.14\% &  6.95\% \\ \hline
    % n=6 &  6.80\% &  6.12\% &  9.72\% &  9.73\% &  7.12\% \\ \hline
    %n=7 &  7.00\% &  6.15\% &  10.39\% &  10.36\% &  7.40\% \\ \hline
  \end{tabular}
\end{center}
 \caption{Relative errors for reconstructing the functions $f_i$, $i=1,\cdots , 5$.}\label{fig:relative error}
\end{table}

\subsection{Accuracy tests with noisy data}
\begin{figure}[ht]
    \centering
    \begin{subfigure}[b]{0.45\textwidth}
        \includegraphics[width=\textwidth]{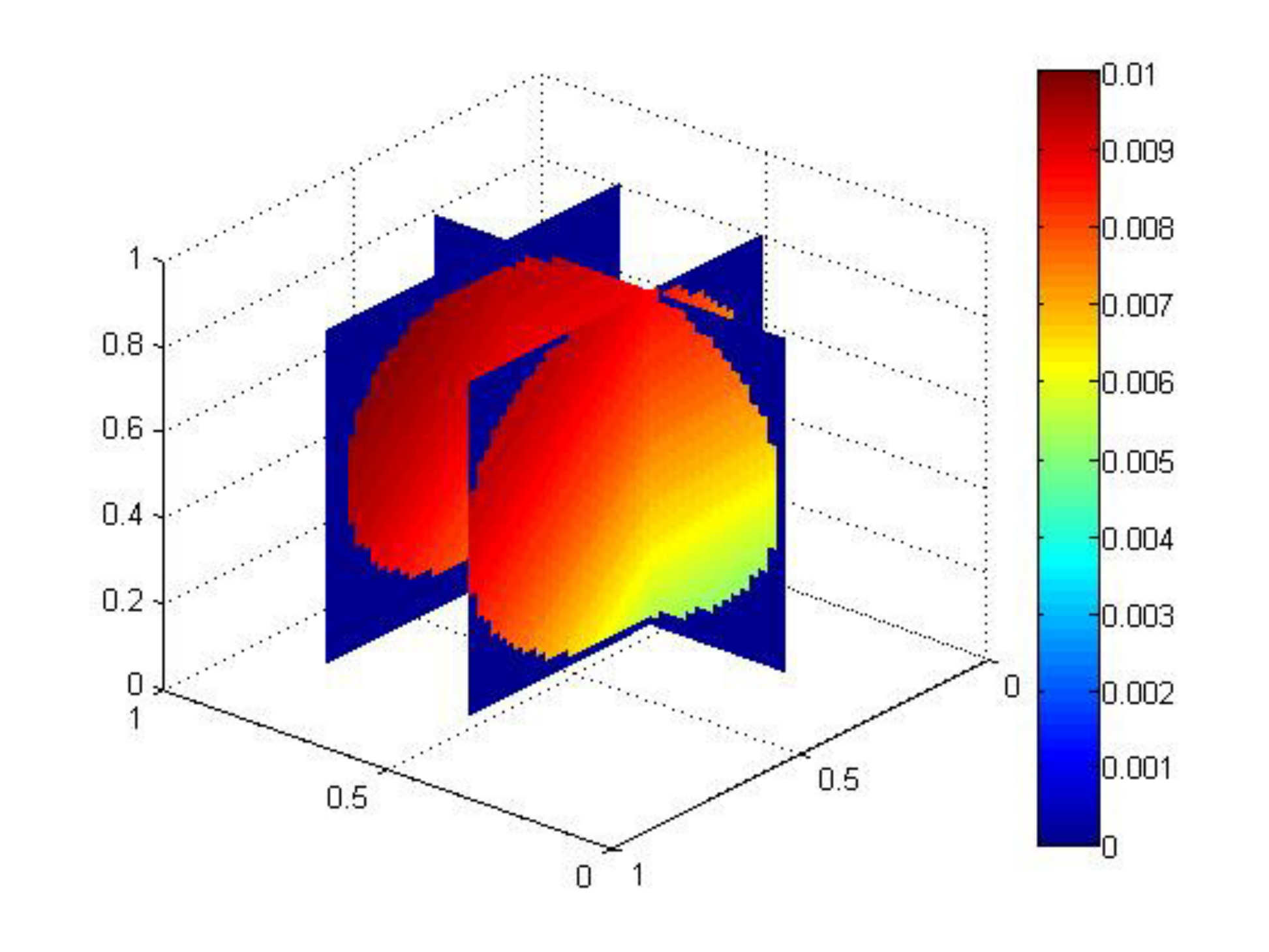}
        \caption{exact solution }
        \label{fig:true solution}
    \end{subfigure}
    ~ %add desired spacing between images, e. g. ~, \quad, \qquad, \hfill etc. 
        \begin{subfigure}[b]{0.45\textwidth}
        \includegraphics[width=\textwidth]{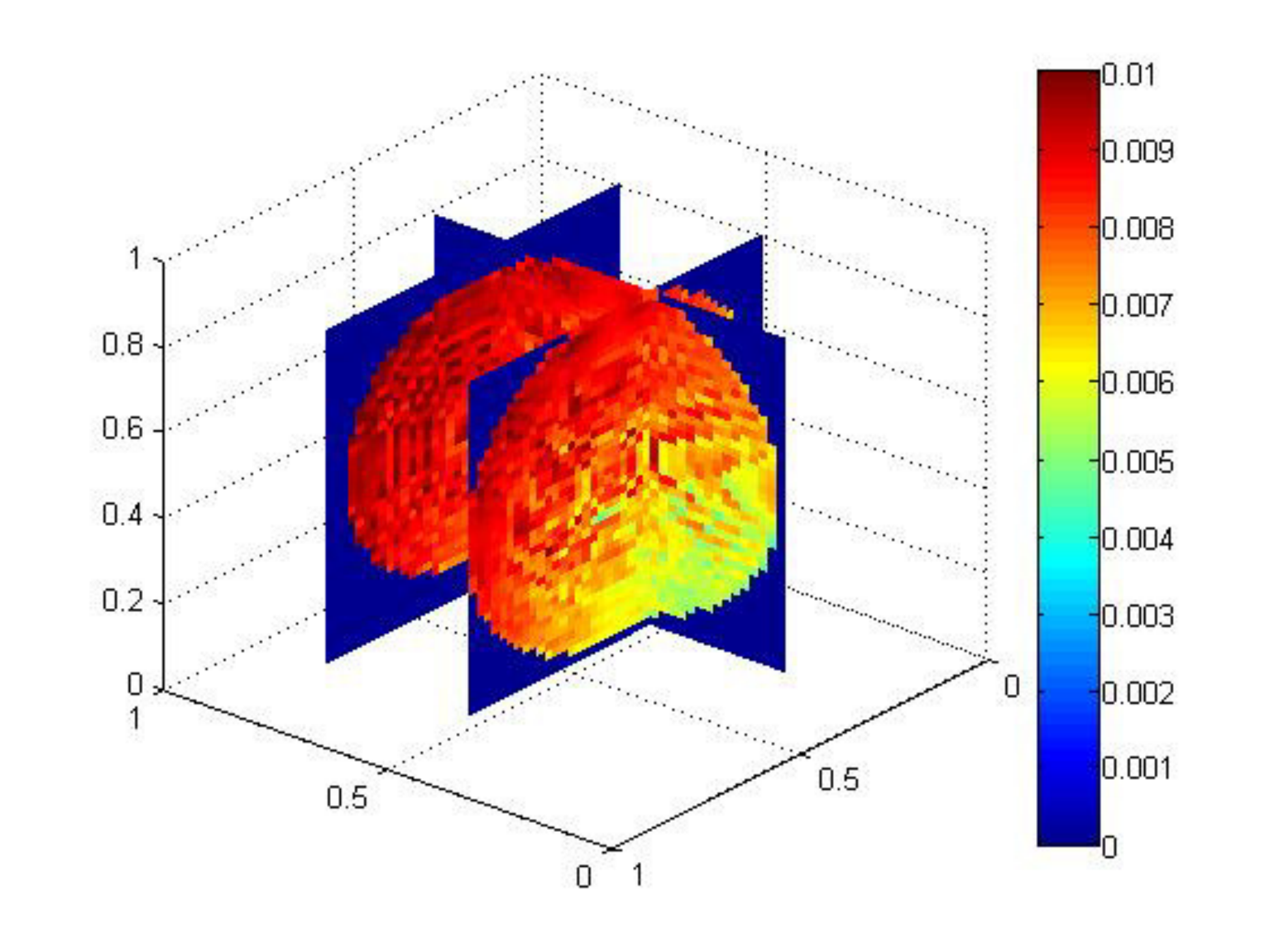}
        \caption{approximate solution for $g_1$}
        \label{fig:approximate solution}
    \end{subfigure}

    \begin{subfigure}[b]{0.45\textwidth}
        \includegraphics[width=\textwidth]{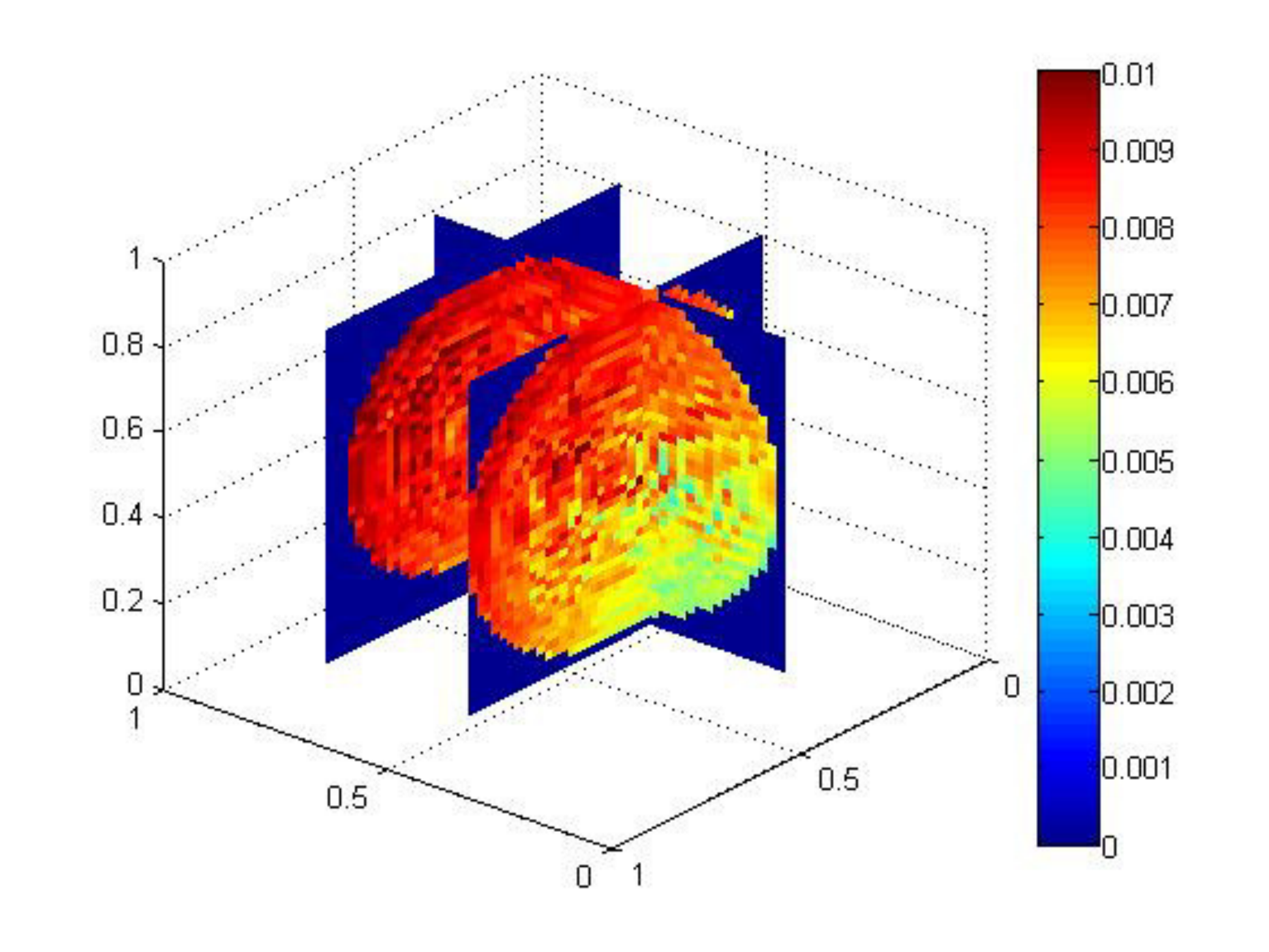}
        \caption{approximate solution for $g_2$}
        \label{fig:true solution}
    \end{subfigure}
   % ~ %add desired spacing between images, e. g. ~, \quad, \qquad, \hfill etc. 
      %  \begin{subfigure}[b]{0.45\textwidth}
       % \includegraphics[width=\textwidth]{randfine.pdf}
        %\caption{approximate solution for $g_3$}
        %\label{fig:approximate solution}
    %\end{subfigure}

        \caption{Graphs of exact and reconstructed solutions using the exact data $g_1$ (error = $6.72\%$)
        and the noisy data $g_2$ (error = $8.50\%$).}\label{fig:solution noisy}
\end{figure}

Next, we show some numerical tests using noise contaminated data $\Lambda^\epsilon g$. We use the same speed $c(x,y,z)$ 
and grid size $h=0.02$ as in Section \ref{experiment} and perform this test by reconstructing the function
$$ g(x,y,z)=0.01\sin(2\pi(x+y+z)/10).$$
The measurement data has been contaminated by uniformly distributed noise $\epsilon$,
$$\Lambda^\epsilon g := \Lambda g +\epsilon$$
with relative error $|\epsilon|/|\Lambda g| = 0.05$ (i.e. $5\%$ noise), where $\epsilon$ is a random function.
In this simulation, we use $n=5$, that is, we use $6$ terms in the Neumann series (\ref{eq:Neumann}). 
Figure \ref{fig:solution noisy} shows the exact and reconstructed solutions using the exact data and noisy data. 
By using the exact data (that is, the data with no noise), we obtain a reconstructed solution with $6.72\%$ error,
while using the noisy data, we obtain a reconstructed solution with $8.50\%$ error.
We observe that our scheme is quite robust with respect to some noise in the data.

%It shows that the noise mostly affect the reconstruction of central region.
%Figure \ref{fig: relative error noisy} shows that the relative error of noisy data is $5\%$ larger than that of exact data. Also, if we choose a finer grid size $h=0.016$, the relative error decreases and the absorption term $\delta$ is changed to 0.16 to get the above result. 

  %     \begin{figure}
    %\centering
 %\begin{center}
  %\begin{tabular}{| l | c | c | c | c | c |}
    %\hline
    %relative error & $g_1$ & $g_2$ & $g_3$  \\ \hline
    %n=5 &  6.72\% &  8.50\% &  8.93\%  \\ \hline
  %\end{tabular}
%\end{center}
 %\caption{Tables of relative errors of test functions reconstructing from exact data $g_1$, noisy data $g_2$ and noisy data with finer grid size $g_3$}\label{fig: relative error noisy}
%\end{figure}

 \subsection{Experiments with different speeds}
 
\begin{figure}[ht]
    \centering
    \begin{subfigure}[b]{0.45\textwidth}
        \includegraphics[width=\textwidth]{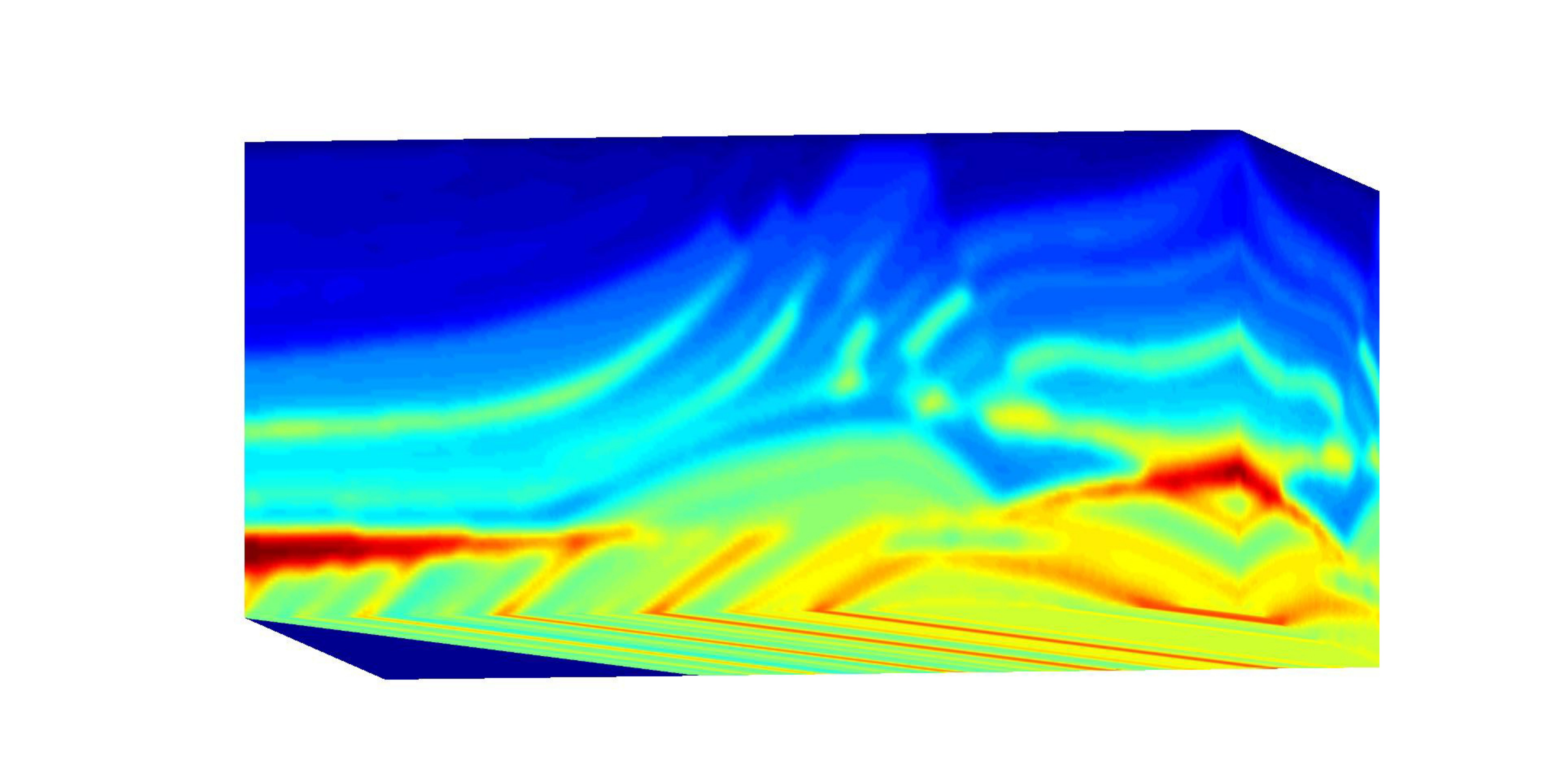}
  \caption{The 3D Marmousi model.}
~
    \end{subfigure}
     \begin{subfigure}[b]{0.45\textwidth}
        \includegraphics[width=\textwidth]{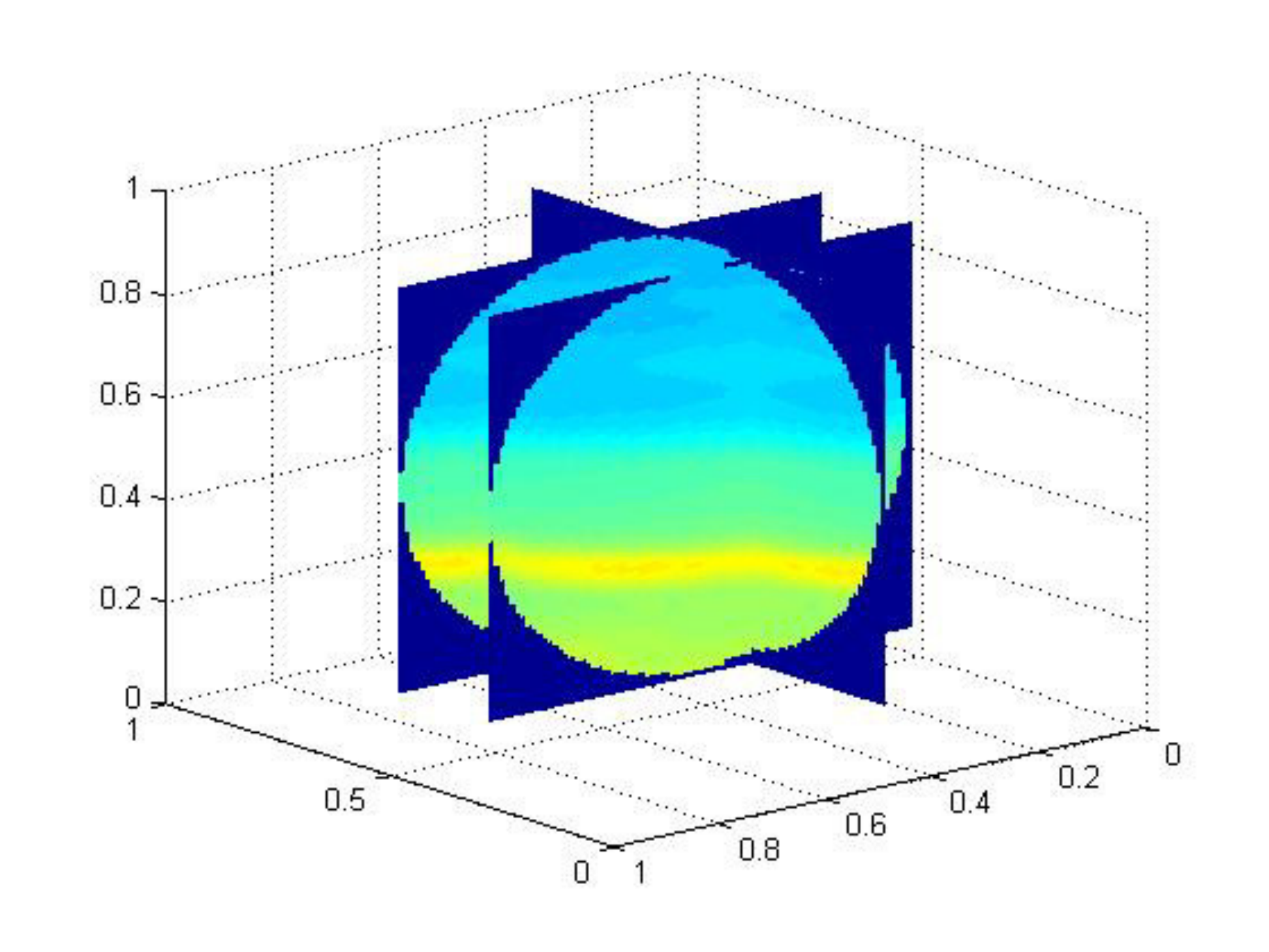}
  
\caption{A spherical section of the Marmousi model as the test speed $c_2$.}
    \end{subfigure}
    ~

     \begin{subfigure}[b]{0.45\textwidth}
        \includegraphics[width=\textwidth]{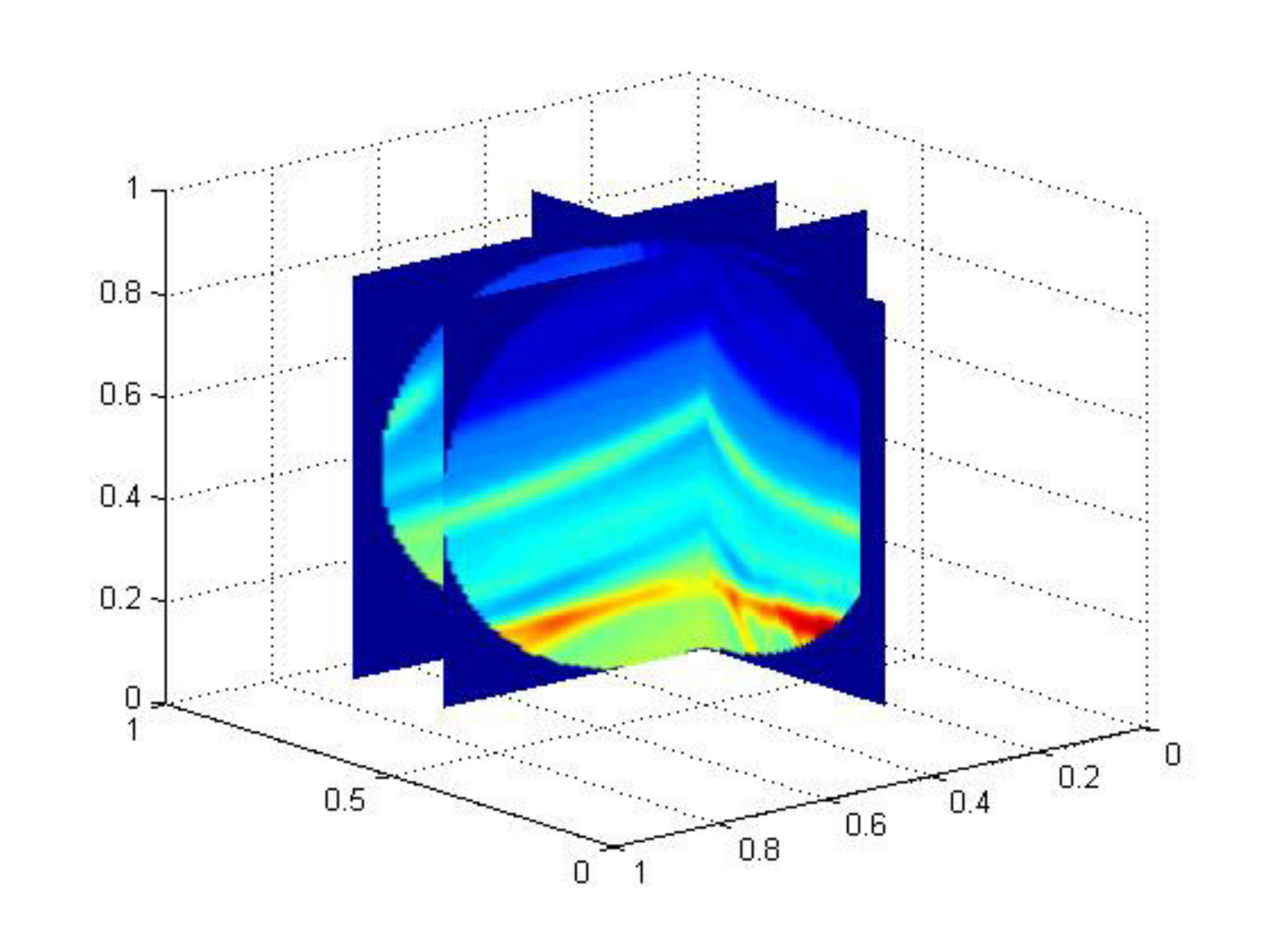}
  
\caption{A spherical section of the Marmousi model as the test speed $c_3$.}
    \end{subfigure}
        \caption{The 3D Marmousi model.}\label{fig: marmous}
\end{figure}

 Finally, we compare the performance of the reconstruction using
different choices of sound speeds. We use grid size $h=0.02$ as above and perform this test by reconstructing the following function
$$ f(x,y,z)=0.01\sin(2\pi(x+y+z)/10).$$
We will test the performance of our algorithm using three choices of sound speeds.
The first speed in our test is defined as 
$$c_1(x,y,z)=1+0.2\sin(3\pi x)\sin(\pi y)\sin(2\pi z).$$
The second and the third
tests are related to the well-known benchmark problem: the Marmousi model, shown in Figure~\ref{fig: marmous}.
In our simulations, we take two spherical sections of the Marmousi model, called $c_2$ and $c_3$, as shown in Figure~\ref{fig: marmous}.

   \begin{figure}[ht]
    \centering
    \begin{subfigure}[b]{0.45\textwidth}
        \includegraphics[width=\textwidth]{realco.pdf}
        \caption{exact solution }
        \label{fig:true solution}
    \end{subfigure}
    ~ %add desired spacing between images, e. g. ~, \quad, \qquad, \hfill etc. 
        \begin{subfigure}[b]{0.45\textwidth}
        \includegraphics[width=\textwidth]{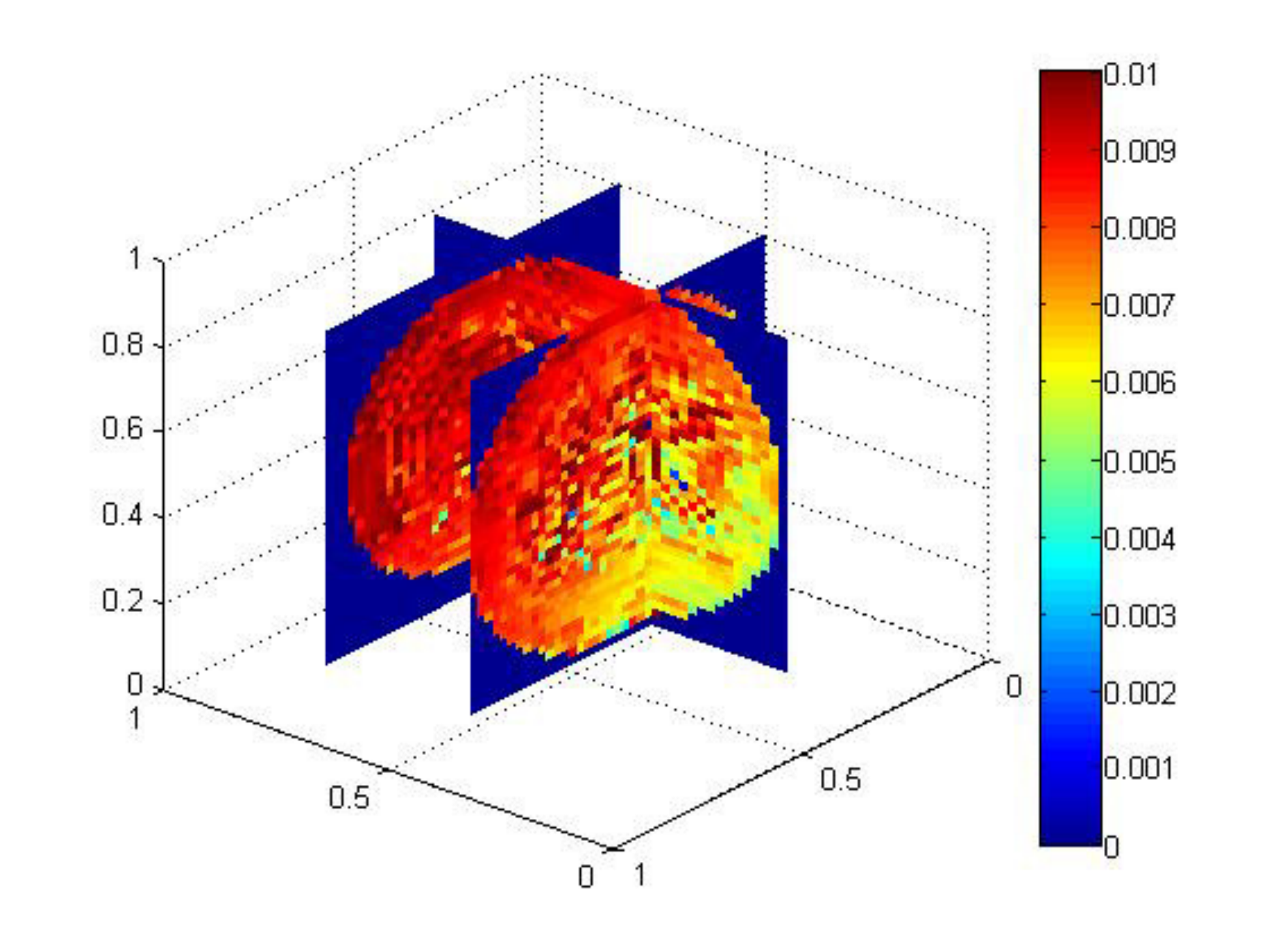}
        \caption{approximate solution for $c_1$}
        \label{fig:approximate solution}
    \end{subfigure}

    \begin{subfigure}[b]{0.45\textwidth}
        \includegraphics[width=\textwidth]{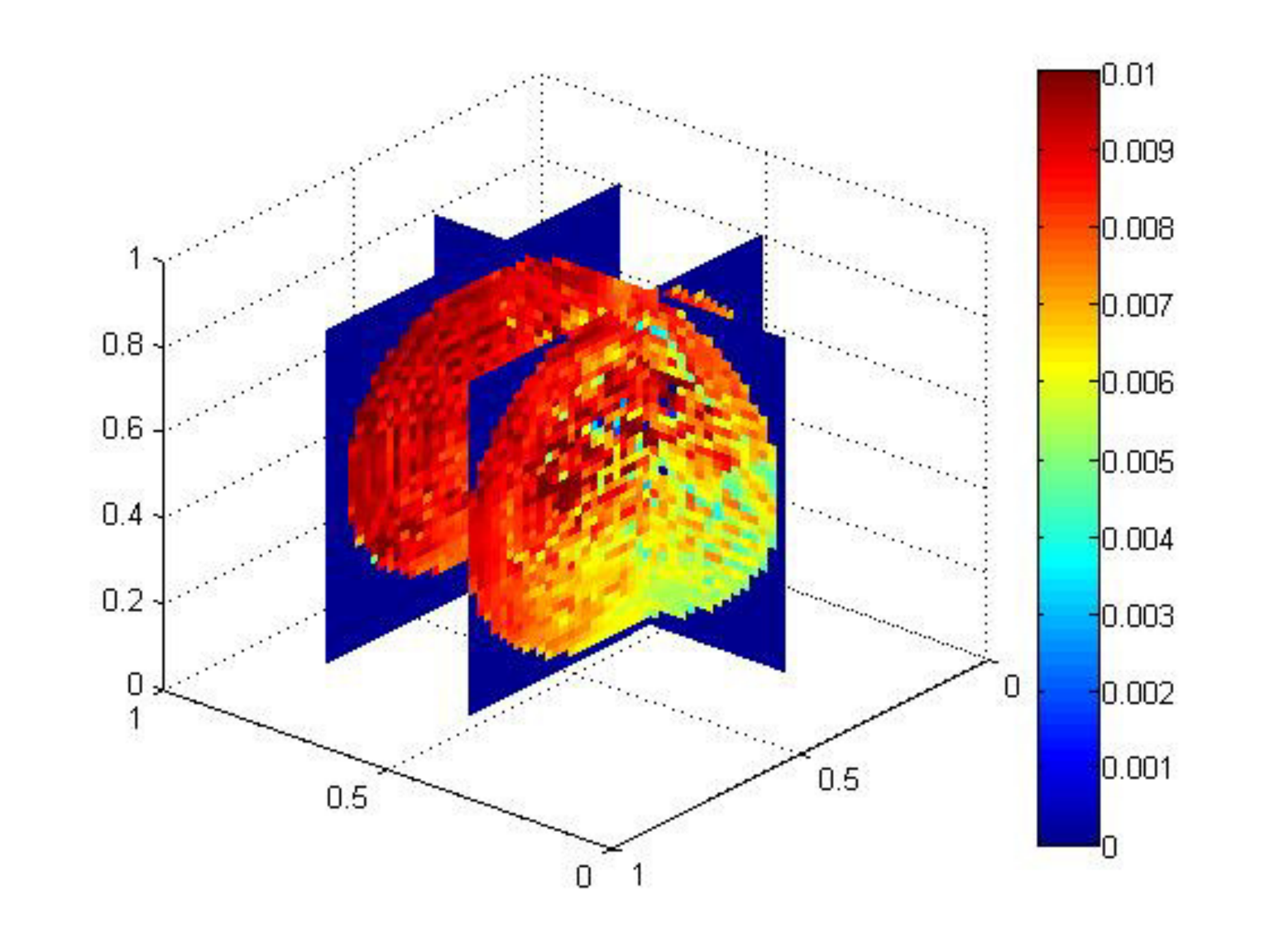}
        \caption{approximate solution for $c_2$}
        \label{fig:true solution}
    \end{subfigure}   
     \begin{subfigure}[b]{0.45\textwidth}
        \includegraphics[width=\textwidth]{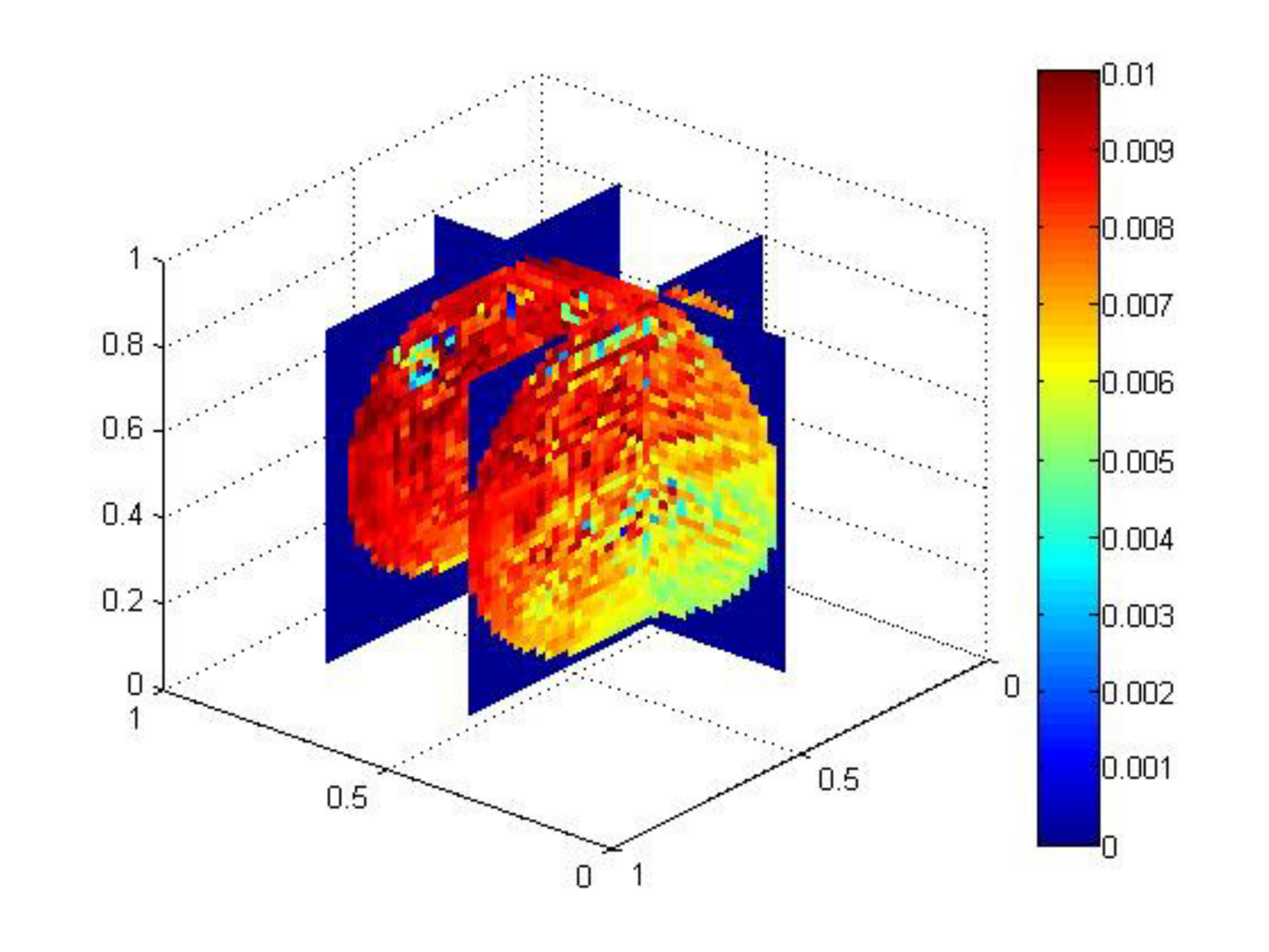}
        \caption{approximate solution for $c_3$}
        \label{fig:true solution}
    \end{subfigure}   

        \caption{Graphs of exact and approximate solutions reconstructing from speed $c_1$ and the Marmousi models $c_2$ and $c_3$.}\label{fig:solution diff}
        \end{figure}
        
               \begin{table}[ht]
    \centering
 \begin{center}
  \begin{tabular}{| l | c | c | c | c | c |}
    \hline
     & $n=0$ & $n=1$ & $n=2$  & $n=3$ & $n=4$ \\ \hline
     relative error for test speed $c_1$ &  44.86\% &  22.88\% &  14.19\% &  11.62\% &  11.47\% \\ \hline %&  12.08\% \\ \hline %12.89 % 13.78
     relative error for test speed $c_2$ &  48.02\% &  25.39\% &  15.71\% &  12.42\% &  11.80\% \\ \hline %&  12.06\% \\ \hline %10.00 %10.25
     relative error for test speed $c_3$ &  58.18\% &  35.41\% &  22.78\% &  16.25\% &  13.39\% \\ \hline %&  12.49\% \\ \hline %12.47 %12.82
  \end{tabular}
\end{center}
 \caption{Relative errors for using the 3 test speeds with grid size $h=0.02$.}\label{fig:converge0}
\end{table}

% \begin{figure}
   % \centering
 %\begin{center}
  %\begin{tabular}{| l | c | c | c | c |}
    %\hline
    % & $h=0.02$ & $h=0.01$ &$h=0.01$ with 5\% noise &$h=0.01$ with function $f_5$  \\ \hline
     %relative error for $n=4$ &  11.80\% &    10.23\% &11.75\% & 10.56\% \\ \hline
  %\end{tabular}
%\end{center}
 %\caption{Tables of relative errors of test functions for Marmousi model with grid size $h=0.01,0.02$}\label{fig:compare}
%\end{figure}

%The 3D model was built by initially
%extending the 2D Marmousi model in the third
%dimension, then shearing this volume successively in three
%directions with a spatially variable degree of shear.
%Successive in-lines look superficially like the original
%model, but they differ in detail, they vary in the cross-line direction. Then, we choose a sphere inside the 3D model as our domain. Figure \ref{fig: marmous} shows an acoustic and a fully elastic 3D shot
%record generated using this model. \\

Figure  \ref{fig:solution diff} shows the exact and reconstructed solutions using the three speeds and $n=4$, that is, five terms in the Neumann series. 
From these figures, we observe very good quality of reconstructions. 
The errors of the reconstructions using the first five terms of the Neumann series are shown in Table~\ref{fig:converge0}.
We observe that the errors decrease as more terms in the Neumann series are used. 
This confirms the theoretical findings. 

%It shows that the method can only recover the boundary region but the central region has difficulty to recover. Also, from figure \ref{fig:converge0}, the approximate solutions coverage for the Neumann series based numerical method when $n$ increase and even we use a more heterogeneous model, we still get  a good result which the $l_2$ relative error is 13.39\% for $n=4$.

Lastly, we perform a similar test for the Marmousi model
with the reconstruction of $f_5=x+e^{y+z/2}$, which is the same as the function in previous section \ref{experiment}.
We also consider the performance with noisy data and the use of a finer grid size $h=0.01$.
The corresponding results using the first $6$ terms in the Neumann series 
are shown in Table \ref{fig:converge}.
We observe that the quality of approximation is improved when more terms are used in the Neumann series.
We also observe that our method is quite robust with respect to $5\%$ noise in the data. 
On the other hand, the performance of the reconstructions of $f$ and $f_5$ is quite similar.

    \begin{table}[ht]
    \centering
 \begin{center}
  \begin{tabular}{| l | c | c | c | c | c |c|}
    \hline
     & $n=0$ & $n=1$ & $n=2$  & $n=3$ & $n=4$ & $n=5$\\ \hline
     for function $f$ &  52.18\% &  28.95\% &  17.53\% &  12.35\% &  10.23\% &  9.42\% \\ \hline %9.15 % 9.09
      for function $f$ and 5\% noisy data &  53.38\% &  30.69\% &  19.47\% &  14.15\% &  11.75\% &  10.68\% \\ \hline %10.21 %10.05
      for function $f_5$ &  52.16\% &  28.95\% &  17.64\% &  12.57\% &  10.56\% &  9.84\% \\ \hline %9.63 %9.63
  \end{tabular}
\end{center}
 \caption{Relative errors using the Marmousi model with finer grid size $h=0.01$.}\label{fig:converge}
\end{table}

%For the points which are away from centre by 0.3, the $l_2$ relative error of Marmousi model is 10.98\% for  grid size $h=0.02$. We also choose grid size $h=0.01$  and absorption term $\delta=0.1$ for the Marmousi model and the $l_2$ relative error is 10.23\% for $n=4$. We also try non-symmetric function $f_5=x+e^{y+z/2}$, which is the same as the function in previous section \ref{experiment}. The resulting $l_2$ relative error is 10.56\% for $n=4$.
%Also, from figure \ref{fig:converge}, the approximate solutions coverage for the Neumann series based numerical method when $n$ increase.

 \section{A new phase space method for traveltime tomography}
 \label{sec:travel}

Our new method is related to the first part of the work and the  Stefanov-€"Uhlmann identity which links two Riemannian metrics with their travel time information. The method first uses those geodesics that produce smaller mismatch with the travel time measurements and continues on in the spirit of layer stripping. We then apply the above reconstruction method to the Stefanov-€"Uhlmann identity in order to solve the inverse problem for reconstructing the index of refraction of a medium. %We demonstrate that our method can recover the isotropic metrics and the three dimensional Marmousi synthetic model. 
The new method is related to the methods in \cite{Chung1,Chung2}. However, we will use the above new method for X-ray transform in solving an integral equation involved
in the inversion process. 
We will discuss the method in detail in the following section. 

\subsection{Mathematical formulation for traveltime tomography}
%We are interested in reconstructing a Riemannian metric in a bounded domain by a set of boundary measurements.
%\subsubsection{Basic setup}
Let $\Omega$ be a  strictly convex bounded domain with smooth boundary $\Gamma=\partial\Omega$ in $\R^d$ and let $g(x) = (g_{ij}(x))$ be a Riemannian metric defined on it. Let $d_g(x,y)$ denote the geodesic distance between $x$ and $y$ for any $x,y\in\Omega$. The inverse problem is whether we can determine the Riemannian metric $g$ up to the natural obstruction of a diffeomorphism by knowing $d_g(x,y)$ for any $x,y\in\Gamma$. While theoretically there are a lot of works addressing this question , we explore the possibility of recovering the Riemannian metric numerically. (see \cite{Beylkin,Croke,Michel,Mukhometov1,Mukhometov2,Pestov2,Sharafutdinov})
We use phase space formulation so that multipathing in physical space is allowed. Rather than the boundary distance function we look at the scattering relation which measures the point and direction of exit of a geodesic plus the travel time if we know the point and direction of entrance of the geodesic. 
The notation mostly follow the previous sections and the previous work \cite{Chung2}.
%Below we formulate the problem more precisely. 
%
%Assume that we have two $C^k$ metrics $g_1$ and $g_2$ with $k\geq 2$ satisfying $g_1=g_2$ in $\R^n-\Omega$. Following \cite{Chung1,Chung2}, we define Hamiltonian $H_g$ by 
%$$ H_g(x,\xi)=\frac{1}{2}(\sum^n_{i,j=1}g^{ij}(x)\xi_i\xi_j-1)$$
%for each $x \in \Omega$ and $ \xi\in\R^n$, where $(g_{ij})^{-1}=(g^{ij})$.
%Let $X^{(0)}=(x^{(0)},\xi^{(0)})$ be a given initial condition, where $x^{(0)}\in\partial\Omega$ and $\xi^{(0)}\in S^{n-1}_g(x^{(0)})$, such that the inflow condition holds,
%$$\sum_{i,j=1}^n g^{ij}(x^{(0)})\xi^{(0)}_i\nu_j (x^{(0)})<0$$
%where $\nu(x)$ is the unit outward normal vector of $\partial\Omega$ at the point $x$ and $\nu_j(x)$ denote the $j$th component of this vector. We define $X_{g_j}(s,X^{(0)})=(x_{g_j}(s,X^{(0)}),\xi_{g_j}(s,X^{(0)}))$ by the solution to the hamiltonian system defined by 
%\begin{flalign*}
%\noindent&\frac{dx}{ds}=\frac{\partial H_{g_j}}{\partial \xi} \qquad , \qquad \frac{d\xi}{ds}=-\frac{\partial H_{g_j}}{\partial x} \\
%\text{with the initial} &\text{ condition}&\\
%&(x^{(0)},\xi^{(0)})=X^{(0)}.
% \end{flalign*}
%The solution $X_g$ defines a geodesic/ray in the phase space, parametrically via $x(s)$ with the co-tangent vector $\xi(s)$ at any point $x(s)$. The parameter $s$ denotes travel time. Thus, we denote the set of geodesics $X_g$ which are contained in $\Omega$ with endpoints at $\partial\Omega$ by $\mathcal{M}_{\Omega}$.

The continuous dependence on the initial data of the solution of the Hamiltonian system is characterized by the Jacobian, 
$$J_{g_j}(s)=J_{g_j}(s,X^{(0)}):=\frac{\partial X_{g_j}}{\partial X^{(0)}}(s,X^{(0)})= \left( \begin{array}{cc}
\frac{\partial x}{\partial x^{(0)}} & \frac{\partial x}{\partial \xi^{(0)}}  \\
\frac{\partial \xi}{\partial x^{(0)}} & \frac{\partial \xi}{\partial \xi^{(0)}} \end{array} \right).$$
It can be shown easily from the definition of $J$ and the corresponding Hamiltonian system (the detail will be derived in the Appendix) that $J_{g_j},j=1,2$, satisfies
\begin{align*}
\frac{dJ}{ds}=MJ, \qquad J(0)=I,
\end{align*}
where, in terms of $H=H_{g_j}$ ($j=1,2$), the matrix $M$ is defined by
$$M=\left( \begin{array}{cc}
H_{\xi,x} & H_{\xi,\xi}\\
-H_{x,x} & -H_{x,\xi} \end{array}\right).$$ 
Following \cite{Chung2}, we use the linearized Stefanov-Uhlmann identity  (the detail will be derived in the Appendix) :
\begin{equation}\label{SUid}
X_{g_1}(t,X^{(0)})-X_{g_2}(t,X^{(0)})\approx\int^t_0J_{g_2}\big(t-s,X_{g_2}(s,X^{(0)})\big)\times\partial_{g_2}V_{g_2}(g_1-g_2)\big(X_{g_2}(s,X^{(0)})\big)ds,
\end{equation}
which is the foundation for our numerical procedure.\\
By the group property of Hamiltonian flows the Jacobian matrix is equal to
\begin{equation*}
J_{g_2}\big(t-s,X_{g_2}(s,X^{(0)}\big)=J_{g_2}\big(t,X^{(0)}\big)\times J_{g_2}\big(s,X^{(0)}\big)^{-1}
\end{equation*}\\
In the case of an isotropic medium,
$$g_{ij}=\frac{1}{c^2}\delta_{ij},$$ 
where $c$ is a function from $\R^d$ to $\R$. Then
$$V_{g_k}=(c^2_k\xi,-(\nabla c_k)c_k|\xi|^2).$$
Hence the derivative of $V$ with respect to $g$, $\partial_{g}V_{g}(\lambda)$ is given by 
$$\partial_{g}V_{g}(\lambda)=(2c\lambda\xi,-(\nabla c\cdot\lambda+\nabla\lambda\cdot c)|\xi|^2).$$
In the case of an isotropic medium with the metric $g_{ij}=\frac{1}{c^2}\delta_{ij}$, we have 
$$M=\left( \begin{array}{cc}
2c\xi(\nabla c)^T & (c)^2I\\
-(\nabla^2 c)c|\xi|^2-(\nabla c)(\nabla c)^T|\xi|^2 & -2c(\nabla c)\xi^T \end{array}\right).$$

\subsection{New phase space method using reconstruction of X-ray transform}
The numerical method is an iterative algorithm based on the linearized Stefanov-Uhlmann identity using a hybrid approach. The metric $g$ is defined on  an underlying Eulerian grid in the physical domain. The integral equation (\ref{SUid}) is recovered by the geodesic X-ray data along a ray for each $X^{(0)}$ in phase space. On the ray, values of $g$ are computed by interpolation from grid point values. Hence, each integral equation along a particular ray yields a linear equation for grid values of $g$ in the neighborhood of the ray in the physical domain. Here is the iterative algorithm for finding $g$.

Let $Z$ be the set  of all grid point $z_j,j=1,2,\dots,p$. Let $X^{(0)}_i\in\mathcal{S}^-,i=1,2,\dots,m$, be the initial locations and directions of those $m$ measurements $X_g(t_i,X^{(0)}_i)\in\mathcal{S}^+$ where $t_i$ is the exit time corresponding to the $i$th geodesic starting at $X^{(0)}_i$. Starting with an initial guess of the metric $g^0$, we construct a sequence $g^n$ as follows.

Define the mismatch vector
$$d^n_i=X_g(t_i,X^{(0)}_i)-X_{g^n}(t_i,X^{(0)}_i)$$
and an operator $\hat{I_i}$ along the $i$th geodesic by
$$\hat{I_i}(g-g^n):=\int^t_0J_{g^n}\big(t-s,X_{g^n}(s,X^{(0)}_i)\big)\times\partial_{g^n}V_{g^n}(g-g^n)\big(X_{g^n}(s,X^{(0)}_i)\big)ds$$
%$$(If^n_i)\lambda=\int^{t_i}_0J_{g^n}\big(t_i-s,X_{g^n}(s,X^{(0)}_i)\big) \times \partial_{g^n}V_{g^n}\big(\lambda)(X_{g^n}(s,X^{(0)}_i)\big)ds.$$
Note that both $d^n_i$ and $\hat{I_i}(g-g^n)$ depend on $X^{(0)}_i$. Then we use the above reconstruction method to recover $\lambda:=g-g^n$ at each grid points by the mismatch vector. Thus, we define a new mismatch vector $\hat{d^n_i}$ from $\hat{I_i}$
$$\hat{d^n_i}:=\hat{I_i}^*d^n_i$$
Also, we define an operator $I$ based on (\ref{SUid}) along the $i$th geodesic
$$I(\lambda)(X_g(t_i,X^{(0)}_i))=\hat{I_i}^*\hat{I_i}(\lambda).$$
Thus, for each $n\geq 0$, we use the reconstruction formula  %(\ref{eq:Neumann}) 
$$\lambda=\sum^{\infty}_{n=0}K^n P(A^*A- \delta \Delta)^{-1}P^*\Lambda (I\lambda),$$
where $P,A$ is defined as the same as the previous section and
$$K = Id - P (A^*A)^{-1} P^*(\Lambda \circ I).$$
Hence we then define, $$g^{n+1}=g^{n}+\lambda.$$
We remark that more details can be found in the Appendix. 

\subsection{Numerical implementations}\label{sec:num}
 In this section, we will briefly explain the details of the numerical implementations of phase space method. First, we will discuss the detail of constructing the mismatch vector $d^n_i$. Then, we will explain the calculation of the line integrals, i.e. the operator $\hat{I_i}$. Finally, the detail of the reconstruction formula and how the metric is updated will be explained. 
 
\subsubsection{Setup of mismatch vector}
The mismatch vector is defined by $d^n_i=X_g(t_i,X^{(0)}_i)-X_{g^n}(t_i,X^{(0)}_i)$. Hence, we need a set of the initial locations and directions $\{X^{(0)}_i\}$. For our settings, we will divide the 3D domain into different layers and thus divide the layer into many small disks. For each disk, we will choose 900 initial locations and directions around the boundary of the disks evenly. From this set of data, we will derive a set of mismatch vectors using the guess $g^n$ and also the observed data $X_g(t_i,X^{(0)}_i)$.

Also, since the method of choosing the initial locations and directions cannot promise the geodesic will remain in the layer. We will eliminate the geodesic which does not remain in the same layer. Then, the remaining set of the mismatch vector will be the data we used to derive the method. 

\subsubsection{Calculation of the line integrals}
The operator $\hat{I_i}$ is defined to calculate the  line integrals, which is defined by
$$\hat{I_i}(g-g^n):=\int^{t_i}_0J_{g^n}\big(t_i-s,X_{g^n}(s,X^{(0)}_i)\big)\times\partial_{g^n}V_{g^n}(g-g^n)\big(X_{g^n}(s,X^{(0)}_i)\big)ds.$$
Since we have the discrete phase space of each geodesic, i.e. $X_{g^n}(s_j,X^{(0)}_i)$ for any $0\leq s_j\leq t_i$, then we can approximate the operator by
$$\hat{I_i}(g-g^n)\approx\sum_{s_j} J_{g^n}\big(t_i-s_j,X_{g^n}(s_j,X^{(0)}_i)\big)\times\partial_{g^n}V_{g^n}(g-g^n)\big(X_{g^n}(s_j,X^{(0)}_i)\big) \,(X'_{g^n}(s_j,X^{(0)}_i)) \, (s_{j}-s_{j-1}).$$
Hence, the operator can be approximated by a matrix and thus the adjoint operator $\hat{I_i}^*$.

\subsubsection{Update of the metric}
After we construct the line integral operators, we will apply the reconstruction formula to compute the update of the metric. Based on the reconstruction formula, it is the infinite sum of Neumann series. In computation, we need to choose some terms of the infinite sum to represent the whole term. Here, we choose the first five terms, since this terms represent the main part of the sum. 

Next, this update of the metric will be completed on each disk. To ensure the metric can recover, we do the update five times which make the successive error to be small. After doing the update for each disks in the same layer, we will compute the final metric of this layer and move on to the next layer. When we compute the final metric of this layer, there are some overlapping regions for different disks. Then we take the mean of this values to calculate the final metric.

\subsection{Numerical experiments}
In this section we demonstrate the performance of our method using several test examples.  The domain $\Omega$ is a sphere with center $(0.5,0.5,0.5)$ and radius $0.4$. To solve the 
system to get the 
geodesic curves, we applied the classical Runge-Kutta method of 4th order. Also, for the calculations of the error operator $K$ , 
the regularization parameter is chosen as $\delta = 0.02$. 

\subsubsection{Constant case and the linear case}
To test our algorithm, we test it on different speeds. First, we test the constant case and the linear case. 
For these two test cases, we divide the 3D domain into 20 layers and we illustrate the performance in the first two outermost  layers. 

For the constant case $g=10$, we have the relative error $0.0004\%$ for first layer and $0.0005\%$ for the second layer. Thus, the speeds are well recovered in this two layer. Also, refer to the Table \ref{fig:converge0} , we can note that the relative error grows after adding more layers.

For the linear case $g=10+0.1\times(x+y+z)$, we have the relative error $0.0727\%$ for first layer and $0.0599\%$ for the second layer. For this case, we have the similar result as the constant case. By observing Table \ref{fig:converge0}, we also note that the first two layers are recovered almost exactly and the errors grow after a few layers.
The fact that there are larger errors in inner layers is because there are less data available for those regions. 

 \begin{table}[ht]
    \centering
 \begin{center}
  \begin{tabular}{| l | c | c | c | c | c |}
    \hline
     & $1^{st}$ layer & $2^{nd}$ layer & $3^{rd}$ layer  & $4^{th}$ layer &$5^{th}$ layer \\ \hline
     relative error for  constant case &  0.0004\% &  0.0005\% &  0.1643\% &  2.5194\% &  12.8080\% \\ \hline 
     relative error for linear case &  0.0727\% &  0.0599\% &  0.3647\% &  2.6736\% &  14.2001\% \\ \hline 
  \end{tabular}
\end{center}
 \caption{Relative errors for different cases.}\label{fig:converge0}
 \end{table}

\subsubsection{Marmousi model}
\begin{figure}
    \centering
    \begin{subfigure}[b]{0.45\textwidth}
        \includegraphics[width=\textwidth]{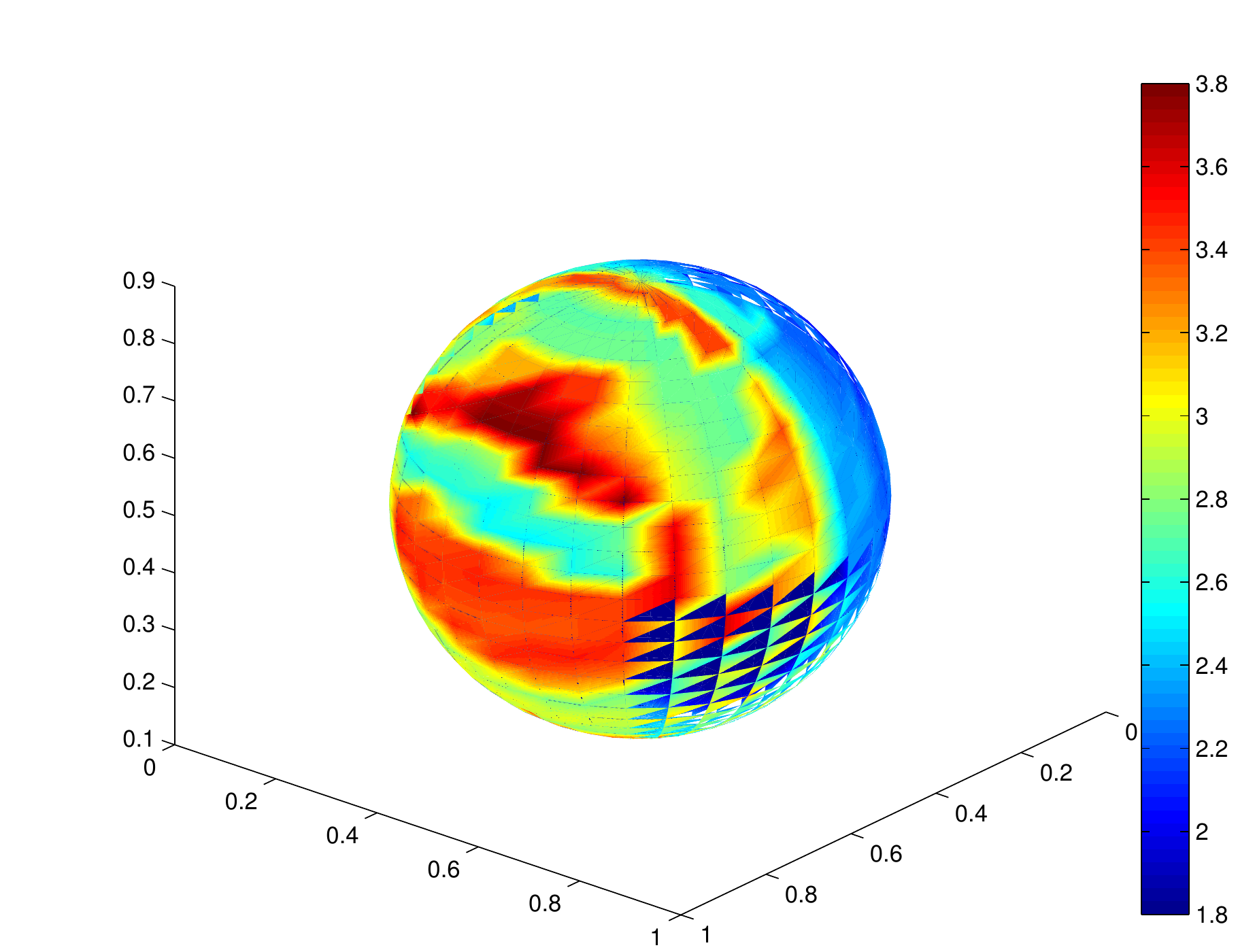}
        \caption{exact solution for first layer (front)}
        \label{fig:true solution}
    \end{subfigure}
    ~ %add desired spacing between images, e. g. ~, \quad, \qquad, \hfill etc. 
        \begin{subfigure}[b]{0.45\textwidth}
        \includegraphics[width=\textwidth]{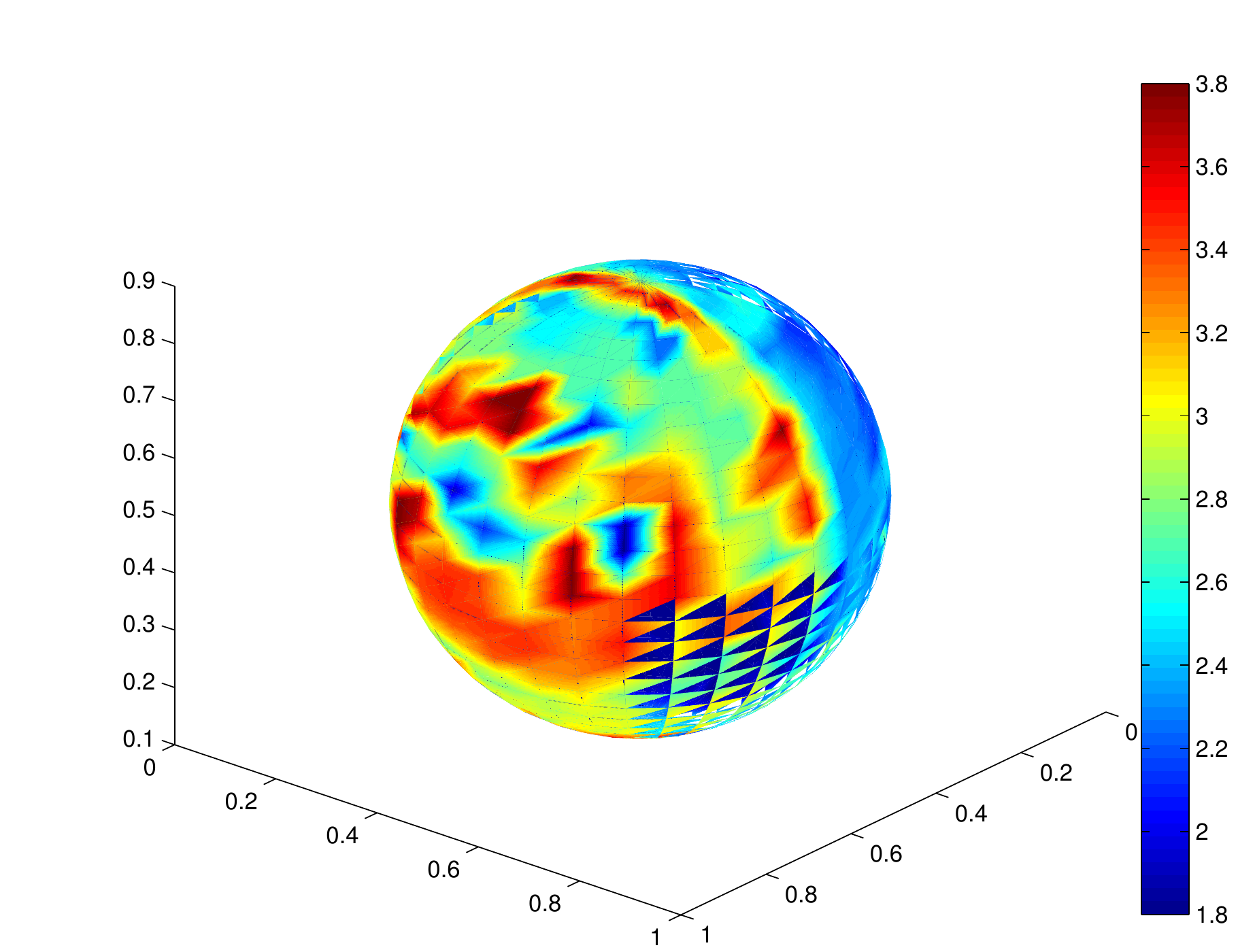}
        \caption{approx. solution for  first layer (front)}
        \label{fig:approximate solution}
    \end{subfigure}

    \begin{subfigure}[b]{0.45\textwidth}
        \includegraphics[width=\textwidth]{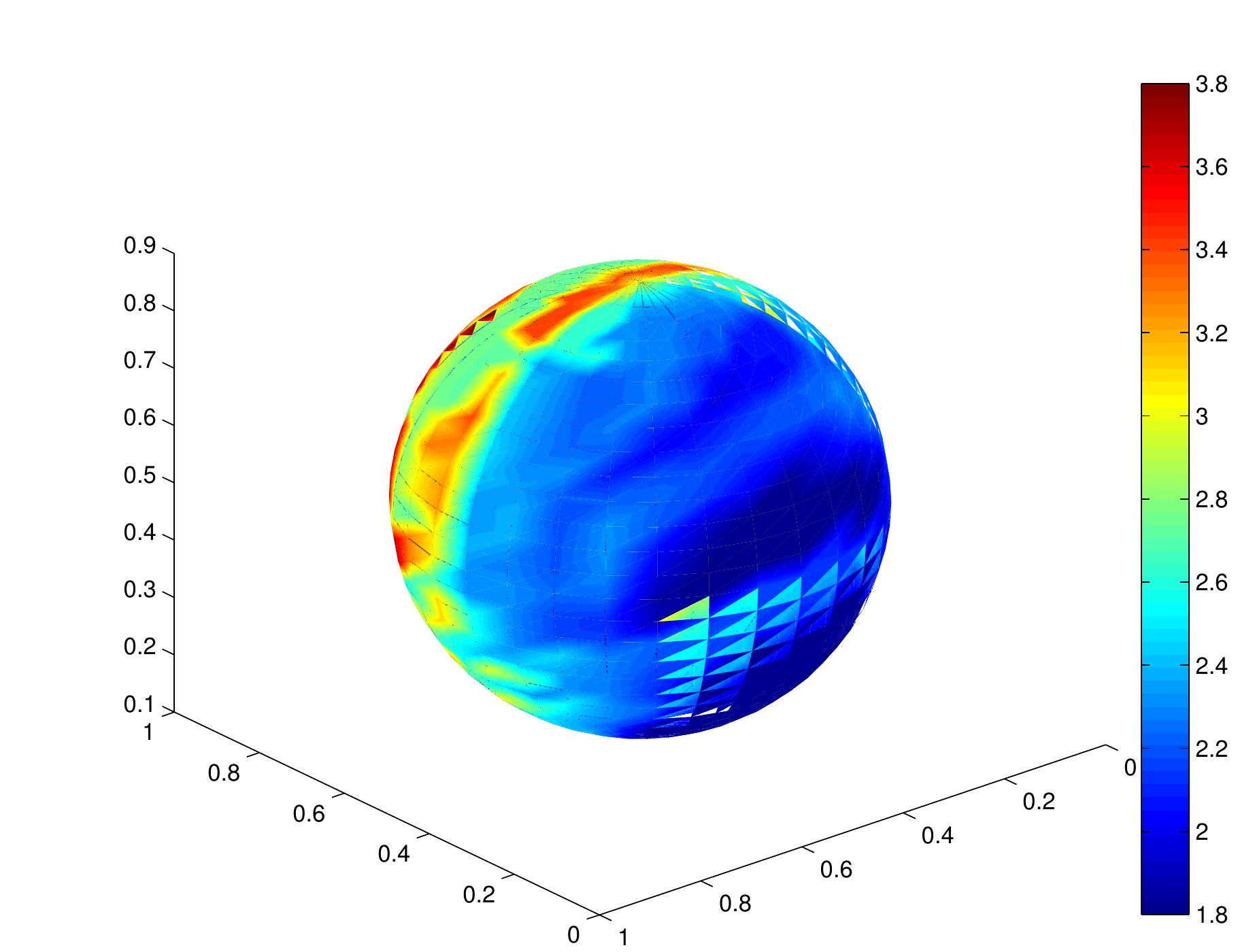}
        \caption{exact solution for first layer (back)}
        \label{fig:true solution}
    \end{subfigure}
    ~ %add desired spacing between images, e. g. ~, \quad, \qquad, \hfill etc. 
        \begin{subfigure}[b]{0.45\textwidth}
        \includegraphics[width=\textwidth]{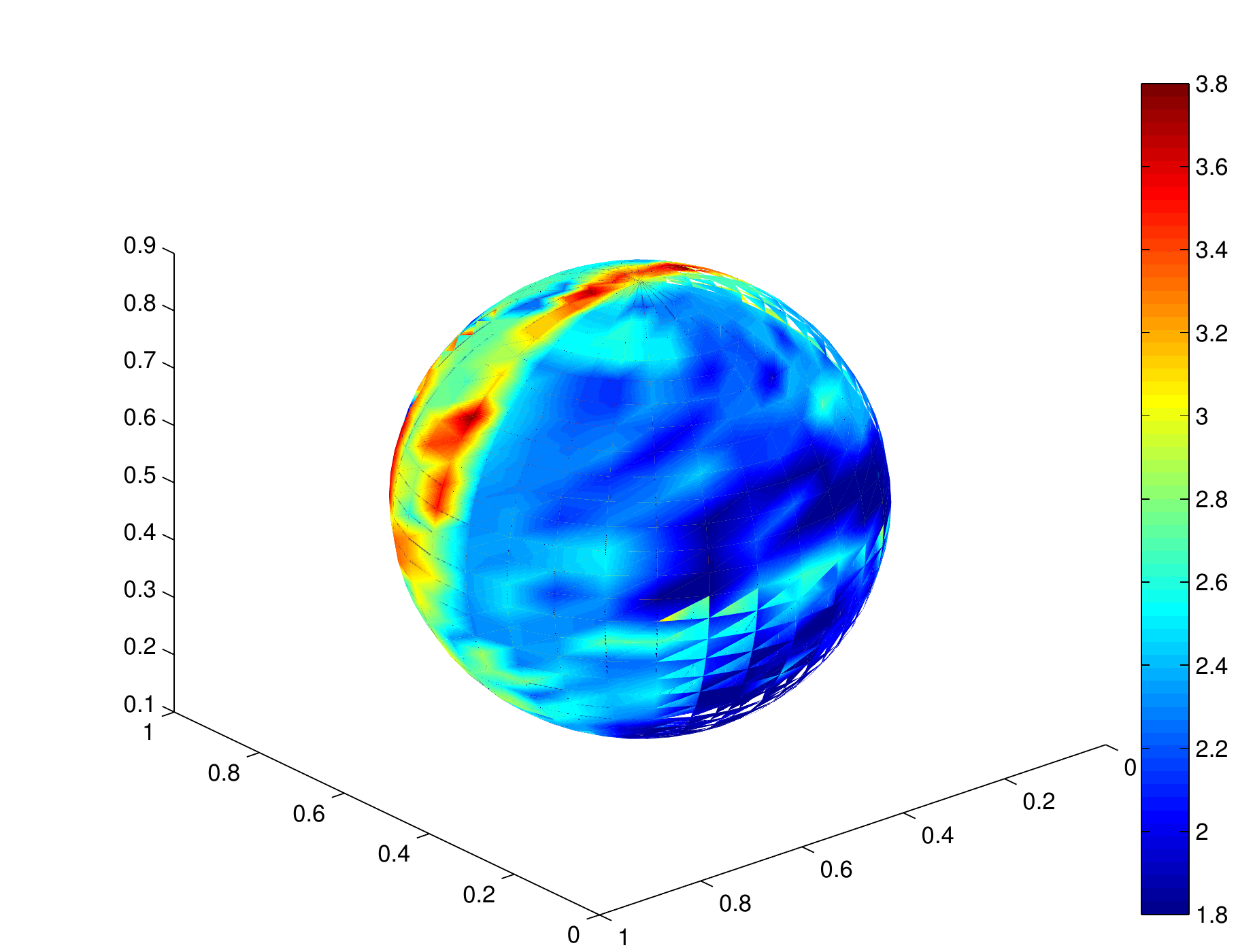}
        \caption{approx. solution for first layer (back)}
        \label{fig:approximate solution}
    \end{subfigure}

         \caption{Graphs of true and approximate solution of first layer of standard marmousi model with fewer layers}
     \label{fig:solution3a}
\end{figure}

\begin{figure}
    \centering

\begin{subfigure}[b]{0.45\textwidth}
        \includegraphics[width=\textwidth]{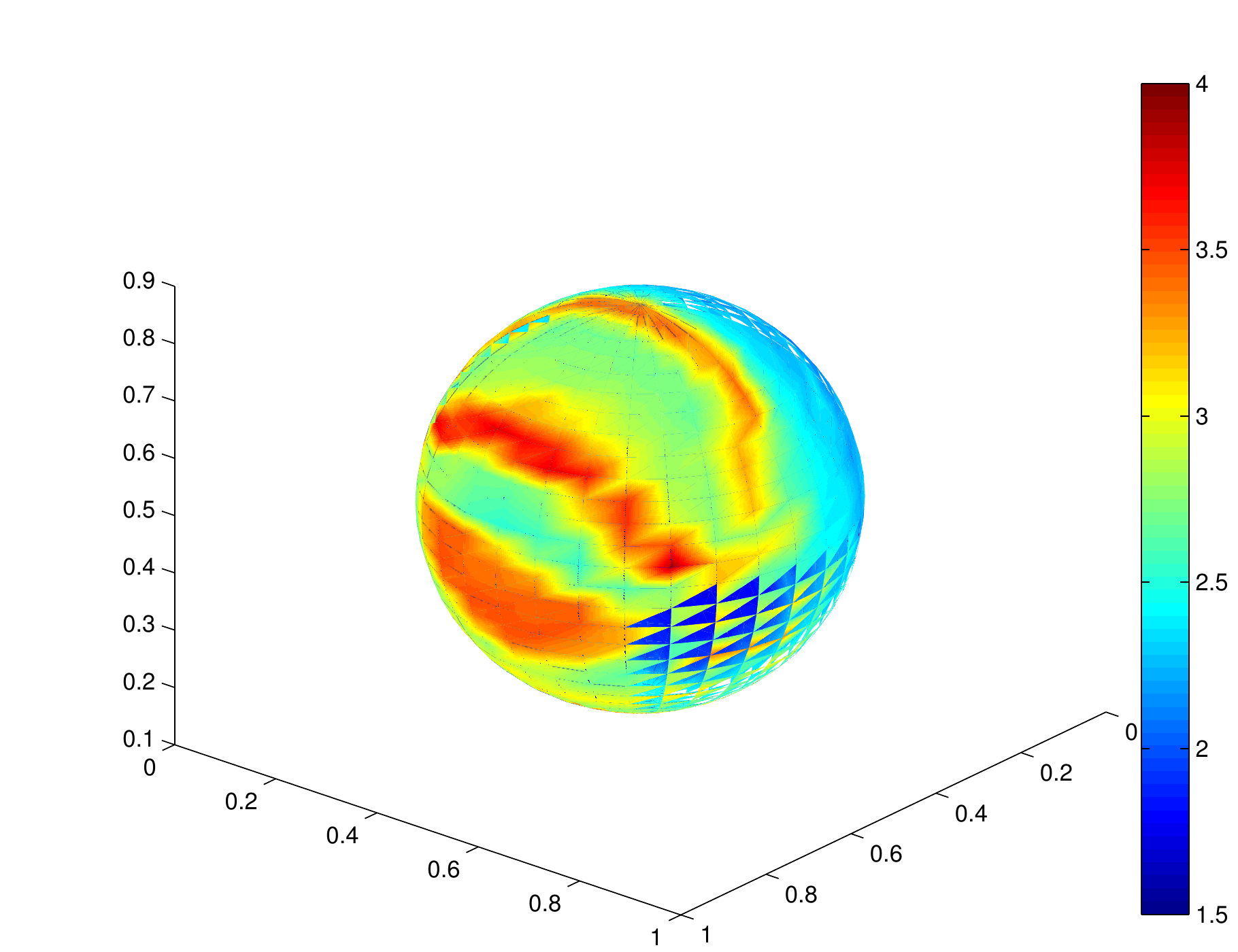}
        \caption{exact solution for second layer (front)}
        \label{fig:true solution}
    \end{subfigure}
    ~ %add desired spacing between images, e. g. ~, \quad, \qquad, \hfill etc. 
        \begin{subfigure}[b]{0.45\textwidth}
        \includegraphics[width=\textwidth]{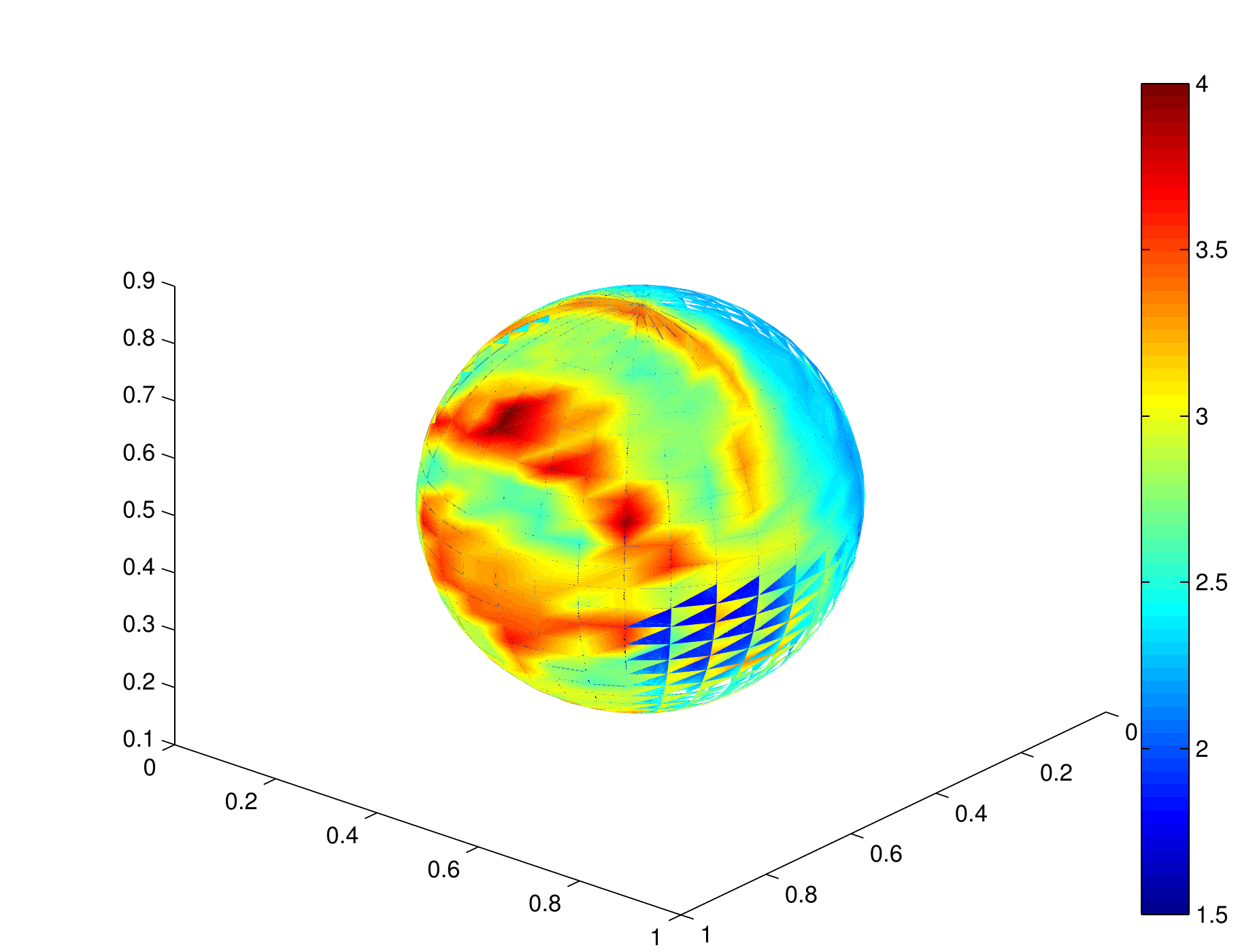}
        \caption{approx. solution for second layer (front)}
        \label{fig:approximate solution}
    \end{subfigure}
      \begin{subfigure}[b]{0.45\textwidth}
        \includegraphics[width=\textwidth]{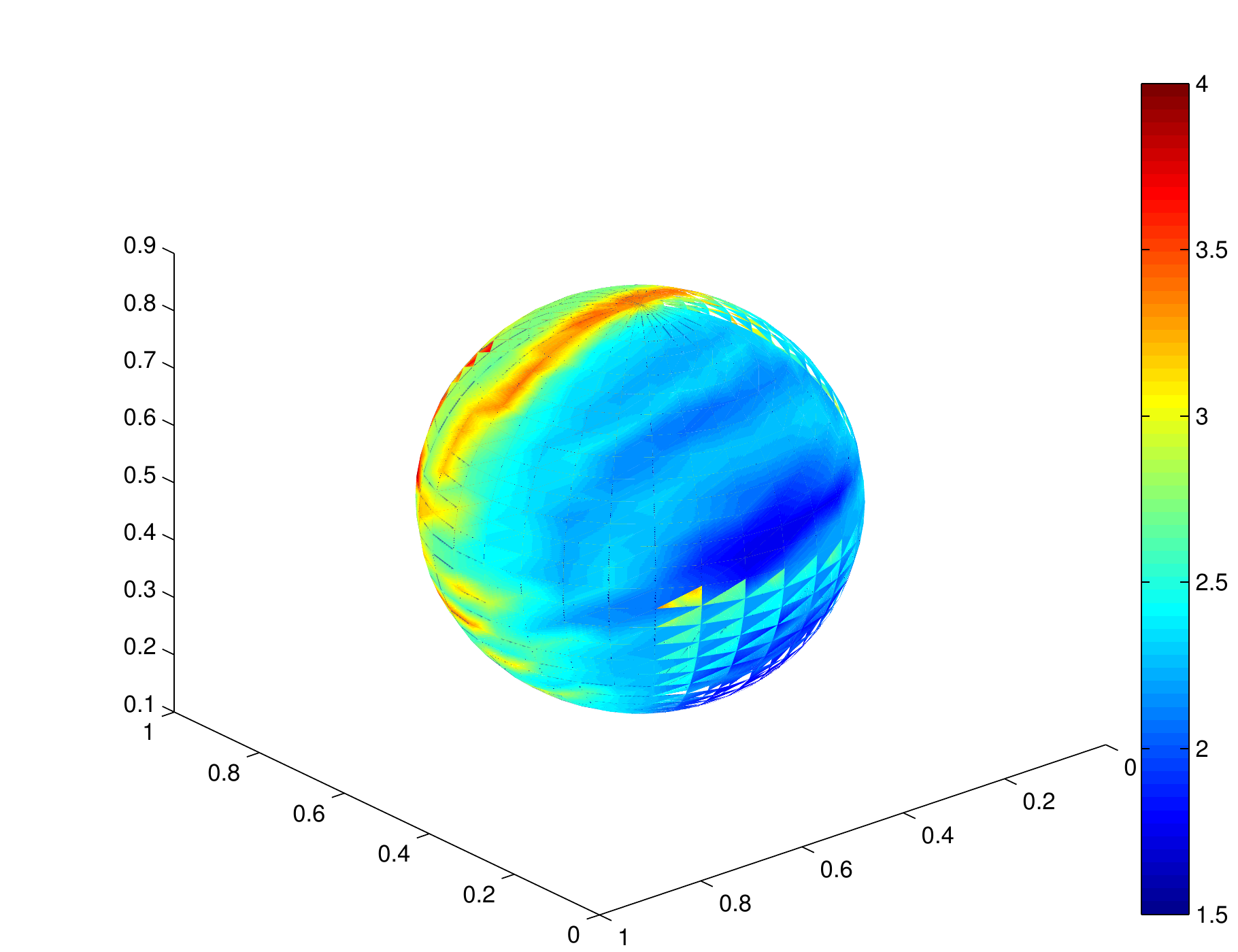}
       \caption{exact solution for second layer (back)}
        \label{fig:true solution}
    \end{subfigure}
    ~ %add desired spacing between images, e. g. ~, \quad, \qquad, \hfill etc. 
        \begin{subfigure}[b]{0.45\textwidth}
        \includegraphics[width=\textwidth]{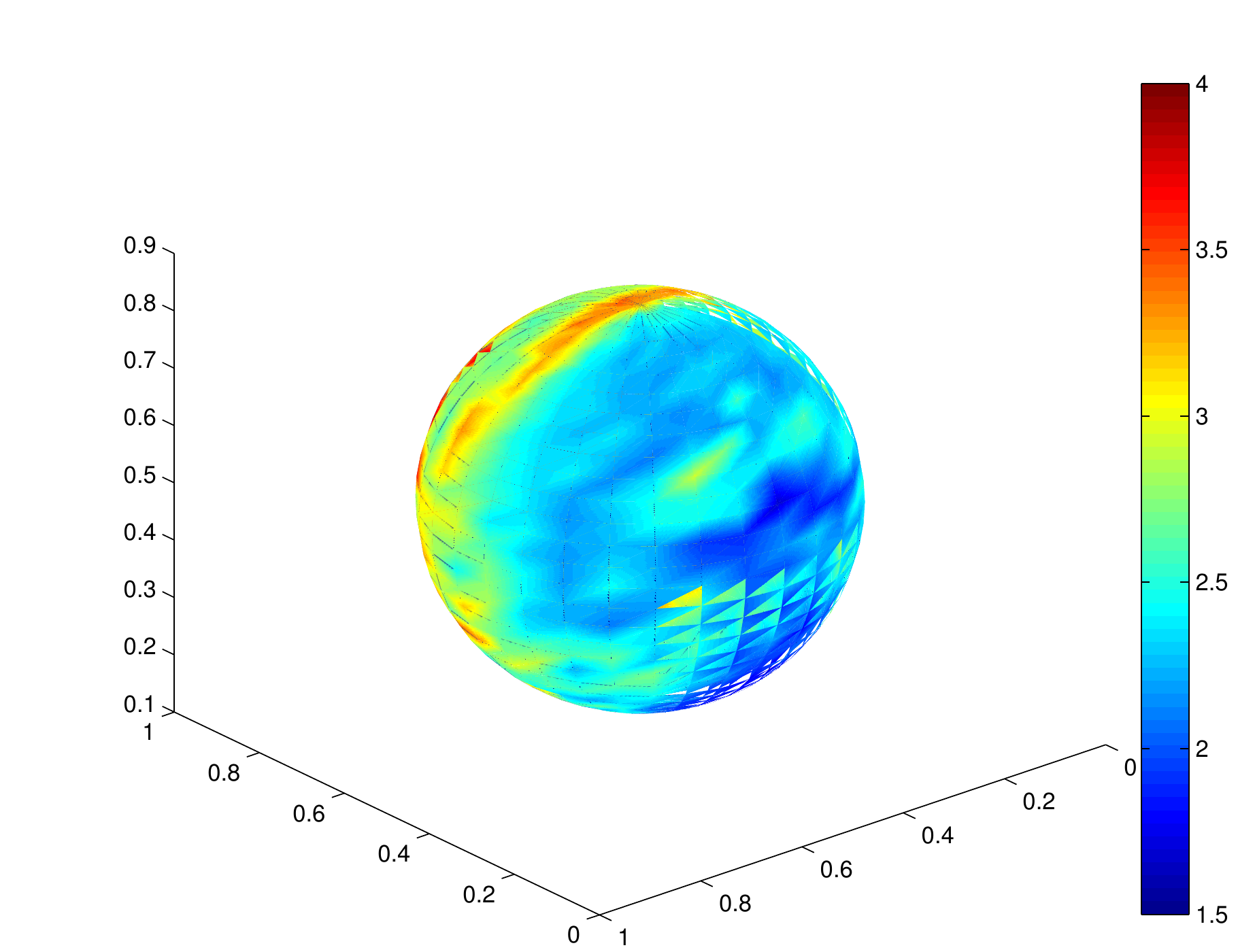}
        \caption{approx. solution for second layer (back)}
        \label{fig:approximate solution}
    \end{subfigure}

\caption{Graphs of true and approximate solution of second layer of standard marmousi model with fewer layers}
     \label{fig:solution3b}
\end{figure}

\begin{figure}
    \centering
\begin{subfigure}[b]{0.45\textwidth}
        \includegraphics[width=\textwidth]{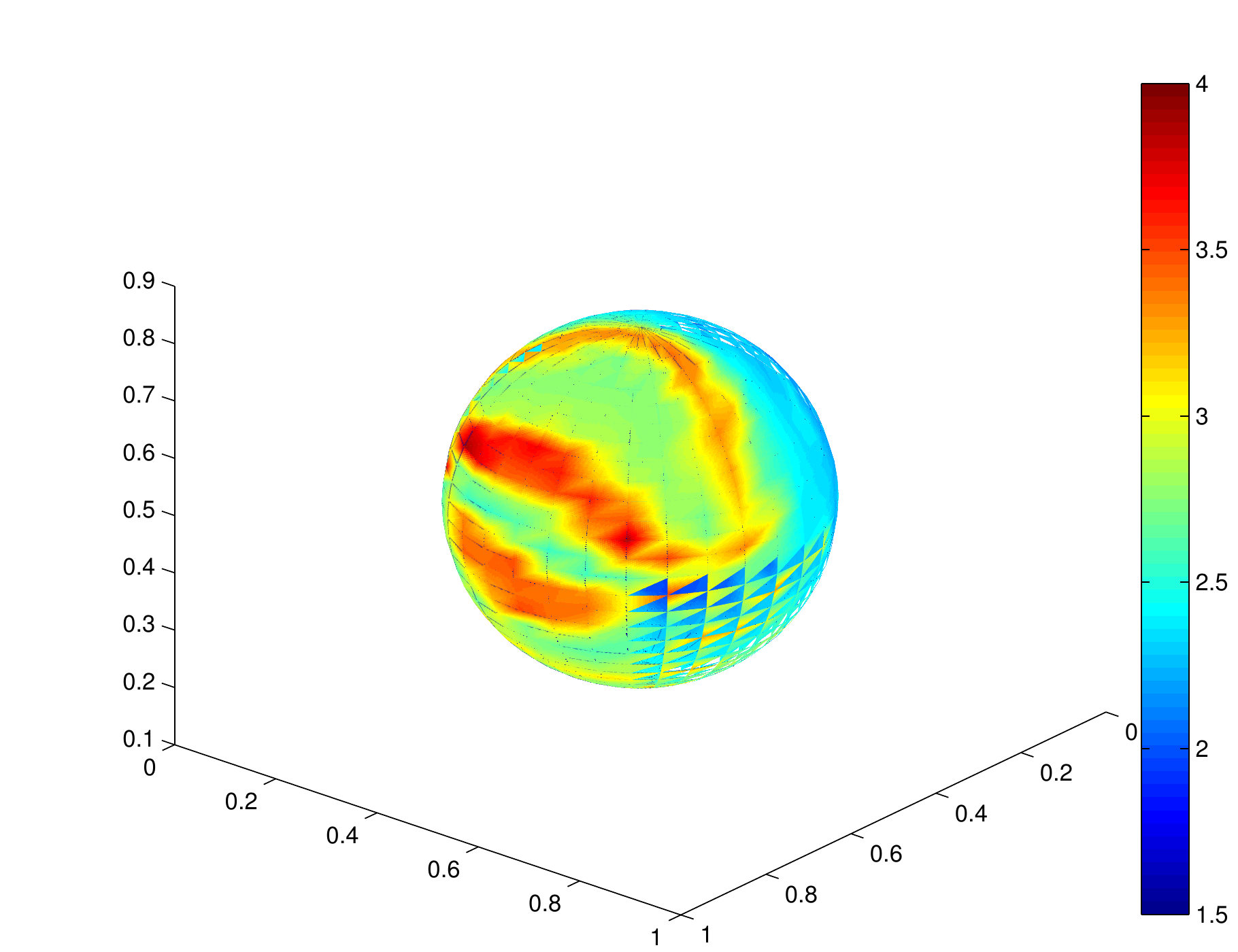}
        \caption{exact solution for third layer (front)}
        \label{fig:true solution}
    \end{subfigure}
    ~ %add desired spacing between images, e. g. ~, \quad, \qquad, \hfill etc. 
        \begin{subfigure}[b]{0.45\textwidth}
        \includegraphics[width=\textwidth]{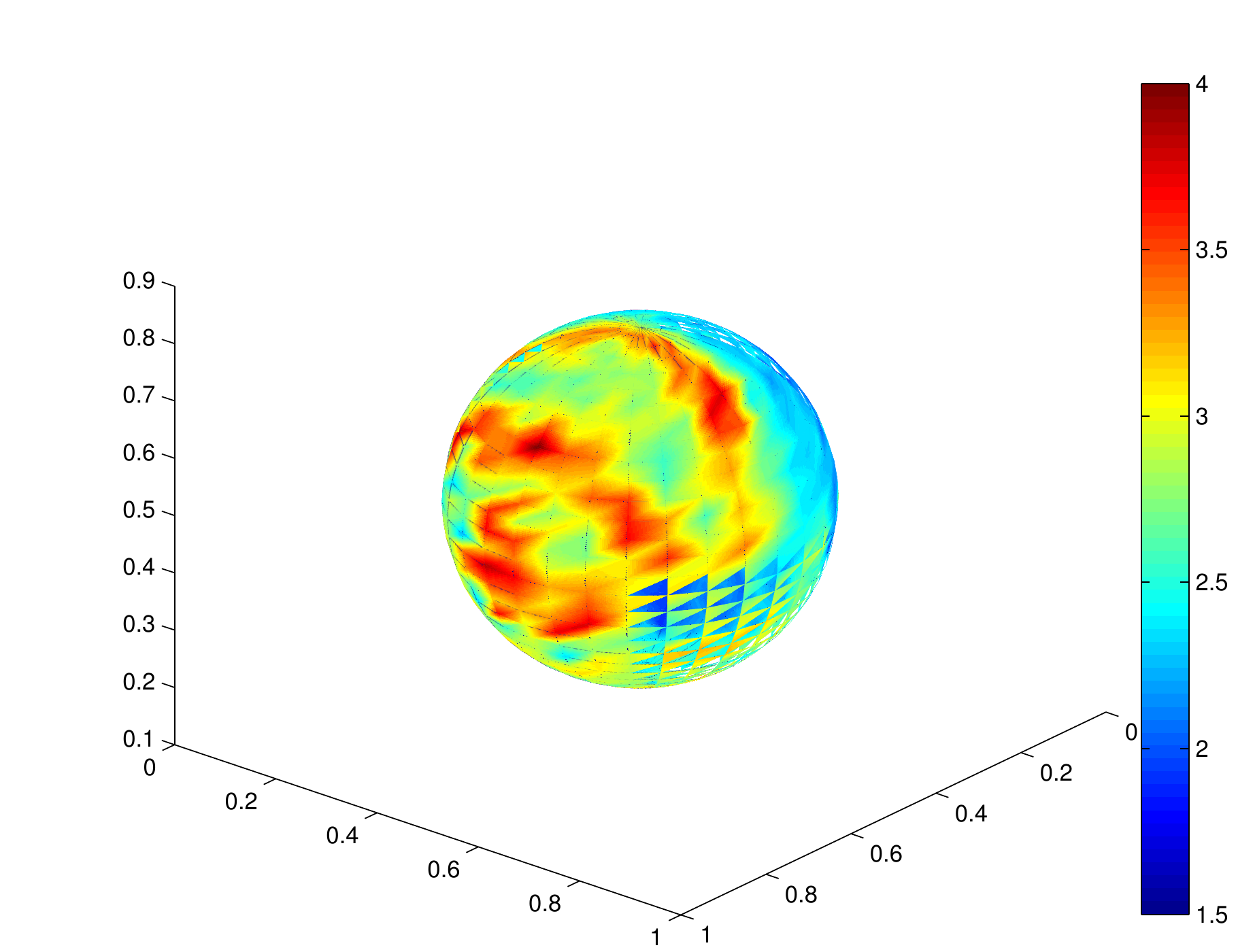}
        \caption{approx. solution for third layer (front)}
        \label{fig:approximate solution}
    \end{subfigure}
      \begin{subfigure}[b]{0.45\textwidth}
        \includegraphics[width=\textwidth]{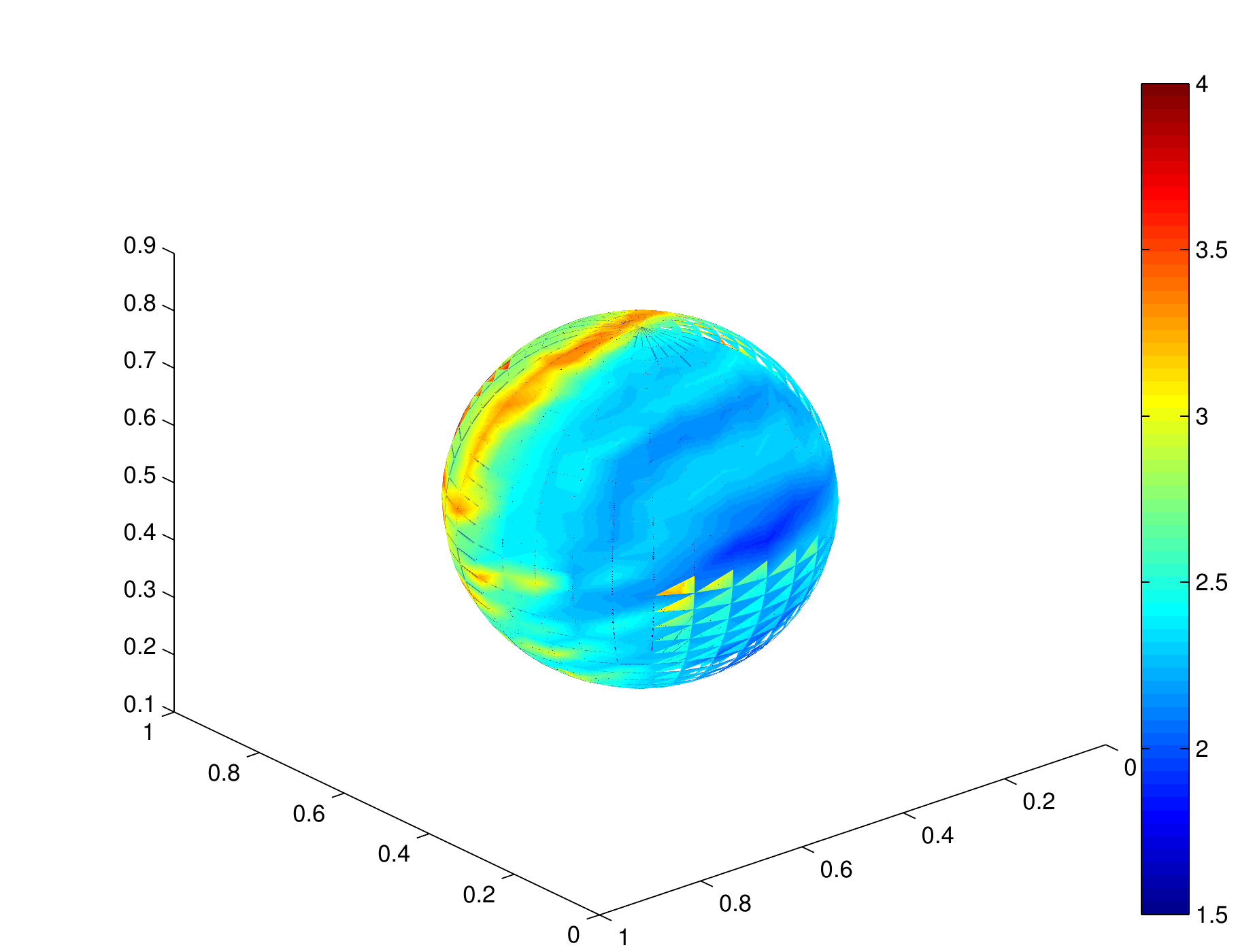}
       \caption{exact solution for third layer (back)}
        \label{fig:true solution}
    \end{subfigure}
    ~ %add desired spacing between images, e. g. ~, \quad, \qquad, \hfill etc. 
        \begin{subfigure}[b]{0.45\textwidth}
        \includegraphics[width=\textwidth]{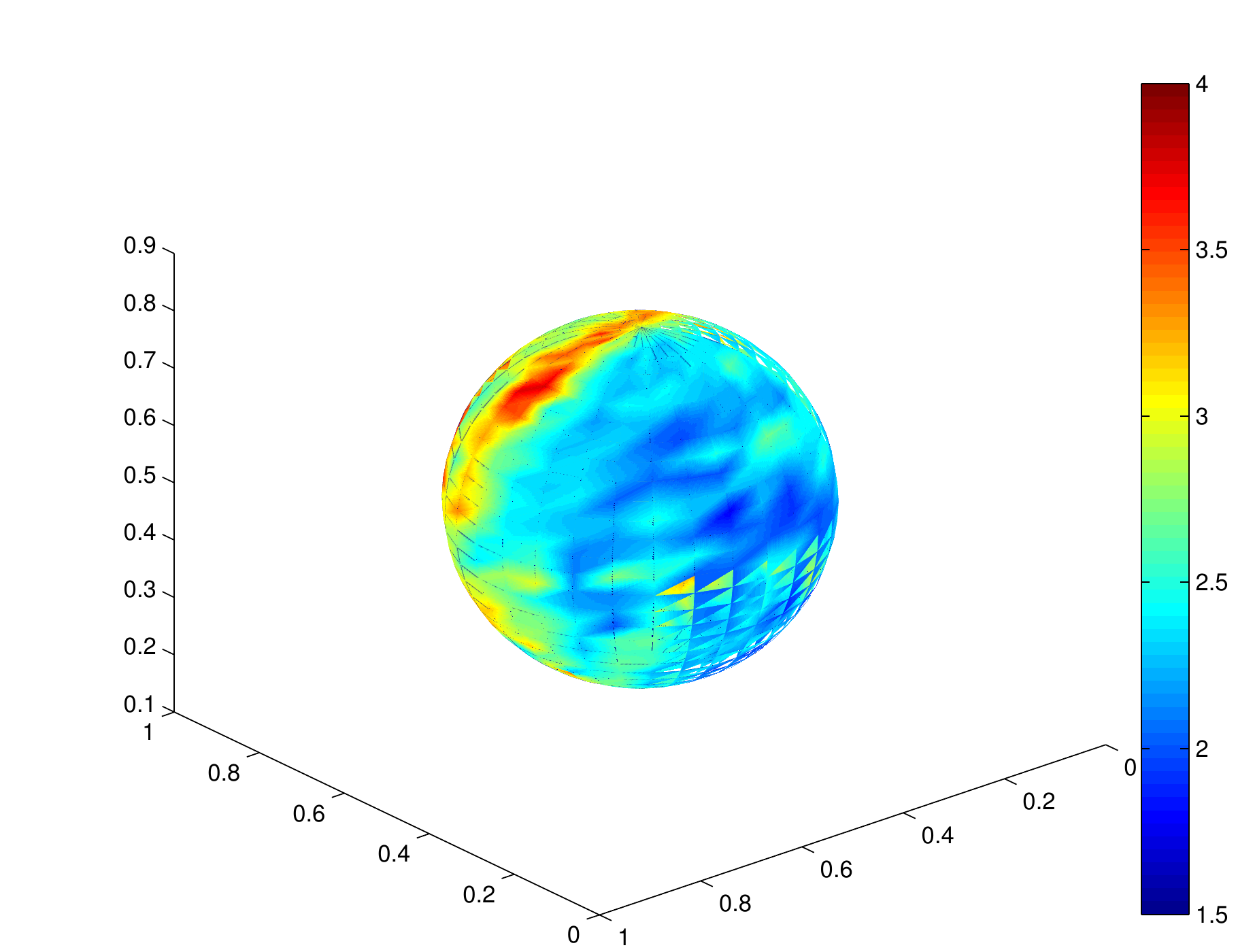}
        \caption{approx. solution for third layer (back)}
        \label{fig:approximate solution}
    \end{subfigure}

\caption{Graphs of true and approximate solution of third layer of standard marmousi model with fewer layers}
     \label{fig:solution3c}
\end{figure}

\begin{figure}
    \centering
\begin{subfigure}[b]{0.45\textwidth}
        \includegraphics[width=\textwidth]{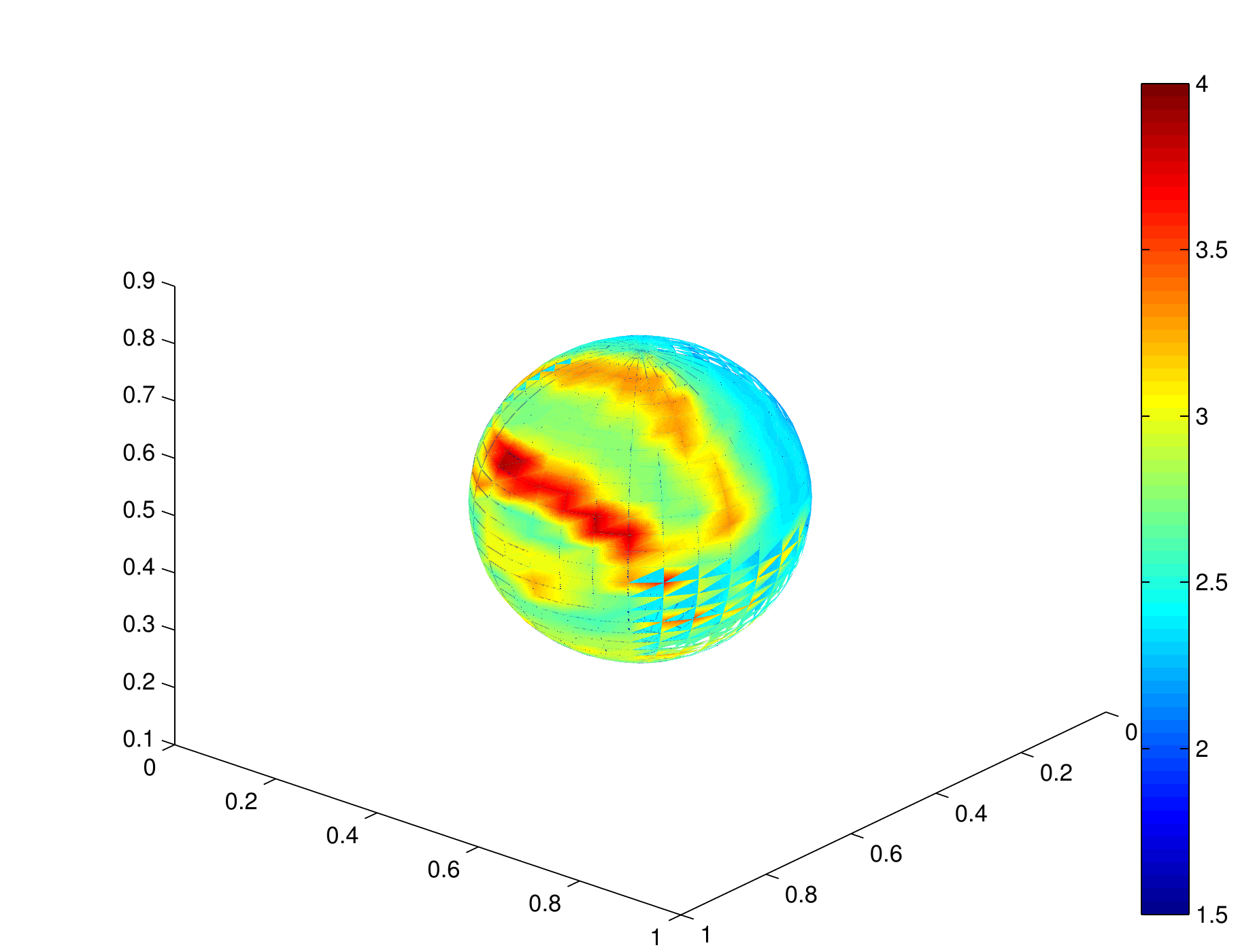}
        \caption{exact solution for forth layer (front)}
        \label{fig:true solution}
    \end{subfigure}
    ~ %add desired spacing between images, e. g. ~, \quad, \qquad, \hfill etc. 
        \begin{subfigure}[b]{0.45\textwidth}
        \includegraphics[width=\textwidth]{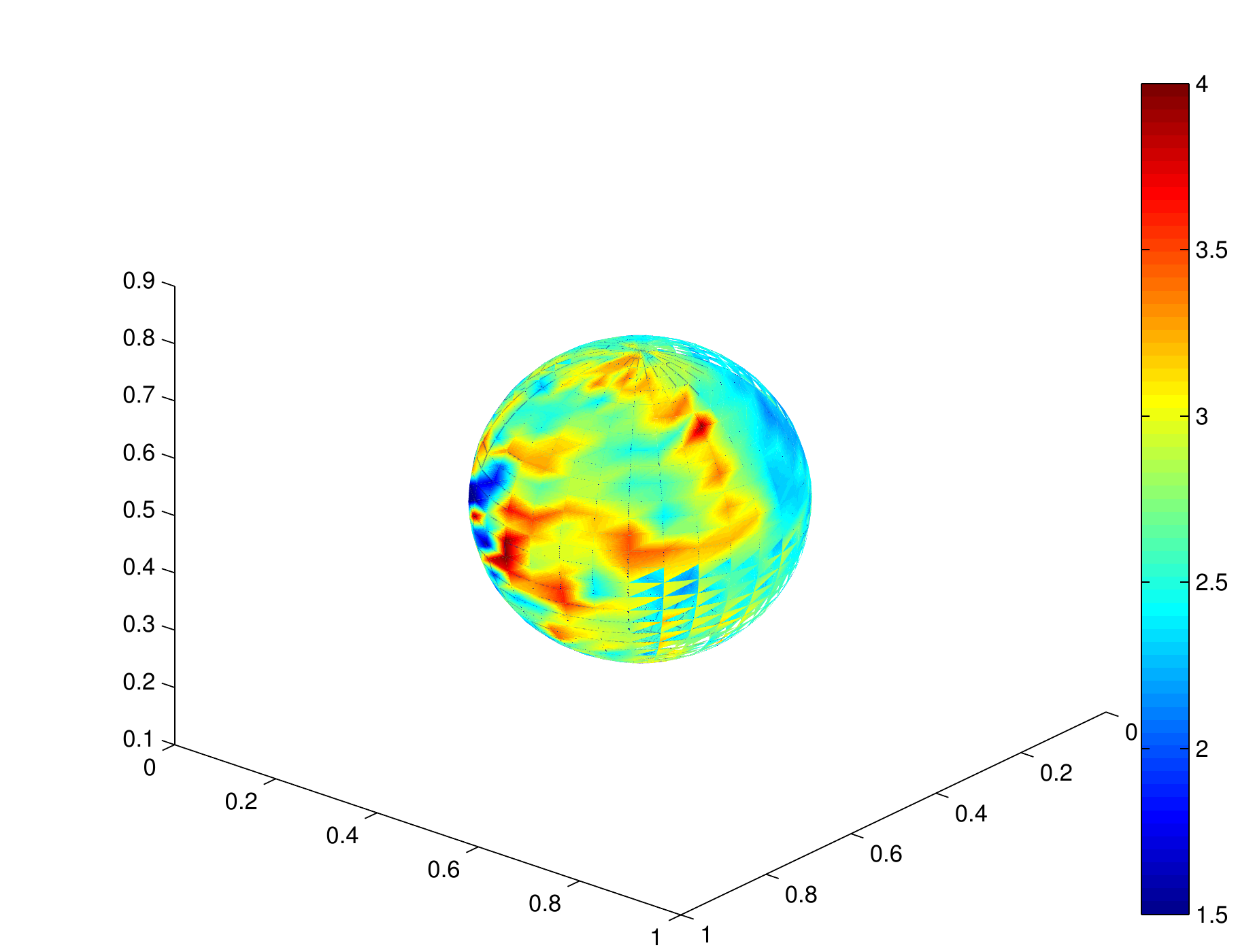}
        \caption{approx. solution for forth layer (front)}
        \label{fig:approximate solution}
    \end{subfigure}
      \begin{subfigure}[b]{0.45\textwidth}
        \includegraphics[width=\textwidth]{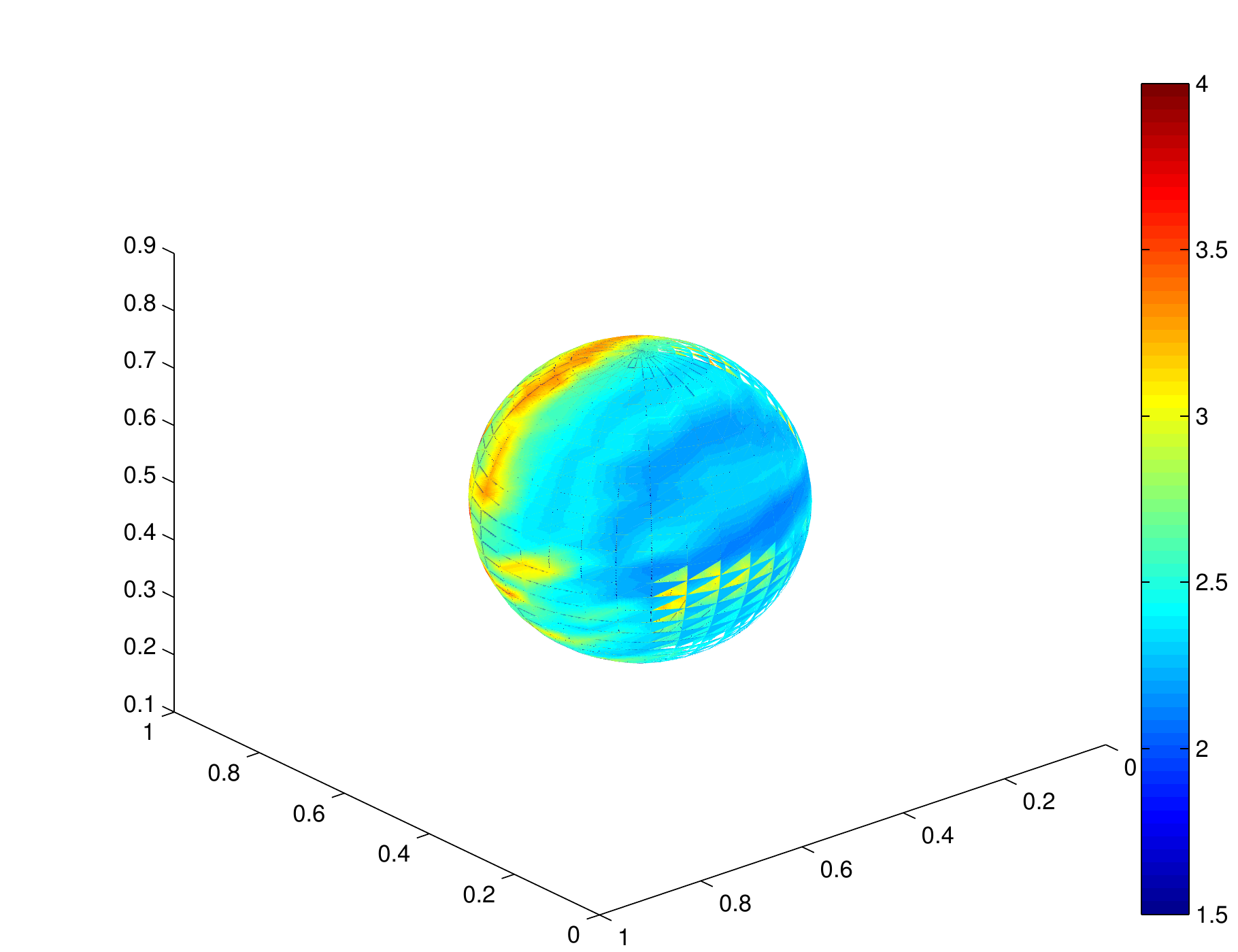}
       \caption{exact solution for forth layer (back)}
        \label{fig:true solution}
    \end{subfigure}
    ~ %add desired spacing between images, e. g. ~, \quad, \qquad, \hfill etc. 
        \begin{subfigure}[b]{0.45\textwidth}
        \includegraphics[width=\textwidth]{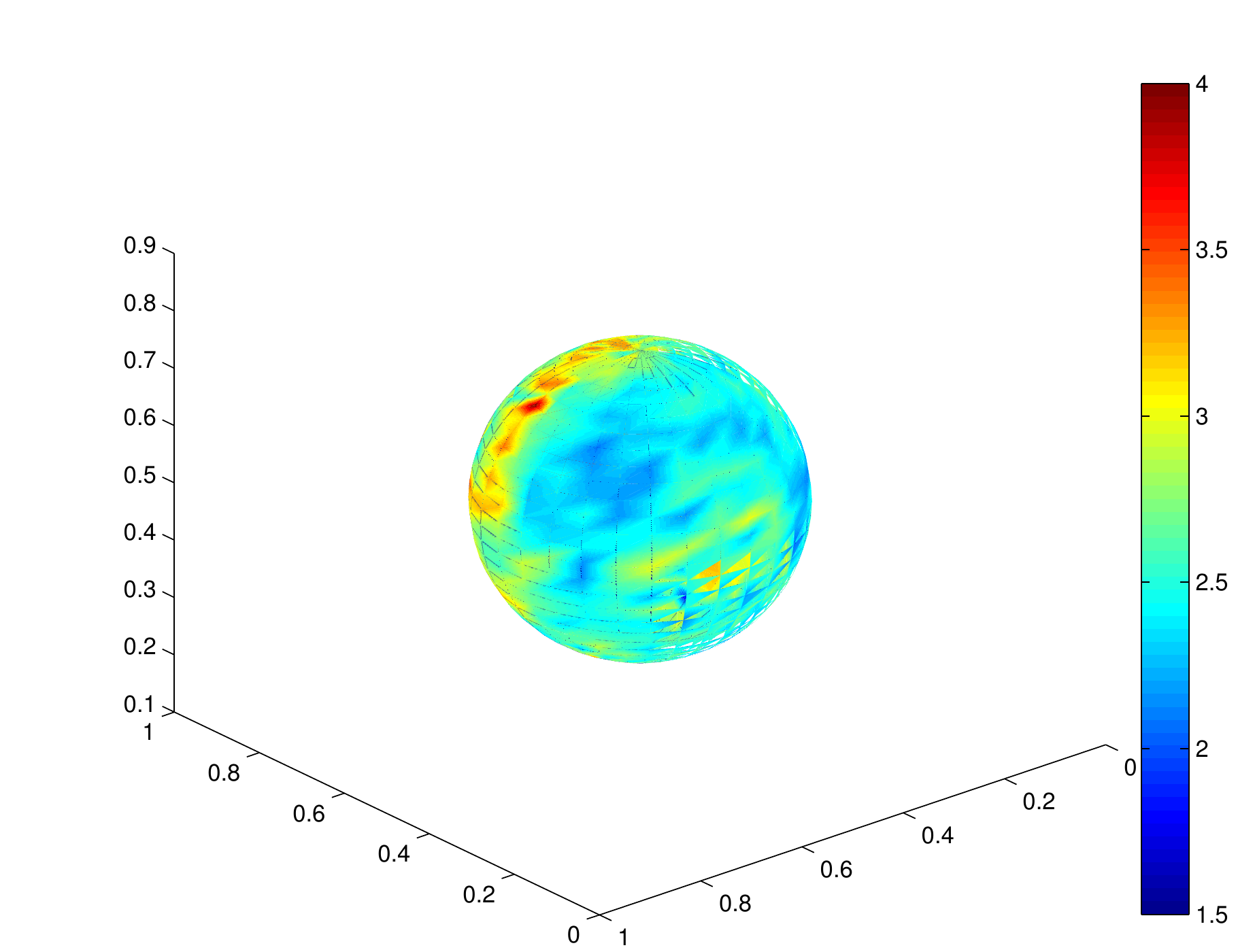}
        \caption{approx. solution for forth layer (back)}
        \label{fig:approximate solution}
    \end{subfigure}

\caption{Graphs of true and approximate solution of forth layer of standard marmousi model with fewer layers}
     \label{fig:solution3d}
\end{figure}

Next, we test the performance using the Marmousi model. We divide the 3D domain into 10 layers and we recover the model in the first few outermost  layers. Then we have the relative error 8.2883\% for first layer , 6.6484\% for the second layer ,9.2633\% for the third layer and 12.8978\% for the forth layer. Figure \ref{fig:solution3a}-\ref{fig:solution3d}, show the graphs of true and approximate solution of first,second , third and forth layers of standard Marmousi model.
We observe that the recovered solutions are in good agreement with the exact solutions. 

 \begin{table}[ht]
    \centering
 \begin{center}
  \begin{tabular}{| l | c | c | c | c | c |}
    \hline
     {}   & $1^{st}$ layer & $2^{nd}$ layer & $3^{rd}$ layer  & $4^{th}$ layer &$5^{th}$ layer \\ \hline
  relative error &  8.2883\% &  6.6484\% &  9.2633\% &  12.8978\% &  13.2901\% \\ \hline 
  \end{tabular}
\end{center}
 \caption{Relative errors for recovering the Marmousi model.}\label{fig:marl}
\end{table}

\section*{Conclusion}

In this paper, we develop a novel numerical scheme for the reconstruction
of an unknown function using its geodesic X-ray transform. 
Our scheme is based on a convergent Neumann series and a layer stripping algorithm. 
In our reconstruction process, we will first compute the unknown function on each small neighborhood
near the boundary of the domain. When the unknown function is reconstructed on a layer near the boundary,
we will reconstruct the unknown function one layer inside the domain, and continue this process until
the function is reconstructed for the whole domain. 
Our method is very efficient, since all computations are performed locally. 
We present a few test cases, including the Marmousi model, to show the performance of the scheme. 
We observe that the quality of approximation will be improved as more terms in the Neumann series are used.
We also observe that the method is quite robust to some noise in the data. 
Finally, we apply this method to a travel time tomography in 3D, in which the inversion of the
geodesic X-ray transform is one important step, and present several numerical results to validate the scheme.

\section*{Acknowledgement}

The research of Eric Chung is partially supported by the Hong Kong RGC General Research Fund
(Project: 14301314) and CUHK Direct Grant for Research 2016-17.
The research of Gunther Uhlmann is partially supported by NSF and Si-Yuan Professorship at IAS, HKUST.

%\clearpage
%\newpage

\appendix
\section {Appendix}
\subsection{Formation of the Hamiltonian system }
The continuous dependence on the initial data of the solution of the Hamiltonian system is characterized by the Jacobian, 
$$J_{g_j}(s)=J_{g_j}(s,X^{(0)}):=\frac{\partial X_{g_j}}{\partial X^{(0)}}(s,X^{(0)})= \left( \begin{array}{cc}
\frac{\partial x}{\partial x^{(0)}} & \frac{\partial x}{\partial \xi^{(0)}}  \\
\frac{\partial \xi}{\partial x^{(0)}} & \frac{\partial \xi}{\partial \xi^{(0)}} \end{array} \right).$$
It can be shown easily from the definition of $J$ and the corresponding Hamiltonian system that $J_{g_j},j=1,2$, satisfies
\begin{align*}
\frac{dJ}{ds}&=\frac{d}{ds}\left(\frac{\partial X_{g_j}}{\partial X^{(0)}}(s,X^{(0)})\right)\\
&=\frac{d}{ds} \left( \begin{array}{cc}
\frac{\partial x}{\partial x^{(0)}} & \frac{\partial x}{\partial \xi^{(0)}}  \\
\frac{\partial \xi}{\partial x^{(0)}} & \frac{\partial \xi}{\partial \xi^{(0)}} \end{array} \right)\\
&= \left( \begin{array}{cc}
\frac{d}{ds}\frac{\partial x}{\partial x^{(0)}} & \frac{d}{ds}\frac{\partial x}{\partial \xi^{(0)}}  \\
\frac{d}{ds}\frac{\partial \xi}{\partial x^{(0)}} & \frac{d}{ds}\frac{\partial \xi}{\partial \xi^{(0)}} \end{array} \right)\\
&= \left( \begin{array}{cc}
\frac{\partial }{\partial x^{(0)}}\frac{dx}{ds} & \frac{\partial }{\partial \xi^{(0)}}\frac{dx}{ds}  \\
\frac{\partial }{\partial x^{(0)}}\frac{d\xi}{ds} & \frac{\partial }{\partial \xi^{(0)}}\frac{d\xi}{ds} \end{array} \right)\\
&= \left( \begin{array}{cc}
\frac{\partial H_{\xi}}{\partial x^{(0)}} & \frac{\partial H_{\xi}}{\partial \xi^{(0)}}  \\
-\frac{\partial H_x}{\partial x^{(0)}} & -\frac{\partial H_x}{\partial \xi^{(0)}} \end{array} \right)\\
&= \left( \begin{array}{cc}
\frac{\partial H_{\xi}}{\partial x}\frac{\partial x}{\partial x^{(0)}} + \frac{\partial H_{\xi}}{\partial \xi}\frac{\partial \xi}{\partial x^{(0)}} & \frac{\partial H_{\xi}}{\partial x}\frac{\partial x}{\partial \xi^{(0)}} + \frac{\partial H_{\xi}}{\partial \xi}\frac{\partial \xi}{\partial \xi^{(0)}}  \\
-(\frac{\partial H_{x}}{\partial x}\frac{\partial x}{\partial x^{(0)}} + \frac{\partial H_{x}}{\partial \xi}\frac{\partial \xi}{\partial x^{(0)}}) & -(\frac{\partial H_{x}}{\partial x}\frac{\partial x}{\partial \xi^{(0)}} + \frac{\partial H_{x}}{\partial \xi}\frac{\partial \xi}{\partial \xi^{(0)}})  \end{array} \right)\\
&=MJ
\end{align*}\\
Then, we have the following equation:
\begin{align*}
\frac{dJ}{ds}=MJ, \qquad J(0)=I,
\end{align*}
where, in terms of $H=H_{g_j}$ ($j=1,2$), the matrix $M$ is defined by
$$M=\left( \begin{array}{cc}
H_{\xi,x} & H_{\xi,\xi}\\
-H_{x,x} & -H_{x,\xi} \end{array}\right).$$

\subsection{Linearizing the Stefanov-Uhlmann identity}
Given boundary measurements for $g_1$, we are interested in recovering the metric $g_1$. We link two metrics by introducing the function
$$F(s):=X_{g_2}\big(t-s,X_{g_1}(s,X^{(0)})\big),$$
where $t=max(t_{g_1},t_{g_2})$ and $t_g=t_g(X^{(0)})$ is the length of the geodesic issued from $X^{(0)}$ with the endpoint on $\Gamma$. Then we first have
\begin{align*}
F(0)&=X_{g_2}\big(t-0,X_{g_1}(0,X^{(0)})\big)\\
&=X_{g_2}(t,X^{(0)})\\
F(t)&=X_{g_2}\big(t-t,X_{g_1}(t,X^{(0)})\big)\\
&=X_{g_2}\big(0,X_{g_1}(t,X^{(0)})\big)\\
&=X_{g_1}(t,X^{(0)})
\end{align*}
Consequently, we have 
\begin{align*}
\int^t_0 F'(s)ds&=F(t)-F(0)\\
&=X_{g_1}(t,X^{(0)})-X_{g_2}(t,X^{(0)}).
\end{align*}
On the other hand, we can denote $V_{g_j}:=\Bigg(\frac{\partial H_{g_j}}{\partial\xi},-\frac{\partial H_{g_j}}{\partial x}\Bigg)$ and then we get
\begin{align*}
F'(s)&=\frac{d}{ds}X_{g_2}\big(t-s,X_{g_1}(s,X^{(0)})\big)\\
&=-\frac{d X_{g_2}}{ds}\big(t-s,X_{g_1}(s,X^{(0)})\big)+\frac{\partial X_{g_2}}{\partial X^{(0)}}\big(t-s,X_{g_1}(s,X^{(0)})\big)\frac{d X_{g_1}}{ds}(s,X^{(0)})\\
&=-V_{g_2}\left(X_{g_2}\big(t-s,X_{g_1}(s,X^{(0)})\big)\right)+\frac{\partial X_{g_2}}{\partial X^{(0)}}\big(t-s,X_{g_1}(s,X^{(0)})\big)V_{g_1} \left(X_{g_1}(s,X^{(0)})\right)
\end{align*}
Then we claim that 
$$V_{g_2}\left(X_{g_2}\big(t-s,X_{g_1}(s,X^{(0)})\big)\right)=\frac{\partial X_{g_2}}{\partial X^{(0)}}\big(t-s,X_{g_1}(s,X^{(0)})\big)V_{g_2} \left(X_{g_1}(s,X^{(0)})\right)$$
To prove the claim, we first need to notice that
\begin{align*}
&\frac{d}{ds}\Big|_{s=0}X(T-s,X(s,X^{(0)}))\\=&\frac{d}{ds}\Big|_{s=0}X(T,X^{(0)})\\
=&0
\end{align*}
Also, we have
\begin{align*}
&\frac{d}{ds}\Big|_{s=0}X(T-s,X(s,X^{(0)}))\\
=&V(X(T,X^{(0)}))+\frac{\partial X}{\partial X^{(0)}}(T,X^{(0)})V(X^{(0)}),~~\forall T
\end{align*}
By setting $T=t-s$ and $X^{(0)}=X_{g_1}(s,X^{(0)})$, thus we prove the claim.\\
Hence, the time integral on the left-hand side is equal to the following \cite{Stefanov5},
\begin{align*}
\int^t_0 F'(s)ds&=\int^t_0\frac{\partial X_{g_2}}{\partial X^{(0)}}\big(t-s,X_{g_1}(s,X^{(0)})\big)\times\big(V_{g_1}-V_{g_2}\big)\big(X_{g_1}(s,X^{(0)})\big)ds\\
&=\int^t_0J_{g_2}\big(t-s,X_{g_1}(s,X^{(0)})\big)\times\big(V_{g_1}-V_{g_2}\big)\big(X_{g_1}(s,X^{(0)})\big)ds,
\end{align*}
where
$$V_{g_j}=\Bigg(\frac{\partial H_{g_j}}{\partial\xi},-\frac{\partial H_{g_j}}{\partial x}\Bigg)=\Bigg(g^{-1}\xi,-\frac{1}{2}\nabla_x(g^{-1}\xi)\cdot \xi\Bigg).$$
This is the so-called Stefanov-Uhlmann identity.

We linearize the above identity about $g_2$ by following a full underlying patin $g_2$,
\begin{equation*}
\int^t_0 F'(s)ds\approx\int^t_0J_{g_2}\big(t-s,X_{g_1}(s,X^{(0)})\big)\times\partial_{g_2}V_{g_2}(g_1-g_2)\big(X_{g_2}(s,X^{(0)})\big)ds,
\end{equation*}
where $\partial_g V_g(\lambda)$ is the derivative of $V_g$ with respect to $g$ at $\lambda$.

Thus, we have the following formula,
\begin{equation}\label{SUid}
X_{g_1}(t,X^{(0)})-X_{g_2}(t,X^{(0)})\approx\int^t_0J_{g_2}\big(t-s,X_{g_2}(s,X^{(0)})\big)\times\partial_{g_2}V_{g_2}(g_1-g_2)\big(X_{g_2}(s,X^{(0)})\big)ds,
\end{equation}
     
% --------------------------------------------------------------
%     You don't have to mess with anything below this line.
% --------------------------------------------------------------
\end{document}